\documentstyle{amsppt}
\magnification=1200
\catcode`\@=11
\redefine\logo@{}
\catcode`\@=13

\define \bn{\Bbb N}
\define \bz{\Bbb Z}
\define \bq{\Bbb Q}
\define \br{\Bbb R}

\define \M{{\Cal M}}
\define\Ha{{\Cal H}}
\define\La{{\Cal L}}
\define \E{{\Cal E}}

\define\rk{\text{rk}~}

%additional to the file

\define\pr{\text{pr}}

\define\discr{\text{discr\ }}

\TagsOnRight
%\NoBlackBoxes

\document

\topmatter
\title
On the classification of hyperbolic root systems of
the rank three. Part II
\endtitle

\author
Viacheslav V. Nikulin \footnote{Supported by
Grant of Russian Fund of Fundamental Research
\hfill\hfill}
\endauthor

\address
Steklov Mathematical Institute,
ul. Gubkina 8, Moscow 117966, GSP-1, Russia
\endaddress
\email
slava\@nikulin.mian.su
\endemail

\abstract
Here we prove classification results announced in Part I
(alg-geom/ 9711032). We classify maximal
hyperbolic root systems of the rank $3$ having
restricted arithmetic type and a generalized lattice Weyl
vector $\rho$ with $\rho^2\ge 0$ (i.e. of elliptic or
parabolic type). We give classification of all reflective of
elliptic or parabolic type elementary hyperbolic lattices of the rank
three.

We apply the same method (narrow places of polyhedra)
which was developed to prove finiteness results on reflective
hyperbolic lattices. We also use some additional arithmetic arguments:
studying of class numbers of central symmetries.

The same methods permit to get similar results for hyperbolic type.
We will consider hyperbolic type in Part III.

These results are important for Theory of Lorentzian Kac--Moody algebras
and some aspects of Mirror Symmetry.
\endabstract

\rightheadtext
{Hyperbolic root systems}
\leftheadtext{V.V. Nikulin}
\endtopmatter

\head
0. Introduction
\endhead

Here we prove results announced in Part I \cite{N14}. We
continue numeration of Sections started in Part I. We also keep
notations of Part I.

We consider main hyperbolic (i. e. of
signature $(1,k)$) lattices (i. e. non-degenerate integral symmetric
bilinear forms) $S$ of the rank three and with square-free
determinant $d$. Here ``main'' means that $S$ should be even for even $d$.
The lattice $S$ is defined uniquely by its determinant $d=\det(S)$ and
some additional invariant $\eta$ (see Sect. 2.2 of Part I or
Sect. 3.2 here) where $0\le \eta <2^t$, here $t$ is the number of odd
prime
divisors of $d$.

Let $W(S)$ be the reflection group of $S$ (generated by reflections
in roots $\alpha \in S$ with $\alpha^2<0$) and $\M$ its fundamental
polyhedron in hyperbolic space $\La (S)$ defined by $S$. We denote by
$A(\M)=\{\phi \in O^+(S)\ |\ \phi(\M)=\M\}$ the group of
symmetries of $\M$. An involution $u\in A(\M)$ is called a
{\it central symmetry} if $u$ acts as a central symmetry in $\La (S)$.
Two central symmetries from $A(\M)$ are called {\it equivalent} if
they are conjugate in $A(\M)$. In Sect. 3, Theorem 3.2.1 we give a
formula for the number $h=hnr(S)$ of classes of central symmetries of the
lattice $S$ with invariants $(d,\eta)$. This formula uses the Legendre
symbol and class numbers of imaginary  quadratic fields.
We use this formula to give the list of all the lattices
$S$ with $h=hnr(S)\le 1$ and $d\le 100000$. This list was announced in
Table 3 of Part I and contains 206 invariants $(d,\eta)$.
This list is important for classification of the reflective lattices $S$.

A hyperbolic lattice $S$ is called {\it reflective} if there exist
a non-zero $\rho \in S$ and a subgroup $A\subset A(\M)$
of finite index such that $A(\rho)=\rho$. Here $\rho$ is called
a {\it generalized lattice Weyl vector.}
If $S$ has a generalized lattice Weyl vector $\rho$ with $\rho^2>0$,
then $S$ is called {\it reflective of elliptic type}
(or {\it elliptically reflective}). If $S$ has a generalized
lattice Weyl vector $\rho$ with $\rho^2=0$ and does not have a generalized
lattice Weyl vector $\rho$ with $\rho^2>0$, then $S$ is called reflective
of {\it parabolic type}. If $S$ has a generalized
lattice Weyl vector $\rho$ with $\rho^2<0$ and does not have a generalized
lattice Weyl vector $\rho$ with $\rho^2\ge 0$, then $S$ is called
reflective
of {\it hyperbolic type}. It is easy to prove that $h\le 1$ if $S$ is
reflective of elliptic type, $h=0$ if $S$ is reflective of parabolic type
and $h=0,2$ if $S$ is reflective of hyperbolic type and $\rk S=3$.
Thus, Table 3 contains all elliptically or parabolically reflective main
hyperbolic lattices $S$ of the rank $3$ and
with square-free determinant $d\le 100000$.

In \cite{N4}, \cite{N5}, \cite{N11} and \cite{N13} we proved
finiteness results about reflective hyperbolic lattices using
some geometrical arguments: studying of narrow places of fundamental
polyhedra. In Sect. 4 we improve and optimize this method 
for 2-dimensional case 
to apply it for classification of reflective hyperbolic lattices of
rank $3$. In Sect. 4.3 we apply these results to fundamental
polygons of the reflection groups of reflective hyperbolic lattices
of elliptic or parabolic type and of rank 3. As a result, we get
estimates of the determinant and some other invariants of these lattices.

In Sect. 5 we apply results of Sects. 4 and 5 to classify all
elliptically or parabolically reflective hyperbolic lattices of
the rank 3 and with square-free determinant.
The list of these lattices was announced in Tables 1 and 2 of Part I
(see \cite{N14}). In particular, we prove that all elliptically or
parabolically reflective main hyperbolic lattices of the rank
$3$ and with square-free determinant are contained in Table 3.
To study reflective type of lattices, we use Vinberg's algorithm
\cite{V2}.

The same methods as developed here permit to classify
reflective hyperbolic lattices of the rank 3 of hyperbolic type.
We hope to do this in Part III.

\head
3. The number of classes of central symmetries of main hyperbolic lattices
with square-free determinant and of the rank $3$.
\endhead

\subhead
3.1. Reminding of some classical results about binary positive
lattices
\endsubhead

Here we remind some classical results about binary positive quadratic
forms (e. g. see \cite{B-Sh} and \cite{C}).

We consider binary (i.e. of the rank two) positive definite lattices $K$
with square-free determinant $d=\det K$. They will be called
{\it fundamental}. If $d$ is odd, the fundamental
lattice $K$ is unimodular over $\bz_2$. Considering $K$ over
$\bz_2$, it is easy to see that $K$ is odd if
$d\equiv 1,\,2\mod 4$.
If $d\equiv -1 \mod 4$, the lattice $K$ may be odd or even.
The number
$$
D=\cases
-\det K, &\text{if $K$ is even,}\\
-4\det K, &\text{if $K$ is odd}
\endcases
\tag{3.1.1}
$$
is called {\it discriminant} of a binary lattice $K$. Discriminants of
fundamental binary lattices are called
{\it fundamental discriminants}. Thus, for a fundamental discriminant
$D$ either $D\equiv 1\mod 4$ or $D\equiv \pm 4,\,8\mod 16$.
A binary positive lattice $K$ is called {\it classical fundamental}
if $d=\det K$ is square-free and additionally $K$ is even if
$d\equiv -1\mod 4$. The discriminant $D$ of a classical fundamental
binary lattice $K$ is called the {\it classical fundamental
discriminant}. Thus, for a classical fundamental
discriminant $D$ either
$D\equiv 1 \mod 4$ (when $K$ is even) or $D\equiv -4,\,8\mod 16$.

It is known

\proclaim{Theorem 3.1.1}
The number $h(D)$
of proper classes (i. e. classes of preserving orientation
isomorphisms of oriented lattices) of classical fundamental binary
lattices $K$ of the discriminant $D<0$ is equal to the class number of
the imaginary quadratic field $\bq (\sqrt{D})$
of the discriminant $D$ and is given by the Dirichlet's formula
$$
h(D)={w\over 2|D|}\sum_{0<r<|D|}\left({D\over r}\right)r
\tag{3.1.2}
$$
where
$$
w=\cases
6, &\text{if $D=-3$,}\\
4, &\text{if $D=-4$,}\\
2, &\text{otherwise,}
\endcases
\tag{3.1.3}
$$
and  $\left({\,D\,\over r}\right)$ is the Kronecker
(i.e. generalized Legendre)  symbol.
\endproclaim

Let $K$ be an odd binary
lattice of an odd square-free determinant $d$
where $d\equiv -1\mod 4$ (i. e. it is fundamental, but it
is not classical fundamental). The maximal even sublattice of $K$ is
the lattice $T(2)$ where $T$ is the fundamental (even) lattice of
the determinant $d$. Thus, $T(2)\subset K$ is an odd overlattice of
index $2$ of the lattice $T(2)$. If $d\equiv -1\mod 8$, then the
overlattice $K$ is unique. If $d\equiv -5\mod 8$, then
there are $3$ overlattices $T(2)\subset K$ (which we should
consider up to proper automorphisms of $T$ to get the number of
proper classes of $K$).
We have $\discr K=-4d$ and denote $h(-4d)$
the number of proper classes of lattices $K$ of that
discriminant. By considerations
above, we get
$$
h(4D)=\cases
h(D),  &\text{if $D\equiv 1\mod 8$,}\\
3h(D), &\text{if $D\equiv 5\mod 8$ and $D<-3$,}\\
1=h(D),&\text{if $D=-3$.}
\endcases
\tag{3.1.4}
$$

Let $K$ be a fundamental positive binary lattice of the
discriminant $D$. Here $D\equiv 0,\,1\mod 4$. If
$D\equiv 0 \mod 4$, then $D\equiv \pm 4,\, 8 \mod 16$.
We use the discriminant form technique \cite{N1}.
The genus of $K$
is defined by the discriminant form
of $K$ (it is quadratic if $K$ is even,
and is bilinear if $K$ is odd) which is defined by the map
$\mu:p\mapsto \mu_p\in \{0,\,1\}$,
$p$ runs through all odd prime $p\vert D$.
For odd (prime) $p\vert D$, the discriminant
form $b_{K_p}=b^{(p)}_{\theta_p}(p)$ where
$\left({\theta_p\over p}\right)=(-1)^{\mu_p}$.
The form $b_{K_2}$
is trivial if either $D\equiv 1\mod 4$ or $D\equiv 4\mod 8$.
If $D\equiv 8 \mod 16$, then $b_{K_2}=b_1^{(2)}(2)$.

We denote by $t$ {\it the number of prime divisors of $D$}.
If $D\equiv 1\mod 4$ (equivalently, $K$ is even),
existence of $K$ with the invariants
$(D,\,\mu)$ is equivalent to
$$
\sum_{p\vert d}{[(1-p)+4\mu_p]} \equiv 2 \mod 8,
\tag{3.1.5}
$$
and there are $2^{t-1}$ genuses.

If $D\equiv -4\mod 16$, (equivalently, $d\equiv 1\mod 4$)
existence of $K$ with the invariants
$(D,\mu)$ is equivalent to
$d\equiv 1 \mod 4$ which is given.
Thus, number of genuses is
equal to $2^{t-1}$.

If $D\equiv 4\mod 16$ (equivalently,
$d\equiv -1\mod 4$ and the lattice $K$ is odd),
existence of $K$ is
equivalent to existence of the corresponding even lattice of the
discriminant $D=-d$ (by the construction above),
and number of genuses is equal to $2^{t-2}$.

If $D\equiv 8\mod 16$, there are no conditions on
$(D, \mu)$ for
existence of $K$, and number of genuses is equal to
$2^{t-1}$.

Thus, {\it the number of different genuses of the discriminant $D$
is equal to}
$2^{\tau (D)}$ where
$$
\tau(D)=
\cases
t-1, &\text{if $D\equiv 1\mod 4$,}\\
t-1, &\text{if $D\equiv -4,\, 8\mod 16$,}\\
t-2, &\text{if $D\equiv 4  \mod 16$,}
\endcases
\tag{3.1.6}
$$
and $t$ is the number of all different prime divisors of $D$.
It is known that each genus contains the same number of
classes. Thus, number of classes $h(D)_g$
in a genus of a fundamental lattice $K$ is
equal to
$$
h(D)_g={h(D)\over 2^{\tau (D)}}.
\tag{3.1.7}
$$

A fundamental binary lattice $K$ is called {\it ambiguous} if
$O(K)$ contains a reflection (equivalently, the oriented lattice $K$ is
proper equivalent to the lattice $K$ with the opposite orientation).
If $O(K)$ contains a
reflection in $\delta_1 \in K$ where $\delta_1$ is primitive,
then it contains reflection in a primitive
$\delta_2\in K$ which is orthogonal to $\delta_1$.
We can suppose that $n_1=\delta_1^2\le n_2=\delta_2^2$.
It follows that a reflective lattice $K$
is generated by either $\{\delta_1,\,\delta_2\}$ or
$\{\delta_1,\delta_2,(\delta_1+\delta_2)/2\}$. We say that
the type is I if we have the first possibility, and
the type is II if we have the second one. It is easy to
see (using classification of $2$-dimensional reflection
groups) that $n_1\le n_2$ and type are invariants of
the ambiguous class $K$ except $K=(1,1;I)=(2,2;II)$ of the
discriminant $-4$. Except this case, two ambiguous binary
lattices are isomorphic iff they have the same invariants
$(n_1,n_2;\text{type})$. We denote by $K(n_1,n_2;\text{type})$ the
ambiguous positive lattice with invariants $(n_1,n_2;\text{type})$.
A lattice $K(n_1,n_2;I)$ exists for any
natural $n_1,\, n_2\in \bn$. A lattice
$K(n_1,n_2;II)$ exists iff $n_1,n_2\in \bn$,
$n_1\equiv n_2\equiv 0\mod 2$ and $n_1+n_2\equiv 0\mod 4$.

We have $\det K(n_1,n_2;I)=n_1n_2$, and
the ambiguous binary lattice
$K(n_1,n_2;I)$ is fundamental if and only if
both $n_1$ and $n_2$ are square-free and
$(n_1,n_2)=1$. This lattice is odd of the discriminant
$D=-4n_1n_2$.

We have $\det K(n_1,n_2;II)=n_1n_2/4$, and
$K(n_1,n_2;II)$ is fundamental iff
both $n_1$, $n_2$ are square-free and $(n_1,n_2)=2$.
The lattice $K(n_1,n_2;II)$ is even iff $n_1+n_2\equiv 0\mod 8$.
Thus, we have
$$
D(K(n_1,n_2;II))=
\cases
-n_1n_2/4, &\text{if $n_1+n_2\equiv 0\mod 8$}\\
-n_1n_2,   &\text{otherwise}.
\endcases
$$

It follows that the number $hr(D)$
of ambiguous classes (proper or improper, it does not matter
because they are improper equivalent to itself)
of the discriminant $D$ is equal to
$$
hr(D)=2^{\tau(D)}=\cases
2^{t-1}, \text{if $D\equiv 1\mod 4$,}\\
2^{t-1}, \text{if $D\equiv -4\mod 16$,}\\
2^{t-2}, \text{if $D\equiv 4\mod 16$,}\\
2^{t-1}, \text{if $D\equiv 8\mod 16$.}
\endcases
\tag{3.1.8}
$$
This number is equal to the number $2^{\tau(D)}$ of genuses
because ambiguous lattices correspond to elements of order
$2$ in the group $A$ of classes of discriminant $D$. We
denote this group as $A_2$. For a class $a\in A$, the
class $-a\in A$ is the class which is
improper equivalent to $a$. The
group of genuses is $A/2A$, and natural homomorphism
$A\to A/2A$ is the genus map.

Further we will be especially interesting in genuses which
do not contain more than one class of non-ambiguous
lattices with respect to the general (i.e.
proper or improper) equivalence.
A genus is called {\it ambiguous} if it
contains an ambiguous lattice. Otherwise it is called
{\it non-ambiguous}.

Using the genus homomorphism, it is easy to see that
a discriminant $D$ contains a
non-ambiguous genus iff
$2^{\tau (D)+1}\vert h(D)$. If $2^{\tau (D)+1}\vert h(D)$,
the number of classes of general equivalence (i.e. proper or
improper) in a non-ambiguous genus is equal to
$h(D)\over 2^{\tau (D)+1}$.
Thus, we have

\proclaim{Lemma 3.1.2} There exists a non-ambiguous genus
of the discriminant $D$ iff
\newline
$2^{\tau (D)+1}\vert h(D)$.
The number of classes of the general equivalence (i.e. proper or
improper)
in a non-ambiguous genus is
equal to $h(D)/2^{\tau (D)+1}$. In particular,
a non-ambiguous genus contains exactly
one class of general equivalence
iff $h(D)=2^{\tau(D)+1}$.
\endproclaim

Number of ambiguous classes in an ambiguous
genus is equal to $\sharp 2A\cap A_2$. By elementary
considerations with finite Abelian groups, we get

\proclaim{Lemma 3.1.3} An ambiguous genus contains
only ambiguous classes iff
$h(D)=2^m$ where $m  \le 2\tau(D)$, and
the class group $A\cong (\bz/4 )^{m-\tau(D)}\oplus (\bz/2)^{2\tau (D)-m}$.
Equivalently, $h(D)=2^m$ where $m\le 2\tau (D)$,
and the number of ambiguous classes in the principal genus
$2A$ is equal to $2^{m-\tau(D)}$.
\endproclaim

\proclaim{Lemma 3.1.4} An ambiguous genus contains
exactly one non-ambiguous class iff
either $h(D)=2^{\tau(D)}\cdot 3$ or
$h(D)=2^{\tau (D)+2}$ and
$A\cong (\bz/8)\oplus (\bz/2)^{\tau(D)-1}$. We
remark that if $h(D)=2^{\tau (D)+2}$,
then either $D$ is of the
type of Lemma 3.1.3 (when an ambiguous
genus contains ambiguous classes only) or
$A\cong (\bz/8)\oplus (\bz/2)^{\tau(D)-1}$ and
an ambiguous genus contains only one class of
general equivalence of non-ambiguous lattices.
The last case is characterized by the property that
the principal genus $2A$ contains
exactly two ambiguous classes.
\endproclaim

Below we calculate numbers
$hnr(D,\mu)$, $hr_I(D,\mu)$,
$hr_{II}(D,\mu)$ and $hr(D,\mu)$
of non-ambiguous classes of general equivalence,
ambiguous classes of type $I$, ambiguous classes
of type $II$ and ambiguous classes respectively
of the genus $(D,\mu)$. We have
$$
hr(-4,\mu)=hr_I(-4,\mu)=hr_{II}(-4,\mu)=1\ \text{and\ } hnr(-4,\mu)=0.
\tag{3.1.9}
$$
If $D\not=-4$, we have
$$
hr(D,\mu)=hr_I(D,\mu)+hr_{II}(D,\mu)
\tag{3.1.10}
$$
and
$$
hnr(D,\mu)=(h(D)/2^{\tau(D)}-hr_I(D,\mu)-hr_{II}(D,\mu ))/2.
\tag{3.1.11}
$$
Below we calculate $hr_I(D,\mu)$ and $hr_{II}(D,\mu)$. Since the
principal genus $pr$ is ambiguous, for the number of
non-ambiguous classes in an ambiguous genus $ambig$ we get
$$
hnr(D, ambig)=(h(D)/2^{\tau(D)}-hr_I(D,pr)-hr_{II}(D,pr))/2
\tag{3.1.12}
$$
We remind (Lemma 3.1.2) that for a non-ambiguous genus
$nambig$ we have
$$
hnr(D,nambig)=
\cases
0,  &\text{if $2^{\tau(D)+1}\not\vert h(D)$,}\\
h(D)\over 2^{\tau(D)+1}, &\text{if $2^{\tau(D)+1}\vert h(D)$.}
\endcases
\tag{3.1.13}
$$

Case $D\equiv 1\mod 4$. Then the determinant
$d=-D\equiv -1\mod 4$. All ambiguous
classes of the genus $(D,\mu)$ are given by $K(2d_1,2d_2;II)$ such
that $d_1d_2=d$ and
$\left({2d_1/p\over p}\right)=(-1)^{\mu_p}$ if $p\vert d_1$ and
$\left({2d_2/p\over p}\right)=(-1)^{\mu_p}$ if $p\vert d_2$. Thus,
$$
\split
hr_{II}(D,\mu)=\sharp \{d_1\vert d \ \mid \
d_1\le d/d_1\ &\&\ \left({2d_1/p\over p }\right)=(-1)^{\mu_p}
\ \forall p\vert d_1\\
&\& \left({2(d/d_1)/p\over p}\right)=(-1)^{\mu_p}\
\forall p\vert (d/d_1)\},
\endsplit
\tag{3.1.14}
$$
and
$$
hr_{I}(D,\mu)=0.
\tag{3.1.15}
$$
The principal genus $pr$ is given by the ambiguous lattice
$K(2,2d;II)$. All ambiguous classes of that genus
are $K(2d_1,2d_2;II)$ such that $d_1d_2=d$ and
$\left({2d_1/p\over p}\right)=\left({2d_1d_2/p\over p}\right)$
if $p\vert d_1$ and $\left({2d_2/p\over p}\right)=
\left({2d_1d_2/p\over p}\right)$ if $p\vert d_2$. Thus,
$$
hr_{II}(D,pr)=\sharp \{d_1\vert d \ \mid \
d_1\le d/d_1\ \&\ \left({d/d_1\over p }\right)=1\ \forall p\vert d_1\ \&
\left({d_1\over p}\right)=1\ \forall p\vert (d/d_1)\}
\tag{3.1.16}
$$
and
$$
hr_{I}(D,pr)=0.
\tag{3.1.17}
$$

Case $D\equiv 4\mod 16$. Then $d=-D/4\equiv -1\mod 4$. All ambiguous
classes of the genus $(D,\mu)$ are given by $K(d_1,d_2;I)$ such
that $d_1d_2=d$ and
$\left({d_1/p\over p}\right)=(-1)^{\mu_p}$ if $p\vert d_1$ and
$\left({d_2/p\over p}\right)=(-1)^{\mu_p}$ if $p\vert d_2$. Thus,
$$
\split
hr_I (D,\mu)=\sharp \{d_1\vert d \ \mid \
d_1\le d/d_1 &\&\ \left({d_1/p\over p }\right)=(-1)^{\mu_p}\
\forall p\vert d_1\\
&\&\left({(d/d_1)/p\over p}\right)=(-1)^{\mu_p}\
\forall p\vert (d/d_1)\},
\endsplit
\tag{3.1.18}
$$
and
$$
hr_{II}(D,\mu)=0.
\tag{3.1.19}
$$
The principal genus is given by the ambiguous lattice
$K(1,d;I)$. All ambiguous classes of that genus
are $K(d_1,d_2;I)$ such that $d_1d_2=d$ and
$\left({d_1/p\over p}\right)=\left({d_1d_2/p\over p}\right)$
if $p\vert d_1$ and
$\left({d_2/p\over p}\right)=\left({d_1d_2/p\over p}\right)$
if $p\vert d_2$. Thus,
$$
hr_I(D,pr)=\sharp \{d_1\vert d\ \mid  \
d_1\le d/d_1\ \&\ \left({d/d_1\over p }\right)=1\ \forall p\vert d_1\ \&
\left({d_1\over p}\right)=1\ \forall p\vert (d/d_1)\},
\tag{3.1.20}
$$
and
$$
hr_{II}(D,pr)=0.
\tag{3.1.21}
$$

Case $D\equiv -4\mod 16$. Then $d=-D/4\equiv 1\mod 4$.
All ambiguous classes of type I of the genus $(D,\mu)$
are given by $K(d_1,d_2;I)$ such
that $d_1d_2=d$ and
$\left({d_1/p\over p}\right)=(-1)^{\mu_p}$ if $p\vert d_1$ and
$\left({d_2/p\over p}\right)=(-1)^{\mu_p}$ if $p\vert d_2$. Thus,
$$
\split
hr_I(D,\mu)=\sharp \{d_1\vert d \ \mid \
d_1\le d/d_1\ &\&\ \left({d_1/p\over p }\right)=(-1)^{\mu_p}\
\forall p\vert d_1\\
&\&\left({(d/d_1)/p\over p}\right)=(-1)^{\mu_p}\ \forall
p\vert (d/d_1)\}.
\endsplit
\tag{3.1.22}
$$
All ambiguous classes of the type II of $(D,\mu)$ are given
by $K(2d_1,2d_2;II)$ such
that $d_1d_2=d$ and
$\left({2d_1/p\over p}\right)=(-1)^{\mu_p}$ if $p\vert d_1$ and
$\left({2d_2/p\over p}\right)=(-1)^{\mu_p}$ if $p\vert d_2$. Thus,
$$
\split
hr_{II}(D,\mu)=\sharp \{d_1\vert d \ \mid \
d_1\le d/d_1\ &\&\ \left({2d_1/p\over p }\right)=(-1)^{\mu_p}\
\forall p\vert d_1\\
&\&\left({2(d/d_1)/p\over p}\right)=(-1)^{\mu_p}\
\forall p\vert (d/d_1)\}.
\endsplit
\tag{3.1.23}
$$
The principal genus is given by the ambiguous lattice
$K(1,d;I)$. All ambiguous classes of the type I of that
genus are $K(d_1,d_2;1)$ such that $d_1d_2=d$ and
$\left({d_1/p\over p}\right)=
\left({d_1d_2/p\over p}\right)$ if $p\vert d_1$ and
$\left({d_2/p\over p}\right)=
\left({d_1d_2/p\over p}\right)$ if $p\vert d_2$. Thus
$$
hr_I(D,pr)= \sharp \{d_1\vert d \ \mid \
d_1\le d/d_1\ \&\ \left({d/d_1\over p }\right)=1\ \forall p\vert d_1
\&\left({d_1\over p}\right)=1\ \forall p\vert (d/d_1)\}.
\tag{3.1.24}
$$
All ambiguous classes of the type II of the principal genus
are $K(2d_1,2d_2;II)$ such that $d_1d_2=d$ and
$\left({2d_1/p\over p}\right)=\left({d_1d_2/p\over p}\right)$
if $p\vert d_1$ and
$\left({2d_2/p\over p}\right)=\left({d_1d_2/p\over p}\right)$
if $p\vert d_2$. Thus
$$
\split
&hr_{II}(D,pr)=\\
&\sharp \{d_1\vert d \ \mid \
d_1\le d/d_1\ \&\ \left({d/d_1\over p }\right)=
\left({2\over p}\right)\ \forall p\vert d_1\ \&
\left({d_1\over p}\right)=\left({2\over p}\right)\
\forall p\vert (d/d_1)\}.
\endsplit
\tag{3.1.25}
$$

Case $D\equiv 8\mod 16$. Then $d=-D/4\equiv 2\mod 4$.
All ambiguous classes of the genus $(D,\mu)$
are given by $K(d_1,2d_2;I)$ such
that $d_1d_2=d/2$ and
$\left({d_1/p\over p}\right)=(-1)^{\mu_p}$ if $p\vert d_1$ and
$\left({2d_2/p\over p}\right)=(-1)^{\mu_p}$ if $p\vert d_2$. Thus,
$$
\split
hr_I (D,\mu)=\sharp \{d_1\vert d/2 \ \mid \
&\&\ \left({d_1/p\over p }\right)=(-1)^{\mu_p}\
\forall p\vert d_1\\
&\&\left({(d/d_1)/p\over p}\right)=(-1)^{\mu_p}\
\forall p\vert (d/2d_1)\},
\endsplit
\tag{3.1.26}
$$
and
$$
hr_{II}(D,\mu)=0.
\tag{3.1.27}
$$

The principal genus is given by the ambiguous lattice
$K(1,d;1)$. All ambiguous classes of that genus
are $K(d_1,2d_2;I)$ such that $d_1d_2=d/2$ and
$\left({d_1/p\over p}\right) =\left({2d_1d_2/p\over p}\right)$
if $p\vert d_1$ and
$\left({2d_2/p\over p}\right)=\left({2d_1d_2/p\over p}\right)$
if $p\vert d_2$. Thus,
$$
hr_I(D,pr)=\sharp \{d_1\vert (d/2)\ \mid  \
\left({d/d_1\over p}\right)=1\ \forall p\vert d_1\ \&
\left({d_1\over p}\right)=1\ \forall p\vert (d/2d_1)\}
\tag{3.1.28}
$$
and
$$
hr_{II}(D,pr)=0.
\tag{3.1.29}
$$

Using considerations above, In Appendix I,
we give Program  $h2$ for "GP/PARI"
calculator which for a fundamental discriminant $D<0$
and the genus $(D,\mu)$ calculates the vector
$$
(hr_I(D,\mu),\,hr_{II}(D,\mu),\,hnr(D,\mu)).
\tag{3.1.30}
$$
The "GP/PARI" calculator uses the Shank's method \cite{Sh}
to calculate the class numbers $h(D)$ of the discriminant $D$.
It is very fast: $O(|D|^{1/4})$ operations. We code the invariant $\mu$
by the non-negative integer $\mu$ having the binary decomposition
$$
\mu=\mu_{p_k}\mu_{p_{k-1}} \dots \mu_{p_1}
\tag{3.1.31}
$$
where $p_1,\,\dots ,\,p_k$ are all odd prime divisors of $D$ in increasing
order.

\subhead
3.2. The number $h$ of non-reflective classes of central symmetries of
main hyperbolic lattices of the rank $3$
\endsubhead

We denote by $S$ a main hyperbolic (i.e. of the signature
$(1,k)$) lattice with square-free determinant
and of the rank 3. Remind that {\it main} means that the lattice $S$
is even if the determinant $d=\det(S)$ is even. If $d$ is odd, then $S$
is necessarily odd.

By Proposition 2.2.4, the lattice $S$ is defined by its invariants
$(d,\eta)$ where $d=\det(S)$ and the invariant $\eta$ is the map
$\eta: p\mapsto \eta_p\in\{0,\,1\}$ of the set of all odd prime
divisors of $d$. The $\eta$ is defined by the condition
$$
b_{S_p}\cong b_{\theta_p}^{(p)}(p),\ \
\left( {\theta_p\over p}\right)=(-1)^{\eta_p}.
\tag{3.2.1}
$$
Here $q_S$ and $b_S$ denote discriminant forms of the lattice $S$.
Here and in what follows we use discriminant forms
technique and notations in \cite{N1}. We will especially often using
theorems of existence of a lattice with a given discriminant form
(Theorems 1.10.1 and 1.16.5 in \cite{N1}).

For even $d$ the lattice $S$ is even and the
discriminant quadratic form $q_{S_2}=q_{\theta_2}^{(2)}(2)$
where $\theta_2\equiv \pm 1 \mod 4$. We denote
$$
\theta_2\mod 8\equiv
\pm 1\mod 8\ \ \text{if\ }\ \theta_2\equiv\pm 1\mod 4.
\tag{3.2.2}
$$
We have
$$
\sum_{odd\ p\vert d}{[(1-p)+4\eta_p]} +\theta_2\mod 8\equiv -1\mod 8.
$$
Here $\text{sign}~q_{\theta_p}^{(p)}(p)\equiv (1-p)+4\eta_p\mod 8$
and $\text{sign}~q_{\theta_2}^{(2)}(2)\equiv \theta_2\mod 8$.
Since $\theta_2\mod 8\equiv \pm 1\mod 8$, we get
$$
\sum_{odd\ p\vert d}{[(1-p)+4\eta_p]}\equiv 0,-2\mod 8
\tag{3.2.3}
$$
and the invariant $\theta_2$ is defined by
$$
\theta_2\mod 8\equiv -\sum_{odd\ p\vert d}{[(1-p)+4\eta_p]}-1\mod 8.
\tag{3.2.4}
$$
Here \thetag{3.2.3} is the condition of existence of a main
hyperbolic lattice $S$ with the invariants $(d,\eta)$ if
$d$ is even. If $d$ is odd, the lattice $S$ is odd, and there are
no condition of existence of $S$. It always does exist.
Thus, {\it for any square-free natural number $d$ and
any map $\eta$ of the set of all odd prime divisors of $d$ into
$\{0,\,1\}$, there exists a main hyperbolic lattice $S$ of
the rank $3$ with invariants $(d,\,\eta)$ if and only if
$$
\sum_{odd\ p\vert d}{[(1-p)+4\eta_p]} \equiv 0,\,-2 \mod 8
\tag{3.2.5}
$$
if $d$ is even.}

We consider primitive elements $f\in S$ such that $n=f^2>0$.
We consider them up to $\pm f$, thus we can suppose that
$f \in V^+(S)$  where $V^+(S)$ is the light-cone of $S$
(see Sect. 1.1). We consider the elements $f$ such that
there exists an automorphism $u_f$ of $S$
which is identical on $f$ and is
$-1$ on the orthogonal negative definite lattice $K=f^\perp_S$.
Then
$$
\text{either}\ S=\bz f\oplus K\ \text{or}\
S=[\bz f\oplus K,\ (f\oplus k)/2],\  k\in K.
\tag{3.2.6}
$$
This automorphism is called a {\it central symmetry} of $S$.
Geometrically, it is the central symmetry at the point
$\br_{++}f$ of the hyperbolic space $\La(S)=V^+(S)/\br_{++}$
defined by $S$. Vice versa, any $\phi\in O^+(S)$ acting as a
central symmetry in $\La(S)$, has the form $\phi=u_f$.
By considerations over $\bz_2$, one can see
(using \thetag{3.2.6}) that either $n$ is odd or
$n\equiv 2 \mod 4$. (One can also use that
$u_f:x\mapsto -x+(2(x,f)/f^2)f$,
$x\in S$, and $f^2\vert 2(S,f)$.) Let $n=2^kn_1$
where $(n_1,2)=1$, $k=0,\,1$. Then by \thetag{3.2.6},
$$
n_1\vert d,\ \text{and\ } k=1\ \text{if\ } d\ \text{is even, and\ }
\left({n/p\over p}\right)=(-1)^{\eta_p}\  \forall\ \text{odd\ }p\vert n.
\tag{3.2.7}
$$
The central symmetry $u_f$ is called {\it reflective} if there exists a
reflection $s_\delta$ of $S$ in $\delta\perp f$ (i.e., $\delta \in K$).
Otherwise, $u_f$ is called {\it non-reflective}. Geometrically,
the central symmetry $u_f$ is reflective if and only if the center
$\br_{++}f$ belongs to a mirror of the reflection group
$W(S)$. We want to calculate
the number $hnr(S)$ of non-reflective
central symmetries $u_f$ of $S$ up to conjugation in
$O(S)$ (i.e., the number of non-reflective classes
of the central symmetries). Since $\pm 1$ belongs to the
center of $O(S)$, it is sufficient to consider $u_f$ up to
$O^+(S)=\{\phi \in O(S)\ \mid \phi(V^+(S))=V^+(S)\}$.
Acting by the reflection group $W(S)$,
we can always suppose that the center $\br_{++}f \in \M$
where $\M$ is a fundamental polyhedron of $W(S)$ (see Sect. 1.1).
If $u_f$ is non-reflective, it is equivalent that
$u_f\in A(P(\M)_{\pr})=\{\phi \in O^+(S)\ \mid \ \phi(\M)=\M\}$.
That two non-reflective central symmetries of $S$ are conjugate
by $O^+(S)$ if and only if they are conjugate by
$A(P(\M)_{\pr})$. Thus, the number $h=hnr(S)$ of
classes of non-reflective central symmetries of $S$ is
the same invariant $h=h(S)$ of $S$ which we have introduced in
Sect. 2.3. The central symmetry $u_f$ is defined by the
element $f$ (it is defined by $u_f$ up to $\pm f$). Thus, it is
sufficient to study all these $f$ up to action of $O(S)$. Obviously,
$(n,K)$ are invariants of
$f$. Here $K=f^\perp_S$ is a negative binary lattice
which we consider up to isomorphism.

First, we calculate invariants of
genus of the central symmetries $u_f$, $f\in S$ above.
Here two elements $f_1 \in S$ and
$f_2 \in S$ have the same {\it genus} if they are
conjugate over $\bz_p$ for all prime $p$.
Then we calculate the number of
non-reflective classes $u_f$ using class numbers of
binary positive lattices. We had considered them in Sect. 3.1.

\smallpagebreak

Case I: $d\equiv 1 \mod 2$, equivalently $S$ is odd.
Then there are two cases.

Case (I,\,II): $K=f^\perp$ is even. Then $n$ is odd because $S$
is unimodular odd over $\bz_2$, and then $f$ is a characteristic
element of $S$ (i.e. $(f,x)\equiv (x,x)\mod 2$ for any $x\in S$).
By \thetag{3.2.6}, we then have
$S=\bz f\oplus K$. It follows that $n\vert d$,
$\det K=d/n$, and the discriminant quadratic form
$q_{K_p}$, $p\vert (d/n)$, is equal to
$$
q_{K_p}=q_{S_p}=q_{\theta_p}^{(p)}(p),\ \
\left( {\theta_p\over p} \right)=(-1)^{\eta_p},
\ \ p\vert (d/n).
$$
By discriminant form technique \cite{N1},
existence of $K$ is equivalent to
$$
\sum_{p\mid (d/n)}[(1-p)+4\eta_p]\equiv -2 \mod 8.
\tag{3.2.8}
$$
Conditions \thetag{3.2.7} and \thetag{3.2.8} are equivalent to
existence of $u_f$ with $f^2=n$. The lattice $K(-1)$ is an even
fundamental positive binary lattice
of the discriminant $\discr K=-\det K=-d/n$.
Thus, $-d/n\equiv 1 \mod 4$. The genus of $K(-1)$ is equal to
$(-d/n,\epsilon(p)+\eta_p)$. Here
$\left({\,-1\,\over p}\right)=(-1)^{\epsilon(p)}$, and it
is known that $\epsilon(p)\equiv (p-1)/2 \mod 2$. To be shorter,
here we denote by $\epsilon(p)+\eta_p$ the map $p\mapsto
\epsilon(p)+\eta_p$
where $p$ runs through all odd prime $p\vert (d/n)$.
It follows that the number of non-reflective classes of
central symmetries $f\in S$ of the type
(I,\,II) is equal to
$$
hnr_{II}(S)=
\sum_{n\mid d\ \&\ (3.2.7)\atop \&\ (3.2.8)}
{hnr(-d/n,\epsilon(p)+\eta_p)},
\tag{3.2.9}
$$
where we consider the sum by all $n$ such that $n\vert d$ and
\thetag{3.2.7}, \thetag{3.2.8} are valid.

Case (I,\,I): $K=f^\perp$ is odd. By \thetag{3.2.6}, \thetag{3.2.7},
then $n\vert 2d$ and $\det K=(n,2)^2d/n$, $\discr K=-4(n,2)^2d/n$.
It follows that $K(-1)$ is a fundamental positive binary lattice.

For odd $p\vert \det K$, we have
$$
q_{K_p}=q_{S_p}=q_{\theta_p}^{(p)}(p),\ \
\left( {\theta_p\over p} \right)=(-1)^{\eta_p},
$$
and
$$
b_{K_2}=\cases
0,            &\text{if $n$ is odd},\\
b_1^{(2)}(2), &\text{if $n$ is even}
\endcases .
$$
If $n$ is odd, existence of $K$ is equivalent to
$$
-2 \mod 8 \in \sum_{p\vert (d/n)}{[(1-p)+4\eta_p]} +
4\omega(d/n)+\{0,\, \pm 2\}\mod 8.
$$
It is equivalent to
$$
\sum_{p\vert (d/n)}{[(1-p)+4\eta_p]}+4\omega(d/n) \not\equiv
2 \mod 8.
\tag{3.2.10}
$$
Here $\omega(k)\equiv (k^2-1)/8\mod 2$, it is known that
$\left({\,2\,\over p}\right)=(-1)^{\omega(p)}$.
For even $n$ there are no conditions of existence of $K$.
The lattice $K$ does always exist. The
genus of $K(-1)$ is equal to $(-4(n,2)^2d/n,\epsilon(p)+\eta_p)$.

Thus, this case is characterized by the condition:
$d$ is odd, $n\vert 2d$, we have
$\bz f\oplus K\subset S$ where $K$ is an odd fundamental
binary lattice of the discriminant $-4(n,2)^2d/n$.
If $n$ is odd, $\bz f \oplus K=S$. If $n$ is even,
the lattice $S$ is generated by $\bz f \oplus K$ and
$u=(f+k)/2$ where $k \in K$ is a primitive element with
the property: $(k,K)\equiv 0 \mod 2$.
The element $k\mod 2K$, and the overlattice
$S$ are defined uniquely. For both these
cases, any automorphism of $K$ can be extended to the automorphism of
$S$ identical on $f$. It follows that the number of non-reflective
classes of $f$ of the type (I,\,I) is equal to
$$
hnr_I(S)=\sum_{n\mid d\ \&\ (3.2.7) \atop \&\ (3.2.10)}
{hnr(-4d/n,\epsilon(p)+\eta_p)}+
\sum_{n=2n_1\mid 2d\ \atop \&  (3.2.7)}{hnr(-16d/n,\epsilon(p)+\eta_p)}.
\tag{3.2.11}
$$

As a result, we get that for the case I (equivalently,
when $d$ is odd) the full number of non-reflective
classes of central symmetries $f\in S$ is equal to
$$
\split
hnr(S)&=
\sum_{n\mid d\ \&\ (3.2.7)\ \&\atop (3.2.8)}
{hnr(-d/n,\epsilon(p)+\eta_p)}\\
&+\sum_{n\mid d\ \&\ (3.2.7)\atop \&\  (3.2.10)}
{\hskip-12pt hnr(-4d/n,\epsilon(p)+\eta_p)}+
\sum_{n=2n_1\mid 2d\ \atop \&\ (3.2.7)}
{hnr(-16d/n,\epsilon(p)+\eta_p)}.
\endsplit
\tag{3.2.12}
$$

\smallpagebreak

Case II: $d$ is even, equivalently, the lattice $S$ is even.
Then $n=2n_1$ where $n_1\vert (d/2)$.
This case is also divided in two cases:

Case (II,\,II): the lattice $K=f^\perp$ has an odd determinant.
By \thetag{3.2.6}, we have $S=\bz f\oplus K$.
Then
$q_{S_2}=q_{n/2}^{(2)}(2)$ and $n_1=n/2\equiv \theta_2 \mod 4$.
By \thetag{3.2.4},
$$
n/2 \equiv -\sum_{odd\ p\vert d}{[(1-p)+4\eta_p]}-1 \mod 4.
\tag{3.2.13}
$$
We have $\det K=d/n$ and,
for $p\vert (d/n)$,
$$
q_{K_p}=q_{S_p}=
q_{\theta_p}^{(p)}(p),\ \
\left( {\theta_p\over p} \right)=(-1)^{\eta_p}.
$$
Existence of $K$ is equivalent to
$\sum_{p\vert (d/n)}{[(1-p)+4\eta_p]}\equiv -2 \mod 8$. This
follows from the condition \thetag{3.2.2} (or \thetag{3.2.3})
of existence of $S$ and \thetag{3.2.7}, \thetag{3.2.13}.
Thus, condition of existence of $K$ for this case is \thetag{3.2.13}
together with \thetag{3.2.7}.

Thus, this case is characterized by the condition:
$d$ is even, $n\vert d$ is even, $S=\bz f\oplus K$, where
$K(-1)$ is an even fundamental binary lattice of
the genus $(-d/n, \epsilon(p)+\eta_p)$. Thus, the number of
non-reflective classes $f$ of the type (II,\,II) is equal to
$$
hnr_{II}(S)=\sum_{n=2n_1\vert d\ \&\ (3.2.7)\ \& \atop (3.2.13)}
{hnr(-d/n, \epsilon(p)+\eta_p)}.
\tag{3.2.14}
$$

Case (II,\,I): the lattice $K=f^\perp$ has even determinant.
The discriminant form $q_{[f]_2}=q_{n/2}^{(2)}(2)$.
By \thetag{3.2.6}, we have
$$
q_{K_2}=q_{-n/2}^{(2)}(2)\oplus q_{\theta_2}^{(2)}(2).
$$
It follows that the lattice $K(1/2)$
is odd fundamental, $K$ has determinant $4d/n$, and for
$p\vert (d/n)$ one has
$$
q_{K_p}=q_{S_p}=
q_{\theta_p}^{(p)}(p),\ \
\left( {\theta_p\ \over\ p\ } \right)=(-1)^{\eta_p}.
$$
Existence of the lattice $K$ is equivalent to
existence of the lattice $S$.
The odd fundamental binary lattice $K(-1/2)$ has the
determinant $d/n$, and  the discriminant $-4d/n$, and
for $p\vert (d/n)$ one has
$$
b_{K(1/2)_p}=b_{2\theta_p}^{(p)}(p),\ \
\left( {2\theta_p\ \over\ p\ } \right)=(-1)^{\eta_p+\omega(p)}.
$$

Thus, for this case,
$d$ is even, $f^2=n\vert d$ is even,
$K=T(-2)$ where $T$ is an odd fundamental positive binary lattice of
the discriminant $-4d/n$ and of the genus
$(-4d/n, \epsilon(p)+\omega(p)+\eta_p)$. The lattice
$S$ is an overlattice $\bz f\oplus K\subset S$ of the index $2$
generated by $(f+k)/2$ where $k\in K$ satisfies
$(k,k)\equiv -n \mod 8$. If $d/n\equiv -1 \mod 4$, then
$k\mod 2K$ and $S$ are unique. If $d/n\equiv 1\mod 4$,
there are exactly two different
elements $k$ which give the same lattice $S$. If $\discr T=-4$,
elements $k\mod 2K$ are conjugate by $O(K)$. If $\discr T<-4$,
elements $k\mod 2K$ are not conjugate and give different
classes of $f \in S$. Further we consider two cases:

Case $d/n \equiv -1 \mod 4$. The number of such
non-reflective classes $f\in S$ of the type (II,\,I) is equal to
$$
hnr_{I,-1}(S)=\sum_{n=2n_1\vert d\ \&\ (3.2.7)\& \atop n_1\equiv -d/2\mod
4}
{hnr(-4d/n,\epsilon(p)+\omega(p)+\eta_p)}.
\tag{3.2.15}
$$

Case $d/n \equiv 1 \mod 4$. Then the class $f\in S$ is non-reflective
if and only if the lattice $K$ does not have a reflection which can
be extended identically on $\bz f$ to give an automorphism of $S$.
If $K$ is non-ambiguous, then $O(K)=\{\pm 1\}$, and we then
get two classes of $f$ corresponding to two
different choices of $k\mod 2K$. If $d/n\not=1$
and $T$ is ambiguous of the type II,
then $O(K)$ has order $4$, only $\pm 1 \in O(K)$ preserve
the element $k\mod 2K$, and reflections of $K$ change places two
possible different elements $k\mod 2K$.
Thus we get exactly one non-reflective
class $f\in S$. If $d/n\not=1$ and $T$ is ambiguous of type I,
then reflections of $T$ are identical on $K/2K$ and can be extended
on $S$ identically on $f$. Thus, they give reflective classes
$f\in S$. If $d/n=1$, the lattice $T=K(1,1;I)=K(2,2;II)$, and it has
reflections of both types I and II. This case does not give non-reflective
classes $f  \in S$. Thus, for this case, the number of non-reflective
classes $f\in S$ of type (II,\,I) is equal to
$$
\split
&hnr_{I,1}(S)=\\
&\sum_{n=2n_1\vert d\ \&\ (3.2.7)\ \&\atop
{n_1\equiv d/2 \mod 4\ \&\ n<d}}
{\hskip-15pt [2 hnr(-4d/n,\epsilon(p)+\omega(p)+\eta_p)+
hr_{II}(-4d/n,\epsilon(p)+\omega(p)+\eta_p)]}.
\endsplit
\tag{3.2.16}
$$

As a result, we get for the case II (equivalently, when $d$ is even) that
the full number of non-reflective classes of central symmetries $f\in S$
is equal to
$$
\split
hnr(S)&=
\hskip-10pt\sum_{n=2n_1\vert d\ \&\ (3.2.7)\ \& \atop (3.2.13)}
{hnr(-d/n, \epsilon(p)+\eta_p)}\\
+&\sum_{n=2n_1\vert d\ \&\ (3.2.7)\& \atop n_1\equiv -d/2\mod 4}
{\hskip-5pt hnr(-4d/n,\epsilon(p)+\omega(p)+\eta_p)}\\
+&\hskip-20pt \sum_{n=2n_1\vert d\ \&\ (3.2.7)\ \& \atop
{n_1\equiv d/2 \mod 4\ \&\ n<d}}
{\hskip-20pt [2 hnr(-4d/n,\epsilon(p)+\omega(p)+\eta_p)+
hr_{II}(-4d/n,\epsilon(p)+\omega(p)+\eta_p)]}.
\endsplit
\tag{3.2.17}
$$
As a final result, we get

\proclaim{Theorem 3.2.1} Let $d$ be a square-free natural number and
$\eta:p\mapsto \{0,\,1\}$ a map of all odd prime divisors
$p\vert d$ into
$\{0,\,1\}$. Then there exists a main hyperbolic lattice $S$ of the
rank $3$ with the square-free determinant $d$ and the invariant $\eta$
(see \thetag{3.2.1})
if and only if for the even $d$ the congruence \thetag{3.2.5} is valid.
The number $h=hnr(S)=hnr(d,\eta)$
of classes of non-reflective central symmetries of the lattice
$S$ with the invariants $(d,\eta)$ is given by
\thetag{3.2.12} for the odd $d$ and by \thetag{3.2.17} for the even $d$.
\endproclaim

Below we will code the invariant $\eta$ by the non-negative
integer $\eta$ having the binary decomposition
$$
\eta=\eta_{p_t} \dots \eta_{p_1}
\tag{3.2.18}
$$
where $p_1,\dots,p_t$ are all odd prime divisors of $d$ in
increasing order.

In Appendix: Programs, we give the Program 2: h3 for
"GP/PARI" calculator which using Theorem 3.2.1
and Program 1: h2 (see Sect. 3.1) calculates the invariant
$h=h(S)=hnr(d,\eta)$ if a main hyperbolic
lattice $S$ of the rank $3$ with the invariants $(d,\eta)$ does
exist (otherwise, the result will be unreasonable). Using the first
statement of the Theorem 3.2.1 and the Program 2, we give Program 3: refh3
which gives all pairs of invariants $(d,\eta)$ such that $d\le N$,
there exists a main hyperbolic lattice of the rank $3$
with the invariants $(d,\eta)$ and
the invariant $hnr(d,\eta)\le 1$. Using Program 3, we found
all these pairs $(d,\eta)$ such that $d\le 100000$. The
result is given in Table 3 (Part I) and contains $206$ lattices.
Thus, we get

\proclaim{Theorem 3.2.2} Table 3 (Part I) gives the complete list
(it has 206 lattices) of main hyperbolic lattices $S$ with
square-free determinant $d\le 100000$ and of the rank $3$ such that
the invariant $hnr(S)\le 1$. In Table 3 we give invariants
$d$, $\eta$, $h=hnr(S)$, the matrix of the lattice $S$ and the
reflective type of the lattice $S$.
\endproclaim

The greatest $d$ of lattices $S$ of the Table 3 is equal to $4466$
in spite we did calculations up to $100000$. Thus, it is very likely
that we have the statement

\proclaim{Conjecture 3.2.3} Table 3 (Part I) gives the complete list
of main hyperbolic lattices $S$ with square-free determinant and
of the rank $3$ such that the invariant $hnr(S)\le 1$.
\endproclaim

In Sect. 5, we shall use Theorem 3.2.2 to find all
$(d,\eta)$ corresponding to elliptically or parabolically
reflective hyperbolic lattices $S$ since for elliptically or
parabolically reflective hyperbolic lattices the invariant $h\le 1$.

Similarly we calculated all pairs $(d,\eta)$ such that
$hnr(d,\eta)=2$ and $d\le 100000$. (One should change in two
places of Program 3 $h\le 1$ by $h=2$.) The list contains
$259$ pairs $(d,\eta)$. The last $10$ pairs having the largest $d$
are: $(4290,1)$, $(4326,2)$, $(4902,4)$, $(4991,7)$,
$(5226,0)$, $(5334,2)$, $(6006,2)$, $(7590,8)$, $(10374,2)$, $(29526,2)$.
It is very likely that this list also contains all pairs $(d,\eta)$
with $h=2$. We shall use the list of these $259$ lattices and lattices
of Table 3 to find in Part III all hyperbolically reflective main
hyperbolic lattices of the rank $3$.  All of them must have the
invariant $h=0,\,2$.

\head
4. Narrow places of elliptic and parabolic convex
polygons on hyperbolic plane, types of polygons.
Application to reflective lattices
\endhead

\subhead
4.1. Narrow places of elliptic convex polygons on the hyperbolic plane
\endsubhead

We remind (see Sect. 1.1) that a convex polyhedron in
a hyperbolic space is called
{\it elliptic} if it is a convex envelope of a finite set of points
(some of them at infinity) and it is non-degenerate. In this section
we shall consider only elliptic (i.e. ordinary finite) convex polyhedra
and often shall omit the word ``elliptic''.

Here we follow the general method (of narrow places of polyhedra)
suggested in \cite{N4}, \cite{N5} for proving finiteness results
about arithmetic reflection groups in hyperbolic spaces.
On the other hand, we shall prove much more delicate and exact
statements, which are important for exact classification.
Our estimates in \cite{N4}, \cite{N5} were universal, they did not
depend on angles of fundamental polyhedra. Here we get estimates
which depend on angles of polyhedra which makes the narrow places of
polyhedra method much more efficient. Our estimates
here are optimal, we belive that one cannot impove them.

We restrict by $2$-dimensional case of narrow places of elliptic
polygons on the hyperbolic plane, but all results can be
easily generalized (like in \cite{N4}, \cite{N5}) on
elliptic polyhedra of arbitrary dimension in hyperbolic spaces.

We shall often use the following trivial but important for us
statement (certainly, it is well-known):

\proclaim{Lemma 4.1.1} ${\sin{y}\over \sin{x}} < {y\over x}$
if $0< x \le \pi/2$ and $x < y$.
\endproclaim

\demo{Proof} The function $sin{x}/x$ is decreasing if $0\le x \le \pi /2$.
In particular, $1\ge sin{x}/x \ge 2/\pi$ if $0\le x\le \pi /2$.
It follows the lemma for $0<x\le y\le \pi/2$. If $0<x \le \pi/2\le y$,
we get
$$
{\sin{y}\over \sin{x}}\le {1\over sin{x}}\le
{1\over x(2/\pi)}\le {\pi /2 \over x} \le {y\over x}.
$$
This proves the statement.
\enddemo

\proclaim{Lemma 4.1.2} Let $(AB)$ and $(CD)$ are two lines on a
hyperbolic plane with terminals $A$, $B$, $C$, $D$ at infinity,
and $O$ a point on the hyperbolic plane which
does not belong to each line $(AB)$ and $(CD)$ and orientations of the
triangles $AOB$ and $COD$ coincide.
We consider angles
$\theta_1=AOB$, $\theta_2=COD$ and $\theta_{12}=BOC$. Let
$\delta_1$ and $\delta_2$ are orthogonal vectors with
square $-2$ to lines $(AB)$ and $(CD)$ respectively such that
$O$ is contained in both half-planes $\Ha_{\delta_1}^+$ and
$\Ha_{\delta_2}^+$.

Then
$$
(\delta_1,\delta_2)=
4{\sin{{\theta_1+\theta_{12}\over 2}}\sin{{\theta_2+\theta_{12}\over 2}}
\over \sin{{\theta_1\over 2}}\sin{{\theta_2\over2}}}-2.
$$

As a corollary, we get :

1) If lines $(AB)$ and $(CD)$ do not intersect each other, then
$$
2 \cosh{\rho }=
(\delta_1,\delta_2)=
4{\sin{{\theta_1+\theta_{12}\over 2}}\sin{{\theta_2+\theta_{12}\over 2}}
\over \sin{{\theta_1\over 2}}\sin{{\theta_2\over2}}}-2
$$
where $\rho$ is the distance between lines $(AB)$ and $(CD)$
(here and in what follows we normalize the curvature $\kappa=-1$).

2) If lines $(AB)$ and $(CD)$ define an angle $\alpha$
containing $O$, then
$$
2 \cos{\alpha}=
(\delta_1,\delta_2)=
4{\sin{{\theta_1-\theta_{21}\over 2}}\sin{{\theta_2-\theta_{21}\over 2}}
\over \sin{{\theta_1\over 2}}\sin{{\theta_2\over2}}}-2
$$
where $\theta_{21}=-\theta_{12}=COB$.
\endproclaim

\demo{Proof} We can correspond to two lines $(AB)$, $(CD)$ and
a connected component (containing $O$) of the complement to these two
lines in the hyperbolic plane two invariants up to motions of the
hyperbolic plane. First invariant is equal to $(\delta_1,\delta_2)$
and uses Klein model of the hyperbolic plane. Second invariant is
equal to the cross ratio $[A:D:C:B]$ where we suppose that
orientations of the triangles $ABO$  and $CDO$ coincide. This invariant
uses Poincar\'e model of the hyperbolic plane. On the other hand,
it is clear that any of these two invariants defines the triplet
($(AB)$, $(CD)$, the connected component to their complement
(containing $O$)) up to motions of the hyperbolic plane.
It follows that there exists a function f(x) such that
$(\delta_1,\delta_2)=f([A:D:C:B])$. One can check that
$f(x)=4x-2$.
\enddemo

\smallpagebreak

Below we shall consider an equation
$$
(u-x)(v-x)=auv
\tag{4.1.1}
$$
where $u,v\ge 0$ and $0\le a=\cos^2{\alpha\over 2}=(1+\cos{\alpha})/2\le
1$.
Equivalently we have the equation
$$
x^2-(u+v)x+uv(1-a)=0.
$$
Its smallest root $x=g(\alpha,u,v)$ equals
$$
x=g(\alpha,u,v)={u+v-\sqrt{a(u+v)^2+(1-a)(u-v)^2}\over 2}
\tag{4.1.2}
$$
where $a=a(\alpha)=\cos^2{\alpha\over 2}=(1+\cos{\alpha})/2$.
In particular, for $u=v$,
$$
g(\alpha,u,u)=u(1-\sqrt{a})=u(1-\cos{\alpha\over 2}).
\tag{4.1.3}
$$

The function $g(\alpha,u,v)$ has the following properties
(we shall not use them further but they are important):

\proclaim{Proposition 4.1.3}
$u\ge g(\alpha,u,v)\ge 0$, $v \ge g(\alpha,u,v)\ge 0$,
$g(\alpha,u,v)'_u\ge 0$,  $g(\alpha,u,v)'_v\ge 0$,
$\left({u-g(\alpha,u,v)\over u}\right)'_u\ge 0$,
$\left({v-g(\alpha,u,v)\over v}\right)'_v\ge 0$.
\endproclaim

\demo{Proof} Suppose that $u\ge v$. The expression
$(u-x)(v-x)$ equals $uv\ge auv$ if $x=0$ and
$(u-x)(v-x)$ equals $0$ if $x=u$. Thus the smallest solution
$g(a,u,v)$ of the equation $(u-x)(v-x)=auv$ satisfies
$v\ge u\ge g(a,u,v)$.

We have
$2g(a,u,v)'_u=1-(a(u+v)+(1-a)(u-v))/\sqrt{a(u+v)^2+(1-a)(u-v)^2}$
$\ge 0$ if $(a(u+v)+(1-a)(u-v))/\sqrt{a(u+v)^2+(1-a)(u-v)^2}\le 1$.
Equivalently, for $-1\le t=(u-v)/(u+v)\le 1$ we should prove that
$(a+(1-a)t)/\sqrt{a+(1-a)t^2}\le 1$. If $-1\le t\le 0$, this is obvious.
For $0\le t\le 1$, we have
$(a+(1-a)t^2)^{3\over 2}\left((a+(1-a)t)/\sqrt{a+(1-a)t^2}\right)'_t=
(1-a)(a+(1-a)t^2)-(1-a)t(a+(1-a)t)=(1-a)a(1-t)\ge 0$. It
follows that $(a+(1-a)t)/\sqrt{a+(1-a)t^2}\le
(a+(1-a))/\sqrt{a+(1-a)}=1$.

By definition, we have
$(u-g(\alpha,u,v))/u=av/(v-g(\alpha,u,v)$. It follows that
\newline
$\left((u-g(\alpha,u,v))/u\right)'_u\ge 0$ because
$g(\alpha,u,v)'_u\ge 0$.

It finishes the proof.
\enddemo

\proclaim{Theorem 4.1.4 (about the narrow place of type (I))}
For any elliptic convex polygon $\M$ on a
hyperbolic plane there exist its four consecutive vertices
$A_0$, $A_1$, $A_2$ and $A_3$ (where $A_0=A_3$ if $\M$ is
a triangle) such that
for orthogonal vectors $\delta_1$, $\delta_2$ and $\delta_3$ to
lines $(A_0A_1)$, $(A_1A_2)$ and $(A_2A_3)$ respectively directed
outwards of $\M$ and with
$\delta_1^2=\delta_2^2=\delta_3^2=-2$
one has $(\delta_1,\delta_2)=2\cos{\alpha_1}$,
$(\delta_2,\delta_3)=2\cos{\alpha_2}$ and
$$
\text{either\ \ } (\delta_1,\delta_3)\le 2\ \ \text{or}\ \
(\delta_1,\delta_3)
<4(\cos{\alpha_1\over 2}+\cos{\alpha_2\over 2})^2-2\le 14
\tag{4.1.4}
$$
where $\alpha_1=A_0A_1A_2$ and $\alpha_2=A_1A_2A_3$.

Moreover, the Gram graph of $\{\delta_1,\,\delta_2,\,\delta_3\}$
is not connected (i.e. this set
is union of two non-empty orthogonal subsets)
if and only if $\alpha_1=\alpha_2={\pi\over 2}$.
\endproclaim

\demo{Proof}
To prove Theorem,
we take a point $O$ inside of $\M=A_1A_2...A_n$.
Let $B_{i1}$ and $B_{i2}$ are terminals at infinity
of the line $l_i=(A_{i-1}A_i)$ where $B_{i1}$, $A_{i-1}$, $A_i$ and
$B_{i2}$ are four consecutive points of the line. We introduce
angles $\alpha_i=A_{i-1}A_iA_{i+1}$,
$\theta_i=B_{i1}OB_{i2}$ and $\theta_{(i+1)i}=B_{(i+1)1}OB_{i2}$.

By Lemma 4.1.2,
$$
2 \cos{\alpha_i}=
4{\sin{{\theta_i-\theta_{(i+1)i}\over 2}}
\sin{{\theta_{i+1}-\theta_{(i+1)i}\over 2}}
\over \sin{{\theta_i\over 2}}\sin{{\theta_{i+1}\over2}}}-2.
$$
Equivalently,
$$
{\sin{{\theta_i\over 2}}\sin{{\theta_{i+1}\over2}}\over
\sin{{\theta_i-\theta_{(i+1)i}\over 2}}
\sin{{\theta_{i+1}-\theta_{(i+1)i}}\over 2}}={2\over 1+\cos{\alpha_i}}.
\tag{4.1.5}
$$
By Lemma 4.1.1, we get
$$
{(\theta_i-\theta_{(i+1)i})(\theta_{i+1}-\theta_{(i+1)i})
\over
\theta_i\theta_{i+1}}<{1+\cos{\alpha_i} \over 2}=
\cos^2{{\alpha_i\over 2}}.
\tag{4.1.6}
$$
It follows that
$$
\theta_{(i+1)i}>g(\alpha_i,\theta_i,\theta_{i+1})
\tag{4.1.7}
$$
(see \thetag{4.1.2}).

To prove Theorem 4.1.4, we choose a line $l_i$ with
the minimal angle $\theta_i$. Let this line be $l_2$. Thus
$$
\theta_2=\min_i{\theta_i}.
\tag{4.1.8}
$$
It follows that $\theta_1\ge \theta_2$ and $\theta_3\ge \theta_2$.
We shall then prove the inequalities \thetag{4.1.4}.

If lines $l_1,l_3$ intersect, we have
$(\delta_1,\delta_3)\le 2$ and \thetag{4.1.4} is valid.
Suppose that the lines $l_1$ and $l_3$ do not intersect.

There exists a line $l_1'=(B_{11}'B_{12}')$ with terminals
$B_{11}'$ and $B_{12}'$
at infinity such that the line
$l_1'$ is contained in $\Ha_{-\delta_1}^+$, points $B_{11}$, $B_{11}'$
and $B_{12}$, $B_{12}'$ are contained in the same half-planes bounded
by $l_2$, the line $l_1'$ has the same angle
$\alpha_1=B_{11}'A_1'B_{22}$ (as $l_1$) with the line $l_2$ (we denote by
$A_1'$ their intersection point),
and the angle $B_{11}'OB_{12}'$ of $l_1'$ is equal to $\theta_2$. We
denote by $\delta_1'$ the orthogonal vector to $l_1'$
directed outwards of $\M$ and with $(\delta_1')^2=-2$, and we denote
$\theta_{21}'=B_{21}OB_{12}'$.

Similarly,
there exists a line $l_3'=(B_{31}'B_{32}')$ with terminals
$B_{31}'$ and $B_{32}'$
at infinity such that the line
$l_3'$ is contained in $\Ha_{-\delta_3}^+$, points $B_{31}$, $B_{31}'$
and $B_{32}$, $B_{32}'$ are contained in the same half-planes bounded
by $l_2$, the line $l_3'$ has the same angle
$\alpha_2=B_{21}A_2'B_{32}'$ (as $l_3$) with the line $l_2$ (we denote by
$A_2'$ their point of intersection),
and the angle $B_{31}'OB_{32}'$ of $l_3'$ is equal to $\theta_2$. We
denote by $\delta_3'$ the orthogonal vector to $l_3'$
directed outwards of $\M$ and with $(\delta_3')^2=-2$, and we denote
$\theta_{32}'=B_{31}'OB_{22}$.

Since the lines $l_1$ and $l_3$ do not intersect, by our construction,
any interval with terminals at $l_1'$ and $l_3'$ intersects
both lines $l_1$ and $l_3$. It follows that distance between lines
$l_1'$ and $l_3'$ is greater than distance between lines $l_1$
and $l_3$. It follows $(\delta_1,\delta_3)\le (\delta_1',\delta_3')$.
It is sufficient to prove \thetag{4.1.4} for lines $l_1'$ and $l_3'$.
(These geometrical considerations are related with properties in
Proposition 4.1.4 of the function $g(\alpha,u,v)$.)

By Lemmas 4.1.1, 4.1.2, and \thetag{4.1.7}, \thetag{4.1.3},
$$
(\delta_1,\,\delta_3)\le (\delta_1',\,\delta_3')=
4{\sin{\theta_2+\theta_2-\theta_{21}'-\theta_{32}'\over 2}
\sin{\theta_2+\theta_2-\theta_{21}'-\theta_{32}'\over 2}
\over
\sin{\theta_2\over 2}\sin{\theta_2\over 2}}-2<
$$
$$
4{(\theta_2+\theta_2-\theta_{21}'-\theta_{32}')
(\theta_2+\theta_2-\theta_{21}'-\theta_{32}')
\over
\theta_2\theta_2}-2\ <
$$
$$
4{(\theta_2 \hskip-2pt +\hskip-2pt\theta_2-g(\alpha_1,\theta_2,\theta_2)-
g(\alpha_2,\theta_2,\theta_2))
(\theta_2\hskip-2pt+\hskip-2pt\theta_2-g(\alpha_1,\theta_2,\theta_2)-
g(\alpha_2,\theta_2,\theta_2))
\over
\theta_2\theta_2}-2=
$$
$$
4(\sqrt{a_1}+\sqrt{a_2})^2-2=
4(\cos{\alpha_1\over 2}+\cos{\alpha_2\over 2})^2-2
$$
where $a_i=cos^2{\alpha_i\over 2}$. It proves \thetag{4.1.4}.

Elements $\delta_1,\delta_2,\delta_3$ generate the
hyperbolic 3-dimensional vector space defining the hyperbolic plane.
Otherwise their lines either have a common point or
are orthogonal to one line which is not the case.
If the Gram graph of these elements is
not connected, two of these elements generate a
$2$-dimensional hyperbolic vector subspace.
Elements $\delta_1$ and $\delta_2$ cannot generate a
hyperbolic $2$-dimensional vector subspace because
their orthogonal lines have a common point $\br_{++}h$ where
$h^2\ge 0$ and $h\not=0$. Then $(h,\,\delta_1)=(h,\,\delta_2)=0$.
The same is valid for $\delta_2$ and $\delta_3$.
Thus the Gram graph of $\{\delta_1,\,\delta_2,\,\delta_3\}$
is not connected if and only if $(\delta_1,\,\delta_2)=
(\delta_3,\delta_2)=0$. Equivalently, $\alpha_1=\alpha_2={\pi\over 2}$.
It finishes the proof of Theorem 4.1.4.
\enddemo

\proclaim{Theorem 4.1.5 (about narrow places of types (II) and (III))}
For any elliptic convex polygon $\M$ having more than $3$
vertices (i.e. it is different from a triangle) on a hyperbolic plane,
one of two possibilities (II) or (III) below is valid:

(II) There exist its five consecutive vertices
$A_0$, $A_1$, $A_2$, $A_3$ and $A_4$
(where $A_0=A_4$ if $\M$ is a quadrangle) such that
for orthogonal vectors $\delta_1$, $\delta_2$, $\delta_3$ and $\delta_4$
to lines $(A_0A_1)$, $(A_1A_2)$, $(A_2A_3)$  and $(A_3A_4)$
respectively directed outwards of $\M$ and with
$\delta_1^2=\delta_2^2=\delta_3^2=\delta_4^2=-2$,
one has $(\delta_i,\delta_{i+1})=2\cos{\alpha_i}$, $i=1,\,2,\,3$, and
$$
\text{either\ \ } (\delta_1,\delta_3)\le 2\ \ \text{or}\ \
(\delta_1,\delta_3)
<4(\cos{\alpha_1\over 2}+\cos{\alpha_2\over 2})^2-2\le 14,
\tag{4.1.9}
$$
$$
(\delta_1,\delta_4)<
$$
$$
\split
&4\,\max_{0\le t\le 1}
{{\left(\sqrt{a_1+(1-a_1)t^2}+\sqrt{a_2+(1-a_2)t^2}+
\sqrt{a_3+a_3t+t^2/4}\right)^2-{t^2\over 4}\over 1+t}}-2=\\
&4\,\max{\left(\left(\cos{\alpha_1\over 2}+\cos{\alpha_2\over 2}+
\cos{\alpha_3\over 2}\right)^2,
{\left(2+\sqrt{2\cos^2{\alpha_3\over 2}+{1\over 4}}\,\right)^2-
{1\over 4}\over 2}\right)}-2\le 34
\endsplit
\tag{4.1.10}
$$
where $\alpha_i=A_{i-1}A_iA_{i+1}$, $i=1,\,\,2,\,3$, and
$a_i=\cos^2{\alpha_i\over 2}$. Moreover,
the set $\{\delta_1,\,\delta_2,$ $\delta_3,\,\delta_4\}$ has
a connected Gram graph.

(III). There exist its six consecutive vertices
$A_0$, $A_1$, $A_2$, $A_3$, $A_4$ and $A_5$
(where $A_0=A_5$ if $\M$ is a pentagon) such that
for orthogonal vectors $\delta_1$, $\delta_2$, $\delta_3$,
$\delta_4$ and $\delta_5$
to lines $(A_0A_1)$, $(A_1A_2)$, $(A_2A_3)$, $(A_3A_4)$ and
$(A_4A_5)$
respectively directed
outwards of $\M$ and with
$\delta_1^2=\delta_2^2=\delta_3^2=\delta_4^2=\delta_5^2=-2$
one has $(\delta_i,\delta_{i+1})=2\cos{\alpha_i}$, $i=1,2,3,4$, and
$$
\text{either\ \ } (\delta_1,\delta_3)\le 2\ \ \text{or}\ \
(\delta_1,\delta_3)
<4(\cos{\alpha_1\over 2}+\cos{\alpha_2\over 2})^2-2\le 14,
\tag{4.1.11}
$$
$$
\text{either\ \ } (\delta_3,\delta_5)\le 2\ \ \text{or}\ \
(\delta_3,\delta_5)
<4(\cos{\alpha_3\over 2}+\cos{\alpha_4\over 2})^2-2\le 14,
\tag{4.1.12}
$$
and
$$
(\delta_1,\delta_5)<
$$
$$
\split
&4\,\max_{0\le t\le s\le 1}
{\left[\left(\left(\sqrt{a_1+(1-a_1)s^2}+\sqrt{a_2+(1-a_2)s^2}+
\right.\right.\right.}\\
&+\left.\sqrt{a_3+a_3(s-t)+{a_3(s-t)^2\over 4}+
{(1-a_3)(s+t)^2\over 4}}+\sqrt{a_4+(1-a_4)t^2}\,\right)^2-\\
&-\left.\left.{(s-t)^2\over 4}\right)/((1+s)(1+t))\right]-2=\\
&4\,\max\left[\left(\cos{\alpha_1\over 2}+\cos{\alpha_2\over 2}+
\cos{\alpha_3\over 2}+\cos{\alpha_4\over 2}\right)^2,\right.\\
&\left.{\left(2+\sqrt{2\cos^2{\alpha_3\over 2}+{1\over 4}}+
\cos{\alpha_4\over 2}\right)^2-{1\over 4}\over 2},\ 4\right]-2\ \le\  62
\endsplit
\tag{4.1.13}
$$
where $\alpha_i=A_{i-1}A_iA_{i+1}$, $i=1,\,2,\,3,\,4$, and
$a_i=\cos^2{\alpha_i\over 2}$. Moreover, the
set $\{\delta_1,\,\delta_2,\,\delta_3,\,\delta_4,\,\delta_5\}$ has a
connected Gram graph.
\endproclaim

\demo{Proof} Like for the proof of Theorem 4.1.4, we take a point
$O$ inside of $\M$ and introduce angles $\theta_i$, $\theta_{(i+1)i}$.
We choose such a consecutive numeration of vertices of $\M$ that
$$
\theta_2+\theta_3 =\min_i{(\theta_i+\theta_{i+1})}\  \text{and}\
\theta_2\le \theta_3.
\tag{4.1.14}
$$
Then $\theta_1\ge \theta_3 \ge \theta_2$
and, like for the proof of Theorem 4.1.4, we get
\thetag{4.1.9} and \thetag{4.1.11} (see considerations below about
cases (II) and (III)).

Below we consider two cases:

Case (ii): $\theta_4\ge (\theta_2+\theta_3)/2$. We then prove
\thetag{4.1.10}.

For $1\ge a_1,\,a_2,\,a_3\ge 0$ and $0\le t$ we introduce a function
$$
f_{p2}(a_1,\,a_2,\,a_3,\,t)=
$$
$$
={\left(\sqrt{a_1+(1-a_1)t^2}+\sqrt{a_2+(1-a_2)t^2}+
\sqrt{a_3+a_3t+t^2/4}\right)^2-{t^2\over 4}\over 1+t}.
\tag{4.1.15}
$$
We prove that
$$
(\delta_1,\delta_4)<4\max_{0\le t\le 1}{f_{p2}(a_1,\,a_2,\,a_3,\,t)}-2.
\tag{4.1.16}
$$
If $(\delta_1,\,\delta_4)\le 2$, it is true because
$4f_{p2}(a_1,\,a_2,\,a_3,\,1)-2=
2((2+\sqrt{2a_3+1/4})^2-1/4)-2\ge
2((2+1/2)^2-1/4)-2=10$.

Like for the proof of Theorem 4.1.4, we can suppose that
$\theta_1=\theta_3$ (instead of $\theta_1\ge \theta_3$) and
$\theta_4=(\theta_2+\theta_3)/2$
(instead of $\theta_4\ge (\theta_2+\theta_3)/2$). We denote
$c=(\theta_2+\theta_3)/2$ and $z=(\theta_3-\theta_2)/2\ge 0$.
Then $\theta_2=c-z$ and $\theta_3=c+z$.

Like for the proof of Theorem 4.1.4, we have
$$
((\delta_1,\,\delta_4)+2)/4<
$$
$$
{(\theta_3+\theta_2+\theta_3-\theta_{21}-\theta_{32}-\theta_{43})
(c+\theta_2+\theta_3-\theta_{21}-\theta_{32}-\theta_{43})
\over \theta_3c}<
$$
$$
(\theta_3+\theta_2+\theta_3-g(\alpha_1,\,\theta_3,\,\theta_2)-
g(\alpha_2,\,\theta_2,\,\theta_3)-g(\alpha_3,\,\theta_3,\,c))\times
$$
$$
(c+\theta_2+\theta_3-g(\alpha_1,\,\theta_3,\,\theta_2)-
g(\alpha_2,\,\theta_2,\,\theta_3)-g(\alpha_3,\,\theta_3,\,c))/
\theta_3c=
$$
$$
(z/2+\sqrt{a_1c^2+(1-a_1)z^2}+\sqrt{a_2c^2+(1-a_2)z^2}+
\sqrt{{z^2\over 4}+a_3cz+a_3c^2})\times
$$
$$
(-z/2+\sqrt{a_1c^2+(1-a_1)z^2}+\sqrt{a_2c^2+(1-a_2)z^2}+
\sqrt{{z^2\over 4}+a_3cz+a_3c^2})/(c+z)c=
$$
$$
{(\sqrt{a_1c^2+(1-a_1)z^2}+\sqrt{a_2c^2+(1-a_2)z^2}+
\sqrt{{z^2\over 4}+a_3cz+a_3c^2})^2-z^2/4
\over
(c+z)c}=
$$
$$
{\left(\sqrt{a_1+(1-a_1)t^2}+\sqrt{a_2+(1-a_2)t^2}+
\sqrt{{t^2\over 4}+a_3t+a_3}\right)^2-t^2/4
\over
1+t}
$$
where $t=z/c$ and $0\le t\le 1$. It proves \thetag{4.1.16}.

If one of angles $\alpha_1$ or $\alpha_2$ is not $\pi/2$,
Gram graph of $\{\delta_1,\,\delta_2,\,\delta_3\}$ is connected.
Then Gram graph of $\{\delta_1,\,\delta_2,\,\delta_3,\,\delta_4\}$ is
connected because $\delta_1,\,\delta_2,\,\delta_3$ generate the full
3-dimensional hyperbolic vector space. If $\alpha_1=\alpha_2=\pi/2$
(equivalently, $(\delta_1,\,\delta_2)=(\delta_3,\,\delta_2)=0$),
then $\delta_1$, $\delta_3$ generate a hyperbolic $2$-dimensional
subspace and $(\delta_1,\,\delta_3)$ $\not=0$. If Gram graph of
$\{\delta_1,\,\delta_2,\,\delta_3,\delta_4\}$ is not connected,
$(\delta_1,\,\delta_4)=(\delta_3,\,\delta_4)=0$. It follows that
$\delta_4=\lambda \delta_2$, $\lambda \in \br$. Then sides
of $\M$ orthogonal to $\delta_2$ and $\delta_4$ should
coincide. We get a contradiction.

Case (ii) gives polygons $\M$ which satisfy the condition
(II) of Theorem 4.1.5. We mention that by Theorem 4.1.4,
any quadrangle $\M$ satisfies the case (II) of Theorem 4.1.5.
Really, for a quadrangle, $(\delta_1,\delta_4)\le 2$.

Case (iii): $\theta_4\le (\theta_2+\theta_3)/2$ and $\M$ has more
than $4$ vertices. We denote $(\theta_2+\theta_3)/2=c$.

By \thetag{4.1.14}, $\theta_1\ge \theta_3$,
$\theta_4+\theta_5\ge (\theta_2+\theta_3)=2c$. It follows,
$\theta_5\ge c$.

Since $\theta_3\ge c$, $\theta_4\le c$ and $\theta_5\ge c$, like
in the proof of Theorem 4.1.4,
we get the inequality \thetag{4.1.12}.

For $1\ge a_1,\,a_2,\,a_3,\,a_4\ge 0$ and $0\le s,t$, we introduce a
function
$$
f_{p3}(a_1,\,a_2,\,a_3,a_4,\,s,\,t)=
$$
$$
\split
&\left(\left(\sqrt{a_1+(1-a_1)s^2}+\sqrt{a_2+(1-a_2)s^2}+
\right.\right.\\
&+\left.\sqrt{a_3+a_3(s-t)+{a_3(s-t)^2\over 4}+
{(1-a_3)(s+t)^2\over 4}}+\sqrt{a_4+(1-a_4)t^2}\,\right)^2-\\
&-\left.{(s-t)^2\over 4}\right)/((1+s)(1+t)).
\endsplit
\tag{4.1.17}
$$

We prove that
$$
(\delta_1,\,\delta_5)<4\max_{0\le t\le s\le 1}
{f_{p3}(a_1,\,a_2,\,a_3,\,a_4,\,s,\,t)}-2.
\tag{4.1.18}
$$
If $(\delta_1,\,\delta_5)\le 2$, it is true because
$f_{p3}(a_1,\,a_2,\,a_3,\,a_4,\,1,\,1)=4$. If $(\delta_1,\,\delta_5)>2$,
like for the proof of Theorem 4.1.4, we can assume that
$\theta_1=\theta_3$ (instead of $\theta_1\ge \theta_3$) and
$\theta_5=2c-\theta_4$ (instead of $\theta_5\ge 2c-\theta_4$).

We denote $c=(\theta_2+\theta_3)/2$, $z=(\theta_3-\theta_2)/2$ and
$w=(\theta_5-\theta_4)/2$. We have $\theta_1=\theta_3=c+z$,
$\theta_2=c-z$,  $\theta_4=c-w$, $\theta_5=c+w$. By definition,
$c\ge 0$, $z\ge 0$, $w\ge 0$. Moreover, $z\ge w$ because
$\theta_3+\theta_4=2c+z-w\ge \theta_2+\theta_3=2c$
(we use \thetag{4.1.14}).

Like for the proof ot Theorem 4.1.4, using Lemmas 4.1.1, 4.1.2, we get
$$
((\delta_1,\delta_5)+2)/4<
$$
$$
{(\theta_3 \hskip-1pt +\hskip-1pt\theta_2\hskip-1pt+
\hskip-1pt\theta_3\hskip-1pt+\hskip-1pt\theta_4\hskip-1pt-\hskip-1pt
\theta_{21}\hskip-1pt-\hskip-1pt\theta_{32}\hskip-1pt-\hskip-1pt
\theta_{43}\hskip-1pt-\hskip-1pt\theta_{54})
(\theta_5\hskip-1pt+\hskip-1pt\theta_2\hskip-1pt+\hskip-1pt\theta_3
\hskip-1pt+\hskip-1pt\theta_4\hskip-1pt-\hskip-1pt
\theta_{21}\hskip-1pt-\hskip-1pt\theta_{32}\hskip-1pt-\hskip-1pt
\theta_{43}\hskip-1pt-\hskip-1pt\theta_{54})
\over \theta_3\theta_5}<
$$
$$
\split
&(\theta_3+\theta_2+\theta_3+\theta_4-
g(\alpha_1,\theta_3,\theta_2)-g(\alpha_2,\theta_2,\theta_3)-
g(\alpha_3,\theta_3,\theta_4)-g(\alpha_4,\theta_4,\theta_5))\times\\
&(\theta_5\hskip-2pt+\hskip-2pt\theta_2\hskip-2pt+\hskip-2pt
\theta_3\hskip-2pt+\hskip-2pt\theta_4\hskip-2pt-\hskip-2pt
g(\alpha_1,\theta_3,\theta_2)\hskip-2pt-\hskip-2pt
g(\alpha_2,\theta_2,\theta_3)\hskip-2pt-\hskip-2pt
g(\alpha_3,\theta_3,\theta_4)\hskip-2pt-\hskip-2pt
g(\alpha_4,\theta_4,\theta_5))
/(\theta_3\theta_5)=
\endsplit
$$
$$
\split
&\left(\left(\sqrt{a_1c^2+(1-a_1)z^2}+\sqrt{a_2c^2+(1-a_2)z^2}+
\right.\right.\\
&\sqrt{a_3c^2+a_3c(z-w)+{a_3(z-w)^2\over 4}+{(1-a_3)(z+w)^2\over 4}}+\\
&\left.\left.\sqrt{a_4c^2+(1-a_4)w^2}\right)^2-{(z-w)^2\over 4}\right)/
((c+z)(c+w)).
\endsplit
$$
Denoting $s=z/c$ and $t=w/c$, we get
$$
((\delta_1,\delta_5)+2)/ 4<f_{p3}(a_1,\,a_2,\,a_3,\,a_4,\,s,\,t)
$$
where $0\le t\le s\le 1$.

Similarly to the case (ii), one can easily prove that Gram graph of
$\{\delta_1,\,\delta_2,\,\delta_3,\,\delta_4,$ $\delta_5\}$ is
connected. Thus, for the case (iii), we get the case (III) of Theorem.

To finish the proof of Theorem, in Lemmas 4.1.6 and 4.1.7 below,
we find
\newline
$\max_{0\le t\le 1}{f_{p2}(a_1,\,a_2,\,a_3,\,t)}$
and
$\max_{0\le t\le s\le 1}{f_{p3}(a_1,\,a_2,\,a_3,\,a_4,\,s,\,t)}$.
\enddemo

\proclaim{Lemma 4.1.6} For $0\le a_1,\,a_2,\,a_3\le 1$
and $t\ge 0$ the function
$$
f_{p2}(a_1,\,a_2,\,a_3,\,t)=
$$
$$
={\left(\sqrt{a_1+(1-a_1)t^2}+\sqrt{a_2+(1-a_2)t^2}+
\sqrt{a_3+a_3t+t^2/4}\right)^2-{t^2\over 4}\over 1+t}
$$
has the maximum
$$
\split
&\max_{0\le t\le 1}{f_{p2}(a_1,\,a_2,\,a_3,\,t)}=
\max \left[f_{p2}(a_1,\,a_2,\,a_3,\,0),\ f_{p2}(a_1,\,a_2,\,a_3,\,1)
\right]=\\
&\max\left[
\left(\sqrt{a_1}+\sqrt{a_2}+\sqrt{a_3}\right)^2,\
{\left(2+\sqrt{2a_3+{1\over 4}}\,\right)^2-
{1\over 4}\over 2}\right].
\endsplit
$$
\endproclaim

\proclaim{Lemma 4.1.7} For $0\le a_1,\,a_2,\,a_3,\,a_4\le 1$
and $s,\,t\ge 0$ the function
$$
f_{p3}(a_1,\,a_2,\,a_3,a_4,\,s,\,t)=
$$
$$
\split
&\left(\left(\sqrt{a_1+(1-a_1)s^2}+\sqrt{a_2+(1-a_2)s^2}+
\right.\right.\\
&+\left.\sqrt{a_3+a_3(s-t)+{a_3(s-t)^2\over 4}+
{(1-a_3)(s+t)^2\over 4}}+\sqrt{a_4+(1-a_4)t^2}\,\right)^2-\\
&-\left.{(s-t)^2\over 4}\right)/((1+s)(1+t))
\endsplit
$$
has the maximum
$$
\split
&\max_{0\le t\le s\le 1}{f_{p3}(a_1,\,a_2,\,a_3,\,a_4,\,s,\,t)}=\\
&\max \left[f_{p3}(a_1,\,a_2,\,a_3,\,a_4,0,\,0),\
f_{p3}(a_1,\,a_2,\,a_3,\,a_4,\,1,\,0),\
f_{p3}(a_1,\,a_2,\,a_3,\,a_4,\,1,\,1)
\right]=\\
&\max\left[\left(\sqrt{a_1}+\sqrt{a_2}+\sqrt{a_3}+\sqrt{a_4}\right)^2,\
{\left(2+\sqrt{2a_3+{1\over 4}}+a_4\right)^2-
{1\over 4}\over 2},\ 4\right]
\endsplit
$$
\endproclaim

\demo{Proof of Lemma 4.1.6 for
$1\ge a_1,\,a_2\ge 0.2723$, and $1\ge a_3\ge 5/32=0.15625$} We shall use
Theorem 4.1.5 only for $1\ge a_i\ge 1/2$ when a polygon $\M$ has
acute angles. Therefore, we restrict proving Lemma 4.1.6 for the
parameters $a_i$ above.

To prove Lemma 4.1.6 for parameters $a_1,\,a_2,\,a_3$,
it is sufficient to show that
$$
f_{p2}(a_1,\,a_2,\,a_3,\,t)''_{tt}\ge 0,\ \ \text{if}\  0\le t\le 1.
\tag{4.1.19}
$$

For $1\ge a\ge 0$ and $0\le t\le 1$ we consider
$$
u={a+(1-a)t^2\over 1+t}.
\tag{4.1.20}
$$
We have
$$
u'_t={(1-a)t^2+2(1-a)t-a\over (1+t)^2}=1-a-{1\over (1+t)^2},\ \ \
u''_{tt}={2\over (1+t)^3}.
\tag{4.1.21}
$$
It follows,
$$
(u^{{1\over 2}})''_{tt}=
u^{-{3\over 2}}\left(-{1\over 4}(u'_t)^2+{1\over 2}u\,u''_{tt}\right)=
$$
$$
\split
&\left({a+(1-a)t^2\over t+1}\right)^{-{3\over 2}}\times\\
&{-(1/4)(1-a)^2t^4 -(1-a)^2t^3+(3/2)(a-a^2)t^2+(a-a^2)t+
a-a^2/4\over (1+t)^4}\ge
\endsplit
$$
$$
\split
&\left({a+(1-a)t^2\over t+1}\right)^{-{3\over 2}}\times\\
&{-(1/4)(1-a)^2t^2 -(1-a)^2t^2+(3/2)(a-a^2)t^2+(a-a^2)t+
a-a^2/4\over (1+t)^4}=
\endsplit
$$
$$
\left({a+(1-a)t^2\over t+1}\right)^{-{3\over 2}}\times
{(1-a)({11\over 4}a-{5\over 4})t^2+(a-a^2)t+
a-{1\over 4}a^2\over (1+t)^4}\ge
$$
$$
\left({a+(1-a)t^2\over t+1}\right)^{-{3\over 2}}\times(1+t)^{-4}\times
\cases
a-{1\over 4}a^2, &\text{if ${1\over 3}\le a\le 1$}\\
a-{1\over 4}a^2, &\text
{if $0\le a \le {1\over 3}$ \& $0\le t\le {4a\over 5-11a}$}\\
-{5\over 4}+6a-4a^2,
&\text{if $0\le a\le {1\over 3}$ \& ${4a\over 5-11a}\le t\le 1$}.
\endcases
$$
It follows,
$$
(u^{{1\over 2}})''_{tt}\ge
$$
$$
\left({a+(1-a)t^2\over t+1}\right)^{-{3\over 2}}\times(1+t)^{-4}\times
\cases
a-{1\over 4}a^2, &\text{if ${1\over 3}\le a\le 1$}\\
a-{1\over 4}a^2, &\text
{if $0\le a\le {1\over 3}$ \& $0\le t\le {4a\over 5-11a}$}\\
-{5\over 4}+6a-4a^2,
&\text{if $0\le a\le {1\over 3}$ \& ${4a\over 5-11a}\le t\le 1$}
\endcases
$$
$$
\ge
\left({a+(1-a)t^2\over t+1}\right)^{-{3\over 2}}\times
\min (a-{1\over 4}a^2,-{5\over 4}+6a-4a^2)/(1+t)^4.
\tag{4.1.22}
$$
Here $\min{(a-{1\over 4}a^2,\ -{5\over 4}+6a-4a^2)}\ge 0$ if
$a\ge 1/4$. If follows
$$
(u^{{1\over 2}})''_{tt}\ge
{\min (a-{1\over 4}a^2,-{5\over 4}+6a-4a^2)\over (1+t)^{5\over 2}}\ge 0
\tag{4.1.23}
$$
if $1\ge a\ge {1\over 4}$ and $0\le t\le 1$.

For $1\ge a\ge 0$ and $0\le t\le 1$ we consider
$$
v={t^2/4+at+a\over 1+t}.
\tag{4.1.24}
$$
We have
$$
v'_t={t^2/4+t/2\over (1+t)^2}={1\over 4}-{1\over 4(1+t)^2},\ \
v''_{tt}={1\over 2(1+t)^3}.
\tag{4.1.25}
$$
It follows
$$
(v^{1\over 2})''_{tt}=
\left({t^2/4+at+a\over 1+t}\right)^{-{3\over2}}\times
{-(1/64)t^4-(1/16)t^3+(1/4)at+(1/4)a\over
(1+t)^4}\ge
$$
$$
\left({t^2/4+at+a\over 1+t}\right)^{-{3\over2}}\times
{-(1/64)t^3-(1/16)t^3+(1/4)at+(1/4)a\over
(1+t)^4}=
$$
$$
\left({t^2/4+at+a\over 1+t}\right)^{-{3\over2}}\times
{-(5/64)t^3+(1/4)at+(1/4)a\over
(1+t)^4}\ge
$$
$$
\left({t^2/4+at+a\over 1+t}\right)^{-{3\over2}}\times
\cases
{1\over 4}a, &\text{if ${5\over 16}\le a\le 1$}\\
{1\over 4}a, &\text
{if $0\le a \le {5\over 16}$ \& $0\le t\le 4\sqrt{{a\over 5}}$}\\
-{5\over 64}+{1\over 2}a,
&\text{if $0\le a\le {5\over 16}$ \& $4\sqrt{{a\over 5}}\le t\le 1$}.
\endcases
$$
It follows
$$
(v^{{1\over 2}})''_{tt}\ge
$$
$$
\left({t^2/4+at+a\over 1+t}\right)^{-{3\over2}}\times
(1+t)^{-4}\times
\cases
{1\over 4}a, &\text{if ${5\over 16}\le a\le 1$}\\
{1\over 4}a, &\text
{if $0\le a \le {5\over 16}$ \& $0\le t\le 4\sqrt{{a\over 5}}$}\\
-{5\over 64}+{1\over 2}a,
&\text{if $0\le a\le {5\over 16}$ \& $4\sqrt{{a\over 5}}\le t\le 1$}
\endcases
$$
$$
\ge
\left({t^2/4+at+a\over 1+t}\right)^{-{3\over2}}\times
{\min{({1\over 4}a,\ -{5\over 64}+{1\over 2}a)}
\over
(1+t)^4}.
\tag{4.1.26}
$$
Here $\min{({1\over 4}a,\ -{5\over 64}+{1\over 2}a)}\ge 0$ if
$a\ge {5\over 32}$.
If follows
$$
(v^{1\over 2})''_{tt}\ge
(1/4+2a)^{-{3\over 2}}\times
{\min{({1\over 4}a,\ -{5\over 64}+{1\over 2}a)}\over (1+t)^{5\over 2}}
\tag{4.1.27}
$$
if $1\ge a\ge 5/32=0.15625$, $0\le t\le 1$.

We have
$$
\left({t^2/4\over 1+t}\right)'_t={1\over 4}-{1\over 4(1+t)^2},\ \
\left({t^2/4\over 1+t}\right)''_{tt}={1\over 2(1+t)^3}.
\tag{4.1.28}
$$

We denote for $1\ge a_i\ge 0$ and $0\le t\le 1$,
$$
u_1={a_1+(1-a_1)t^2\over 1+t},\ \ u_2={a_2+(1-a_2)t^2\over 1+t}\ \
u_3={t^2/4+a_3t+a_3\over 1+t}
$$
and
$$
w=\sqrt{u_1}+\sqrt{u_2}+\sqrt{u_3}.
$$
We have
$$
(w^2)''_{tt}=(w'_t)^2+2ww''_{tt}\ge 2ww''_{tt}.
\tag{4.1.29}
$$
From \thetag{4.1.23}, \thetag{4.1.27}, \thetag{4.1.28} and
\thetag{4.1.29},
for $1\ge a_1,\,a_2,\ge 1/4$ and
$1\ge a_3\ge 5/32$, we get for $f_{p2}=w^2-(1/4)t^2/(1+t)$,
$$
f_{p2}(a_1,\,a_2,\,a_3,\,t)''_{tt}\ge
$$
$$
\split
&\left[2\left((a_1+(1-a_1)t^2)^{1\over 2}+(a_2+(1-a_2)t^2)^{1\over 2}+
(t^2/4+a_3t+a_3)^{1\over 2}\right)\times\right.\\
&\left(\min{(a_1-{1\over 4}a_1^2,\ -{5\over 4}+6a_1-4a_1^2)}+
\min{(a_2-{1\over 4}a_2^2,\ -{5\over 4}+6a_2-4a_2^2)}+\right.\\
&\left.\left.(1/4+2a_3)^{-{3\over 2}}
\min{({1\over 4}a_3,\ -{5\over 64}+{1\over 2}a_3)}\right)-
{1\over 2}\right]/(1+t)^3\ge
\endsplit
$$
$$
\split
&\left[2\left(\sqrt{a_1}+\sqrt{a_2}+\sqrt{a_3}\right)\times\right.\\
&\left(\min{(a_1-{1\over 4}a_1^2,\ -{5\over 4}+6a_1-4a_1^2)}+
\min{(a_2-{1\over 4}a_2^2,\ -{5\over 4}+6a_2-4a_2^2)}+\right.\\
&\left.\left.(1/4+2a_3)^{-{3\over 2}}
\min{({1\over 4}a_3,\ -{5\over 64}+{1\over 2}a_3)}\right)-
{1\over 2}\right]/(1+t)^3\ge 0
\endsplit
\tag{4.1.30}
$$
if
$$
\split
&\left(\sqrt{a_1}+\sqrt{a_2}+\sqrt{a_3}\right)\times\\
&\left(\min{(a_1-{1\over 4}a_1^2,\ -{5\over 4}+6a_1-4a_1^2)}+
\min{(a_2-{1\over 4}a_2^2,\ -{5\over 4}+6a_2-4a_2^2)}+\right.\\
&\left.(1/4+2a_3)^{-{3\over 2}}
\min{({1\over 4}a_3,\ -{5\over 64}+{1\over 2}a_3)}\right)
\ge {1\over 4}.
\endsplit
\tag{4.1.31}
$$
Thus, \thetag{4.1.19} is true if $1\ge a_1,\,a_2\ge 1/4$,
$1\ge a_3\ge 5/32$ and \thetag{4.1.31} is valid.
For $1\ge a_1,\,a_2\ge a\ge 1/4$ and arbirary $1\ge a_3\ge 5/32$,
\thetag{4.1.31} is valid if
$$
2(2\sqrt{a}+\sqrt{5\over 32})(-{5\over 4}+6a-4a^2)\ge {1\over 4}.
\tag{4.1.32}
$$
It is true if $1\ge a\ge 0.2723$.

It proves Lemma 4.1.6 for $1\ge a_1,\,a_2\ge 0.2723$ and
$1\ge a_3\ge 5/32=0.15625$.
\enddemo

\smallpagebreak

\demo{Proof of Lemma 4.1.7 for $1\ge a_1,\,a_2,\,a_3,\,a_4\ge 0.37646$}
We shall use
Theorem 4.1.5 only for $1\ge a_i\ge 1/2$ when a polygon $\M$ has
acute angles. Therefore, we restrict proving Lemma 4.1.7 for the
parameters $a_i$ above.

To prove Lemma 4.1.7 for parameters $a_1,\,a_2,\,a_3,\,a_4$,
it is sufficient to show that
$$
f_{p3}(a_1,\,a_2,\,a_3,\,a_4,\,s,\,t)''_{tt}\ge 0,\ \ \text{if}\
0\le t\le s\le 1,
\tag{4.1.33}
$$
$$
f_{p3}(a_1,\,a_2,\,a_3,\,a_4,\,s,\,0)''_{ss}\ge 0,\ \ \text{if}\
0\le s\le 1;
\tag{4.1.34}
$$
and
$$
f_{p3}(a_1,\,a_2,\,a_3,\,a_4,\,s,\,s)''_{ss}\ge 0,\ \ \text{if}\
0\le s\le 1.
\tag{4.1.35}
$$

We have
$$
f_{p3}(a_1,\,a_2,\,a_3,\,a_4,\,s,\,0)=
$$
$$
{\left(\sqrt{a_1+(1-a_1)s^2}+\sqrt{a_2+(1-a_2)s^2}+
\sqrt{a_3+a_3s+{s^2\over 4}}+\sqrt{a_4}\right)^2-
{s^2\over 4}\over 1+s}.
\tag{4.1.36}
$$
Like for the proof of Lemma 4.1.6, it follows for
$1\ge a_1,\,a_2\ge {1\over 4}$, $1\ge a_3\ge {5\over 32}$,
$1\ge a_4\ge 0$ that
$$
f_{p3}(a_1,\,a_2,\,a_3,\,a_4,\,s,\,0)''_{ss}\ge
$$
$$
\split
&\left[2\left((a_1+(1-a_1)s^2)^{1\over 2}+(a_2+(1-a_2)s^2)^{1\over 2}+
(s^2/4+a_3s+a_3)^{1\over 2}+a_4^{1\over 2}\right)\times\right.\\
&\left(\min{(a_1-{1\over 4}a_1^2,\ -{5\over 4}+6a_1-4a_1^2)}+
\min{(a_2-{1\over 4}a_2^2,\ -{5\over 4}+6a_2-4a_2^2)}+\right.\\
&\left.\left.(1/4+2a_3)^{-{3\over 2}}
\min{({1\over 4}a_3,\ -{5\over 64}+{1\over 2}a_3)}\right)-
{1\over 2}\right]/(1+s)^3\ge
\endsplit
$$
$$
\split
&\left[2\left(\sqrt{a_1}+\sqrt{a_2}+\sqrt{a_3}+\sqrt{a_4}
\right)\times\right.\\
&\left(\min{(a_1-{1\over 4}a_1^2,\ -{5\over 4}+6a_1-4a_1^2)}+
\min{(a_2-{1\over 4}a_2^2,\ -{5\over 4}+6a_2-4a_2^2)}+\right.\\
&\left.\left.(1/4+2a_3)^{-{3\over 2}}
\min{({1\over 4}a_3,\ -{5\over 64}+{1\over 2}a_3)}\right)-
{1\over 2}\right]/(1+s)^3\ge 0
\endsplit
\tag{4.1.37}
$$
if
$$
\split
&\left(\sqrt{a_1}+\sqrt{a_2}+\sqrt{a_3}+\sqrt{a_4}\right)\times\\
&\left(\min{(a_1-{1\over 4}a_1^2,\ -{5\over 4}+6a_1-4a_1^2)}+
\min{(a_2-{1\over 4}a_2^2,\ -{5\over 4}+6a_2-4a_2^2)}+\right.\\
&\left.(1/4+2a_3)^{-{3\over 2}}
\min{({1\over 4}a_3,\ -{5\over 64}+{1\over 2}a_3)}\right)
\ge {1\over 4}.
\endsplit
\tag{4.1.38}
$$
Thus, \thetag{4.1.34} is true if $1\ge a_1,\,a_2\ge 1/4$,
$1\ge a_3\ge 5/32$, $1\ge a_4\ge 0$  and \thetag{4.1.38} is valid.

For $1\ge a_1,\,a_2\ge a\ge 1/4$, arbirary $1\ge a_3\ge 5/32$ and
arbitrary $1\ge a_4\ge 0$, we have
\thetag{4.1.38} if \thetag{4.1.32} is valid. It takes place
if $1\ge a\ge 0.2723$.

We have
$$
f_{p3}(a_1,\,a_2,\,a_3,\,a_4,\,s,\,s)=
$$
$$
{\left(\sqrt{a_1+(1\hskip-1pt -\hskip-1pt a_1)s^2}+
\sqrt{a_2+(1\hskip-1pt -\hskip-1pt a_2)s^2}+
\sqrt{a_3+(1\hskip-1pt-\hskip-1pta_3)s^2}+
\sqrt{a_4+(1\hskip-1pt-\hskip-1pta_4)s^2}\right)^2
\over (1+s)^2}.
\tag{4.1.39}
$$
For $1\ge a\ge 0$ and $0\le s\le 1$ we denote
$$
e={a+(1-a)s^2\over (1+s)^2}
\tag{4.1.40}
$$
We have
$$
\split
&e'_s={-2a+2s(1-a)\over (1+s)^3}={2s\over (1+s)^3}-{2a\over (1+s)^2},\\
&e''_{ss}={2+4a+(4a-4)s\over (1+s)^4}={2-4s\over (1+s)^4}+{4a\over
(1+s)^3}.
\endsplit
\tag{4.1.41}
$$
It follows,
$$
(e^{{1\over 2}})''_{ss}=
e^{-{3\over 2}}\left(-{1\over 4}(e'_s)^2+{1\over 2}e\,e''_{ss}\right)=
$$
$$
e^{-{3\over 2}}\times
{(-2a^2+4a-2)s^3+(-3a^2+3a)s^2+(a^2+a)\over (1+s)^6}\ge
$$
$$
e^{-{3\over 2}}\times
{(-2a^2+4a-2)s^2+(-3a^2+3a)s^2+(a^2+a)\over (1+s)^6}=
$$
$$
\split
&e^{-{3\over 2}}\times
{(-5a^2+7a-2)s^2+(a^2+a)\over (1+s)^6}\ge\\
&e^{-{3\over 2}}\times
\cases
a^2+a, &\text{if $1\ge a\ge {2\over 5}$}\\
-4a^2+8a-2, &\text{if ${2\over 5}\ge a\ge 0$}
\endcases\ \ge\\
&e^{-{3\over 2}}\times{\min(a^2+a,\ -4a^2+8a-2)\over (1+s)^6}.
\endsplit
\tag{4.1.42}
$$
Here $\min(a^2+a,\ -4a^2+8a-2)\ge 0$ if
$a\ge (2-\sqrt{2})/2=0.292893...$ .
It follows
$$
(e^{{1\over 2}})''_{ss}\ \ge\
{\min(a^2+a,\ -4a^2+8a-2)\over (1+s)^3}\ge 0\ \
\text{if}\ 1\ge a\ge (2-\sqrt{2})/2.
\tag{4.1.43}
$$
For $1\ge a_1,\,a_2,\,a_3,\,a_4\ge 0$ we denote
$$
e_i={a_i+(1-a_i)s^2\over (1+s)^2}.
\tag{4.1.44}
$$
We get
$$
f_{p3}(a_1,\,a_2,\,a_3,\,a_4,\,s,\,s)''_{ss} \ge
$$
$$
2(\sqrt{u_1}+\sqrt{u_2}+\sqrt{u_3}+\sqrt{u_4})
((u_1^{1\over 2})''_{ss}+(u_2^{1\over 2})''_{ss}+
(u_3^{1\over 2})''_{ss}+(u_4^{1\over 2})''_{ss})\ge
$$
$$
\split
&2\left(\sqrt{a_1\hskip-1pt+\hskip-1pt(1\hskip-1pt -\hskip-1pt a_1)s^2}+
\sqrt{a_2\hskip-1pt+\hskip-1pt(1\hskip-1pt -\hskip-1pt a_2)s^2}+
\sqrt{a_3\hskip-1pt+\hskip-1pt(1\hskip-1pt-\hskip-1pta_3)s^2}+
\sqrt{a_4\hskip-1pt+\hskip-1pt(1\hskip-1pt-\hskip-1pta_4)s^2}\right)
\times\\
&\left((\min(a_1^2+a_1,\,
-4a_1^2+8a_1-2)+\min(a_2^2+a_2,\,-4a_2^2+8a_2-2)+\right.\\
&\left.+\min(a_3^2+a_3,\,-4a_3^2+8a_3-2)+
\min(a_4^2+a_4,\,-4a_4^2+8a_4-2)\right)/(1+s)^4\ge
\endsplit
$$
$$
\split
&2\left(\sqrt{a_1}+
\sqrt{a_2}+
\sqrt{a_3}+
\sqrt{a_4}\right)
\times\\
&\left((\min(a_1^2+a_1,\,
-4a_1^2+8a_1-2)+\min(a_2^2+a_2,\,-4a_2^2+8a_2-2)+\right.\\
&\left.+\min(a_3^2+a_3,\,-4a_3^2+8a_3-2)+
\min(a_4^2+a_4,\,-4a_4^2+8a_4-2)\right)/(1+s)^4\ge 0
\endsplit
\tag{4.1.45}
$$
if $1\ge a_i\ge (2-\sqrt{2})/2=0.292893...$ .

Let us prove \thetag{4.1.33} for $1\ge a_1,\,a_2,\,a_3,\,a_4\ge 0.37235$.
We denote
$$
f={a+a(s-t)+{a(s-t)^2\over 4}+{(1-a)(s+t)^2\over 4}\over t+1}=
{a+as+{s^2\over 4}+t(-a-as+{s\over 2})+{t^2\over 4}\over t+1}.
$$
We have
$$
f'_t={-2a-2as+{s\over 2}-{s^2\over 4}+{t\over 2}+{t^2\over 4}\over
(t+1)^2}=
{1\over 4}-{(1-s)^2+8a(s+1)\over 4(t+1)^2},
$$
$$
f''_{tt}={(1-s)^2+8a(s+1)\over 2(t+1)^3}.
$$
For $1\ge a\ge 0$ and $0\le t\le s\le 1$, we have
$$
(f^{1\over 2})''_{tt}=
f^{-{3\over 2}}\left(-{1\over 4}(f'_t)^2+{1\over 2}f\,f''_{tt}\right)=
$$
$$
\split
&f^{-{3\over 2}}(-{1\over 64}t^4-{1\over 16}t^3+
({3\over 32}s^2+({3\over 4}a-{3\over 16})s+{3\over 4}a)t^2+\\
&((-{1\over 4}a+{1\over 8})s^3+(-2a^2+{5\over 4}a-{3\over 16})s^2+
(-4a^2+{7\over 4}a)s+(-2a^2+{1\over 4}a))t+\\
&({3\over 64}s^4+({1\over 2}a-{1\over 16})s^3+(a^2+{1\over 2}a)s^2+
(2a^2+{1\over 4}a)s+(a^2+{1\over 4}a))/(t+1)^4\ge
\endsplit
$$
$$
\split
&f^{-{3\over 2}}({3\over 32}s^2+({3\over 4}a-{3\over 16})s+{3\over 4}a-
{5\over 64})t^2+\\
&((-{1\over 4}a+{1\over 8})s^3+(-2a^2+{5\over 4}a-{3\over 16})s^2+
(-4a^2+{7\over 4}a)s+(-2a^2+{1\over 4}a))t+\\
&({3\over 64}s^4+({1\over 2}a-{1\over 16})s^3+(a^2+{1\over 2}a)s^2+
(2a^2+{1\over 4}a)s+(a^2+{1\over 4}a))/(t+1)^4.
\endsplit
$$
For $1\ge a\ge {1\over 4}$ and $0\le t\le s\le 1$ we have
$$
\split
&({3\over 32}s^2+({3\over 4}a-{3\over 16})s+{3\over 4}a-{5\over 64})t^2+\\
&((-{1\over 4}a+{1\over 8})s^3+(-2a^2+{5\over 4}a-{3\over 16})s^2+
(-4a^2+{7\over 4}a)s+(-2a^2+{1\over 4}a))t+\\
&({3\over 64}s^4+({1\over 2}a-{1\over 16})s^3+(a^2+{1\over 2}a)s^2+
(2a^2+{1\over 4}a)s+(a^2+{1\over 4}a))\ge
\endsplit
$$
$$
\split
&((-{1\over 4}a+{1\over 8})s^3+(-2a^2+{5\over 4}a-{3\over 16})s^2+
(-4a^2+{7\over 4}a)s+(-2a^2+{1\over 4}a))t+\\
&({3\over 64}s^4+({1\over 2}a-{1\over 16})s^3+(a^2+{1\over 2}a)s^2+
(2a^2+{1\over 4}a)s+(a^2+{1\over 4}a))\ge
\endsplit
$$
$$
\split
&\min\left({3\over 64}s^4+({1\over 2}a-{1\over 16})s^3+(a^2+{1\over
2}a)s^2+
(2a^2+{1\over 4}a)s+(a^2+{1\over 4}a),\right.\\
&\left.(-{1\over 4}a+{11\over 64})s^4+(-2a^2+{7\over 4}a-{1\over 4})s^3+
(-3a^2+{9\over 4}a)s^2+{1\over 2}as+(a^2+{1\over 4}a)\right).\\
\endsplit
\tag{4.1.46}
$$
Here ${3\over 64}s^4+({1\over 2}a-{1\over 16})s^3+(a^2+{1\over 2}a)s^2+
(2a^2+{1\over 4}a)s+(a^2+{1\over 4}a)\ge a^2+{1\over 4}a$.
Let us prove that
$(-{1\over 4}a+{11\over 64})s^4+(-2a^2+{7\over 4}a-{1\over 4})s^3+
(-3a^2+{9\over 4}a)s^2+{1\over 2}as+(a^2+{1\over 4}a)\ge
{11\over 64}a^2+{1\over 4}a$. Equivalently, we should prove that
$(-{1\over 4}a+{11\over 64})s^4+(-2a^2+{7\over 4}a-{1\over 4})s^3+
(-3a^2+{9\over 4}a)s^2+{1\over 2}as+{53\over 64}a^2\ge 0$.

Assume that  $1\ge a\ge {11\over 16}$.
Then $-{1\over 4}a+{11\over 64}\le 0$ and
$(-{1\over 4}a+{11\over 64})s^4+(-2a^2+{7\over 4}a-{1\over 4})s^3+
(-3a^2+{9\over 4}a)s^2+{1\over 2}as+{53\over 64}a^2\ge
r(a,s)=(-2a^2+{3\over 2}a-{5\over 64})s^3+
(-3a^2+{9\over 4}a)s^2+{1\over 2}as+{53\over 64}a^2$.
We have $r(a,0)\ge 0$ and
$r(a,1)=-{267\over 64}a^2+{17\over 4}a-{5\over 64}\ge 0$.
It follows that $r(a,s)\ge 0$ if
${1\over 2}r''_{ss}=3(-2a^2+{3\over 2}a-{5\over 64})s+
(-3a^2+{9\over 4}a)\le 0$.
The last is
valid if ${1\over 2}r''_{ss}(a,0)=-3a^2+{9\over 4}a\le 0$ and
${1\over 2}r''_{ss}(a,1)=3(-2a^2+{3\over 2}a-{5\over 64})+
(-3a^2+{9\over 4}a)\le 0$.
It is true if $1 \ge a\ge {3\over 4}$ when $-3a^2+{9\over 4}a\le 0$.
Suppose that ${3\over 4}\ge a\ge {11\over 16}$.
Then $-3a^2+{9\over 4}a\ge 0$.
If also $-2a^2+{3\over 2}a-{5\over 64}\ge 0$, then $r(a,s)\ge 0$.
If $-2a^2+{3\over 2}a-{5\over 64} \le 0$, we have
$r(a,s)\ge ((-2a^2+{3\over 2}a-{5\over 64})+
(-3a^2+{9\over 4}a))s^2+{1\over 2}as+{53\over 64}a^2$.
If $(-2a^2+{3\over 2}a-{5\over 64})+
(-3a^2+{9\over 4}a)\le 0$, we again have that
$r(a,s)\ge \min(r(a,0),r(a,1))\ge 0$. If
$(-2a^2+{3\over 2}a-{5\over 64})+
(-3a^2+{9\over 4}a)\ge 0$,
we obviously have $r(a,s)\ge 0$.

Assume that ${11\over 64}\ge a\ge {1/4}$. Then
$(-{1\over 4}a+{11\over 64})s^4+(-2a^2+{7\over 4}a-{1\over 4})s^3+
(-3a^2+{9\over 4}a)s^2+{1\over 2}as+{53\over 64}a^2\ge
r(a,s)=(-2a^2+{7\over 4}a-{1\over 4})s^3+
(-3a^2+{9\over 4}a)s^2+{1\over 2}as+{53\over 64}a^2$. Arguing
with the function $r(a,s)$ like above, we prove that
$r(a,s)\ge 0$. (Here one can get even better result.)

Thus, for $1\ge a\ge {1\over 4}$ and $0\le t\le s\le 1$,
we have
$$
\split
&(f^{1\over 2})''_{tt}\ge\\
&f^{-{3\over 2}}\times {{11\over 64}a^2+{1\over 4}a\over (t+1)^4}
\endsplit
\tag{4.1.47}
$$
We have $f\le (1+(5/4)a)/(t+1)$. It follows
$$
(f^{1\over 2})''_{tt}\ge (1+{5\over 4}a)^{-{3\over 2}}\times
{{11\over 64}a^2+{1\over 4}a\over (t+1)^{5\over 2}}
\tag{4.1.48}
$$
if $1\ge a\ge {1\over 4}$ and $1\ge s\ge t\ge 0$.

We have $(t-s)^2/(t+1)=t+1-2(s+1)+(s+1)^2/(t+1)$, it follows
$$
\left({(t-s)^2\over (t+1)(s+1)}\right)''_{tt}=
{2(s+1)^2\over (s+1)(t+1)^3}\le
{8\over (s+1)(t+1)^3}.
\tag{4.1.49}
$$

Like for the proof of Lemma 4.1.6, from \thetag{4.1.23} and
\thetag{4.1.48}, it follows that for
$1\ge a_1,\,a_2\ge 0$, $1\ge a_3,\,a_4\ge {1\over 4}$,
$$
f_{p3}(a_1,\,a_2,\,a_3,\,a_4,\,s,\,t)''_{tt}\ge
$$
$$
\split
&(1+s)^{-1}(1+t)^{-3}
\left[2\left(\sqrt{a_1+(1-a_1)s^2}+\sqrt{a_2+(1-a_2)s^2}+
\right.\right.\\
&+\left.\sqrt{a_3+a_3(s-t)+{a_3(s-t)^2\over 4}+
{(1-a_3)(s+t)^2\over 4}}+\sqrt{a_4+(1-a_4)t^2}\right)\times\\
&\left((1+{5\over 4}a_3)^{-{3\over 2}}
({11\over 64}a_3^2+{1\over 4}a_3)+
\min{(a_4-{1\over 4}a_4^2,\ -{5\over 4}+6a_4-4a_4^2)}\right)-
\left.{(s+1)^2\over 2}\right]\ge
\endsplit
$$
$$
\split
&(1+s)^{-1}(1+t)^{-3}
\left[2\left(\sqrt{a_1}+\sqrt{a_2}+\sqrt{a_3}+
\sqrt{a_4}\right)\right.\times\\
&\left((1+{5\over 4}a_3)^{-{3\over 2}}
({11\over 64}a_3^2+{1\over 4}a_3)+
\min{(a_4-{1\over 4}a_4^2,\ -{5\over 4}+6a_4-4a_4^2)}\right)-
2\,].
\endsplit
\tag{4.1.50}
$$
Thus, $f_{p3}(a_1,\,a_2,\,a_3,\,a_4,\,s,\,t)''_{tt}\ge 0$
if $1\ge a_1,\,a_2\ge 0$, $1\ge a_3,\,a_4\ge {1\over 4}$ and
$$
\split
&\left(\sqrt{a_1}+\sqrt{a_2}+\sqrt{a_3}+
\sqrt{a_4}\right)\times\\
&\left((1+{5\over 4}a_3)^{-{3\over 2}}
({11\over 64}a_3^2+{1\over 4}a_3)+
\min{(a_4-{1\over 4}a_4^2,\ -{5\over 4}+6a_4-4a_4^2)}\right)\ \ge\
1.
\endsplit
\tag{4.1.51}
$$
For arbitrary $1\ge a_1,\,a_2,\,a_3,\,a_4\ge a\ge {1\over 4}$ it
is valid if
$$
\sqrt{a}
\left((1+{5\over 4}a)^{-{3\over 2}}
({11\over 64}a^2+{1\over 4}a)+
a-{1\over 4}a^2\right)
\ \ge {1\over 4}.
\tag{4.1.52}
$$
It is true if $a\ge 0.37646$.

It proves Lemma 4.1.7 for
$1\ge a_1,\,a_2,\,a_3,\,a_4\ge 0.37646$
when all three conditions \thetag{4.1.33}, \thetag{4.1.34} and
\thetag{4.1.35} are valid.
\enddemo

Using Theorems 4.1.4 and 4.1.5, we can divide all convex elliptic
polygons on the hyperbolic plane in three types. This subdivision is
very useful for fundamental polygons of finite volume of
reflection groups, and we shall use it in Sects. 4.3 and 5.

\proclaim{Theorem 4.1.8} Let $\M$ be an elliptic convex polygon on
a hyperbolic plane. Then $\M$ has one of types (I), (II) or (III)
(or the type (I), (II) or (III) of its narrow place) below:

{\bf Type (I):} There exist its four consecutive vertices
$A_0$, $A_1$, $A_2$, $A_3$ (where $A_0=A_3$ if $\M$ is
a triangle) with angles $\alpha_1=A_0A_1A_2$,
$\alpha_2=A_1A_2A_3\not={\pi\over 2}$ and such that
for orthogonal vectors $\delta_1$, $\delta_2$ and $\delta_3$ to
lines $(A_0A_1)$, $(A_1A_2)$ and $(A_2A_3)$ respectively directed
outwards of $\M$ and with
$\delta_1^2=\delta_2^2=\delta_3^2=-2$
one has $(\delta_1,\delta_2)=2\cos{\alpha_1}$,
$(\delta_2,\delta_3)=2\cos{\alpha_2}\not=0$ and
$$
\text{either\ \ } (\delta_1,\delta_3)\le 2\ \ \text{or}\ \
(\delta_1,\delta_3)
<4(\cos{\alpha_1\over 2}+\cos{\alpha_2\over 2})^2-2\le 14
\tag{4.1.53}
$$
The set $\{\delta_1,\,\delta_2,\,\delta_3\}$ generates the $3$-dimensional
hyperbolic vector space and has a connected Gram graph.
Any triangle or quadrangle $\M$ has the type (I).

{\bf Type (II):} There exist its five consecutive vertices
$A_0$, $A_1$, $A_2$, $A_3$, $A_4$
(where $A_0=A_4$ if $\M$ is a quadrangle) with angles
$\alpha_1=A_0A_1A_2={\pi\over 2}$, $\alpha_2=A_1A_2A_3={\pi\over 2}$,
$\alpha_3=A_2A_3A_4$ and orthogonal vectors
$\delta_1$, $\delta_2$, $\delta_3$ and $\delta_4$
to lines $(A_0A_1)$, $(A_1A_2)$, $(A_2A_3)$  and $(A_3A_4)$
respectively directed outwards of $\M$ and with
$\delta_1^2=\delta_2^2=\delta_3^2=\delta_4^2=-2$ such that
$(\delta_1,\,\delta_2)=(\delta_2,\,\delta_3)=0$,
$(\delta_3,\,\delta_4)=2cos{\alpha_3}$,
$$
(\delta_1,\,\delta_3)<6
\tag{4.1.54}
$$
and
$$
\split
&(\delta_1,\,\delta_4)\ <\
4\,\max{\left(\left(\sqrt{2}+\cos{\alpha_3\over 2}\right)^2,
{\left(2+\sqrt{2\cos^2{\alpha_3\over 2}+{1\over 4}}\,\right)^2-
{1\over 4}\over 2}\right)}-2\le\\
&\le \, 10+8\sqrt{2}=21.313708...\ .
\endsplit
\tag{4.1.55}
$$
Moreover, the set $\{\delta_1,\delta_3,\delta_4\}$ generates the
$3$-dimensional hyperbolic vector space and has
a connected Gram graph.

{\bf Type (III)}.
There exist its six consecutive vertices
$A_0$, $A_1$, $A_2$, $A_3$, $A_4$, $A_5$
(where $A_0=A_5$ if $\M$ is a pentagon) with right angles
$\alpha_1=A_0A_1A_2={\pi\over 2}$, $\alpha_2=A_1A_2A_3={\pi\over 2}$,
$\alpha_3=A_2A_3A_4={\pi\over 2}$, $\alpha_4=A_3A_4A_5={\pi\over 2}$
such that for orthogonal vectors $\delta_1$, $\delta_2$, $\delta_3$,
$\delta_4$ and $\delta_5$ to lines $(A_0A_1)$, $(A_1A_2)$, $(A_2A_3)$,
$(A_3A_4)$ and $(A_4A_5)$ respectively directed
outwards of $\M$ and with
$\delta_1^2=\delta_2^2=\delta_3^2=\delta_4^2=\delta_5^2=-2$
one has $(\delta_1,\,\delta_2)=(\delta_2,\,\delta_3)=
(\delta_3,\,\delta_4)=(\delta_4,\,\delta_5)=0$ and
$$
(\delta_1,\delta_3)<6,
\tag{4.1.56}
$$
$$
(\delta_3,\delta_5)<6,
\tag{4.1.57}
$$
and
$$
(\delta_1,\delta_5)<30.
\tag{4.1.58}
$$
Moreover, the set $\{\delta_1,\,\delta_3,\,\delta_5\}$ generates the
$3$-dimensional hyperbolic vector space and has a connected Gram graph.
\endproclaim

\demo{Proof} If $\M$ is a triangle, at least one angle of $\M$ is
not right. It follows that $\M$ has a sequence of consecutive
vertices of type (I). If $\M$ is a quadrangle,
we consider its narrow place of type (I) from Theorem 4.1.4. If both
angles $\alpha_1=\alpha_2={\pi\over 2}$, we replace $\delta_2$ by
$\delta_4$ with $\delta_4^2=-2$ which is orthogonal to the line
$(A_0A_3)$ and directed outwards of $\M$. Since
a quadrangle has at least one non-right angle (we can assume that it is
$\alpha_0=A_3A_0A_1$), then the sequence of vertices $A_2$, $A_3$, $A_0$,
$A_1$ with orthogonal vectors $\delta_3$, $\delta_4$, $\delta_1$
to sides $(A_2A_3)$, $(A_3A_0)$, $(A_0A_1)$ respectively
satisfies the condition (I) of Theorem 4.1.8.

We now suppose that $\M$ does not have a sequence of vertices of
type (I) (we always suppose that it is a sequence of consecutive
vertices). Then $\M$ has more than $4$ vertices.

By Theorem 4.1.5, $\M$ has a narrow place of type (II) or (III).

Suppose that $\M$ has the narrow place of type (II). We have
$\alpha_1=\alpha_2={\pi\over 2}$ since otherwise the sequence
$A_0$, $A_1$, $A_2$, $A_3$ gives a sequence of type (I)
(possibly after changing the numeration). We then get
estimates \thetag{4.1.54} and \thetag{4.1.55}. Obviously, the
set $\{\delta_1,\,\delta_2,\,\delta_3\}$ generates the 3-dimensional
hyperbolic vector space and has connected components
$\{\delta_1,\,\delta_3\}$ and $\{\delta_2\}$ since
$\alpha_1=\alpha_2={\pi\over 2}$. If the set
$\{\delta_1,\,\delta_3,\,\delta_4\}$ is not connected, it then
follows that $\delta_4$ is orthogonal to the set
$\{\delta_1,\,\delta_3\}$ and $\delta_4=\lambda\delta_2$
which is impossible since the lines
$(A_1A_2)$ and $(A_3A_4)$ are different. It follows that the
Gram graph of the set
$\{\delta_1,\,\delta_3,\,\delta_4\}$ is connected. This set generates
the 3-dimensional hyperbolic vector space because the orthogonal lines
$(A_0A_1)$, $(A_2A_3)$, $(A_3A_4)$ do not have a common point and are
not orthogonal to one line. Thus, we have proved that the $\M$ has
type (II) of Theorem 4.1.8.

Assume that $\M$ has a narrow place of type (III) of
Theorem 4.1.5. If the sequence
$A_0$, $A_1$, $A_2$, $A_3$, $A_4$, $A_5$
does not have a subsequence of the type (I) of Theorem 4.1.8, then
all angles $\alpha_1=\alpha_2=\alpha_3=\alpha_4={\pi\over 2}$ are
right. We then get estimates
\thetag{4.1.56}, \thetag{4.1.57} and \thetag{4.1.58}. The same
consideration as above, shows that the Gram graph of
the $\{\delta_1,\,\delta_3,\,\delta_5\}$ is connected and this
set generates the 3-dimensional hyperbolic vector space. Thus,
the $\M$ has type (III) of Theorem 4.1.8.

It finishes the proof.
\enddemo

\subhead
4.2. Narrow places of restricted parabolic convex polygons
on the hyperbolic plane
\endsubhead

Compare \cite{N9}, \cite{N11}.

We use notation of Sect. 1.1. Fix a point $O=\br_{++}c$ at infinity
of $\La=\La(\Phi)$. Thus, $c\in \Phi$, $(c,\,c)=0$ and $(c,\,V^+)>0$.

Recall that a {\it horosphere} $\E_O$ with the center $O$ is the set
of all lines in $\La$ containing $O$. The line
$l=O\br_{++}h\in \E_O$, $\br_{++}h\in \La$ is the set
$l=\{\br_{++}(tc+h)\ |\ t \in \br\ \text{and}\ (tc+h,\,tc+h)>0\}$.
Fix a constant $R>0$. Then there is a unique $\br_{++}h\in l$
such that $(h,\,c)=R$ and $(h,h)=1$. Given $l_1,\,l_2\in \E_O$, we
denote the corresponding $h$'s by $h_1$, $h_2$ and put
$$
\rho(l_1,l_2)=\sqrt{-(h_1-h_2)^2}.
\tag{4.2.1}
$$
Endowed with this distance, the horosphere $\E_O$ becomes an
affine Euclidean space. If one changes $R$, the distance
$\rho$ is multiplied by a constant. We denote
$$
\E_{O,R}=\{\br_{++}h\in \La\ |\ (h,c)=R\ \text{and}\
(h,h)=1\}.
\tag{4.2.2}
$$
The set $\E_{O,R}\cup\{O\}$ is a sphere in $\overline{\La}$,
which is tangent to $\La_\infty$ at $O$.
Moreover, the set $\E_{O,R}$ is orthogonal to every line $l\in \E_{O}$
at the point $\br_{++}h$, $h\in \E_{O,R}$ that corresponds to $l$.
The distance in $\La$ induces Euclidean distance in
$\E_{O,R}$ which is homothetic to the distance \thetag{4.2.1}.
The set $\E_{O,R}$ is identified with $\E_O$ and is also called a
{\it horosphere}.

Let $K\subset \E_O$. The set
$$
C_K=\bigcup_{\text{line\ }l\in K}{l}
\tag{4.2.3}
$$
is called the {\it cone with vertex $O$ and base $K$}.

A non-degenerate convex locally finite polyhedron $\M$ in $\La$ is
called {\it parabolic} (relative to the point $O\in \La_{\infty}$ if
1) and 2) below are valid:

1) $\M$ is finite at the point $O$, that is, the set
$\{\delta \in P(\M)\ |\ (c,\,\delta)=0\}$ is finite;

2) for every elliptic polyhedron ${\Cal N}\subset \E_O$ (that is
$\Cal N$ is the convex hull of a finite subset of $\E_O$), the
polyhedron $\M\cap C_{\Cal N}$ is elliptic.

A parabolic polyhedron $\M$ is called {\it restricted parabolic} if

3) the set
$$
r(\M)=\{(c,\delta/\sqrt{-\delta^2})\ |\ \delta \in P(\M)\}
\ \ \text{is finite.}
\tag{4.2.4}
$$
Geometrically this means that all hyperplanes $\Ha_\delta$,
$\delta\in P(\M)$, of faces of $\M$ are tangent to a finite set of
horospheres with the center $O$.

We remark that if $O\notin \M$ for a parabolic polyhedron $\M$,
then $\M$ is elliptic. Thus, it is only interesting to consider
parabolic polyhedra relative to a point
$O\in \La_\infty\cap \M$. Moreover, if $\M$ is
parabolic and has a finite number of
faces (or it has a finite set $P(\M)$) then $\M$ is elliptic. Thus,
it is only interesting to consider {\it infinite} (i.e. having infinite
number of faces) parabolic polyhedra.

Theorems 4.1.4, 4.1.5 and 4.1.8 can be generalized
on restricted parabolic polygons with strong inequalities
replaced by non-strong ones.

\proclaim{Theorem 4.2.1 (about the narrow place of type (I))}
For any restricted parabolic convex polygon $\M$ on a
hyperbolic plane there exist its four consecutive vertices
$A_0$, $A_1$, $A_2$ and $A_3$ (where $A_0=A_3$ if $\M$ is
a triangle) such that
for orthogonal vectors $\delta_1$, $\delta_2$ and $\delta_3$ to
lines $(A_0A_1)$, $(A_1A_2)$ and $(A_2A_3)$ respectively directed
outwards of $\M$ and with
$\delta_1^2=\delta_2^2=\delta_3^2=-2$
one has $(\delta_1,\delta_2)=2\cos{\alpha_1}$,
$(\delta_2,\delta_3)=2\cos{\alpha_2}$ and
$$
\text{either\ \ } (\delta_1,\delta_3)\le 2\ \ \text{or}\ \
(\delta_1,\delta_3)
\le 4(\cos{\alpha_1\over 2}+\cos{\alpha_2\over 2})^2-2\le 14
\tag{4.2.5}
$$
where $\alpha_1=A_0A_1A_2$ and $\alpha_2=A_1A_2A_3$.

Moreover, the Gram graph of $\{\delta_1,\,\delta_2,\,\delta_3\}$
is not connected (i.e. this set
is union of two non-empty orthogonal subsets)
if and only if $\alpha_1=\alpha_2={\pi\over 2}$.
\endproclaim

\proclaim{Theorem 4.2.2 (about narrow places of types (II) and (III))}
For any restricted parabolic convex polygon $\M$ having more than $3$
vertices (i.e. it is different from a triangle) on a hyperbolic plane,
one of two possibilities (II) or (III) below is valid:

(II) There exist its five consecutive vertices
$A_0$, $A_1$, $A_2$ , $A_3$ and $A_4$
(where $A_0=A_4$ if $\M$ is a quadrangle) such that
for orthogonal vectors $\delta_1$, $\delta_2$, $\delta_3$ and $\delta_4$
to lines $(A_0A_1)$, $(A_1A_2)$, $(A_2A_3)$  and $(A_3A_4)$
respectively directed outwards of $\M$ and with
$\delta_1^2=\delta_2^2=\delta_3^2=\delta_4^2=-2$,
one has $(\delta_i,\delta_{i+1})=2\cos{\alpha_i}$, $i=1,2,3$,
$$
\text{either\ \ } (\delta_1,\delta_3)\le 2\ \ \text{or}\ \
(\delta_1,\delta_3)
\le 4(\cos{\alpha_1\over 2}+\cos{\alpha_2\over 2})^2-2\le 14,
\tag{4.2.6}
$$
$$
(\delta_1,\delta_4)\le
$$
$$
\split
&4\,\max_{0\le t\le 1}
{{\left(\sqrt{a_1+(1-a_1)t^2}+\sqrt{a_2+(1-a_2)t^2}+
\sqrt{a_3+a_3t+t^2/4}\right)^2-{t^2\over 4}\over 1+t}}-2=\\
&4\,\max{\left(\left(\cos{\alpha_1\over 2}+\cos{\alpha_2\over 2}+
\cos{\alpha_3\over 2}\right)^2,
{\left(2+\sqrt{2\cos^2{\alpha_3\over 2}+{1\over 4}}\,\right)^2-
{1\over 4}\over 2}\right)}-2\le 34
\endsplit
\tag{4.2.7}
$$
where $\alpha_i=A_{i-1}A_iA_{i+1}$, $i=1,\,\,2,\,3,$ and
$a_i=\cos^2{\alpha_i\over 2}$. Moreover,
the set $\{\delta_1,\,\delta_2,$ $\delta_3,\,\delta_4\}$ has
a connected Gram graph.

(III). There exist its six consecutive vertices
$A_0$, $A_1$, $A_2$, $A_3$, $A_4$ and $A_5$
(where $A_0=A_5$ if $\M$ is a pentagon) such that
for orthogonal vectors $\delta_1$, $\delta_2$, $\delta_3$,
$\delta_4$ and $\delta_5$
to lines $(A_0A_1)$, $(A_1A_2)$, $(A_2A_3)$, $(A_3A_4)$ and
$(A_4A_5)$
respectively directed
outwards of $\M$ and with
$\delta_1^2=\delta_2^2=\delta_3^2=\delta_4^2=\delta_5^2=-2$
one has $(\delta_i,\delta_{i+1})=2\cos{\alpha_i}$, $i=1,2,3,4$,
$$
\text{either\ \ } (\delta_1,\delta_3)\le 2\ \ \text{or}\ \
(\delta_1,\delta_3)
\le 4(\cos{\alpha_1\over 2}+\cos{\alpha_2\over 2})^2-2\le 14,
\tag{4.2.8}
$$
$$
\text{either\ \ } (\delta_3,\delta_5)\le 2\ \ \text{or}\ \
(\delta_3,\delta_5)
\le 4(\cos{\alpha_3\over 2}+\cos{\alpha_4\over 2})^2-2\le 14,
\tag{4.2.9}
$$
and
$$
(\delta_1,\delta_5)\le
$$
$$
\split
&4\,\max_{0\le t\le s\le 1}
{\left[\left(\left(\sqrt{a_1+(1-a_1)s^2}+\sqrt{a_2+(1-a_2)s^2}+
\right.\right.\right.}\\
&+\left.\sqrt{a_3+a_3(s-t)+{a_3(s-t)^2\over 4}+
{(1-a_3)(s+t)^2\over 4}}+\sqrt{a_4+(1-a_4)t^2}\,\right)^2-\\
&-\left.\left.{(s-t)^2\over 4}\right)/((1+s)(1+t))\right]-2=\\
&4\,\max\left[\left(\cos{\alpha_1\over 2}+\cos{\alpha_2\over 2}+
\cos{\alpha_3\over 2}+\cos{\alpha_4\over 2}\right)^2,\right.\\
&\left.{\left(2+\sqrt{2\cos^2{\alpha_3\over 2}+{1\over 4}}+
\cos{\alpha_4\over 2}\right)^2-{1\over 4}\over 2},\ 4\right]-2\ \le\  62
\endsplit
\tag{4.2.10}
$$
where $\alpha_i=A_{i-1}A_iA_{i+1}$, $i=1,\,2,\,3,\,4$, and
$a_i=\cos^2{\alpha_i\over 2}$. Moreover, the
set $\{\delta_1,\,\delta_2,\,\delta_3,\,\delta_4,$
$\delta_5\}$ has a
connected Gram graph.
\endproclaim

\proclaim{Theorem 4.2.3} Let $\M$ be a restricted parabolic
convex polygon on
a hyperbolic plane. Then $\M$ has one of types (I), (II) or (III)
(or the type (I), (II) or (III) of its narrow place) below:

{\bf Type (I):} There exist its four consecutive vertices
$A_0$, $A_1$, $A_2$, $A_3$ (where $A_0=A_3$ if $\M$ is
a triangle) with angles $\alpha_1=A_0A_1A_2$,
$\alpha_2=A_1A_2A_3\not={\pi\over 2}$ and such that
for orthogonal vectors $\delta_1$, $\delta_2$ and $\delta_3$ to
lines $(A_0A_1)$, $(A_1A_2)$ and $(A_2A_3)$ respectively directed
outwards of $\M$ and with
$\delta_1^2=\delta_2^2=\delta_3^2=-2$
one has $(\delta_1,\,\delta_2)=2\cos{\alpha_1}$,
$(\delta_2,\,\delta_3)=2\cos{\alpha_2}\not=0$ and
$$
\text{either\ \ } (\delta_1,\,\delta_3)\le 2\ \ \text{or}\ \
(\delta_1,\,\delta_3)
\le 4(\cos{\alpha_1\over 2}+\cos{\alpha_2\over 2})^2-2\le 14
\tag{4.2.11}
$$
The set $\{\delta_1,\,\delta_2,\,\delta_3\}$ generates the $3$-dimensional
hyperbolic vector space and has a connected Gram graph.
Any triangle or quadrangle $\M$ has the type (I).

{\bf Type (II):} There exist its five consecutive vertices
$A_0$, $A_1$, $A_2$, $A_3$, $A_4$
(where $A_0=A_4$ if $\M$ is a quadrangle) with angles
$\alpha_1=A_0A_1A_2={\pi\over 2}$, $\alpha_2=A_1A_2A_3={\pi\over 2}$,
$\alpha_3=A_2A_3A_4$ and orthogonal vectors
$\delta_1$, $\delta_2$, $\delta_3$ and $\delta_4$
to lines $(A_0A_1)$, $(A_1A_2)$, $(A_2A_3)$  and $(A_3A_4)$
respectively directed outwards of $\M$ and with
$\delta_1^2=\delta_2^2=\delta_3^2=\delta_4^2=-2$ such that
$(\delta_1,\,\delta_2)=(\delta_2,\,\delta_3)=0$,
$(\delta_3,\,\delta_4)=2cos{\alpha_3}$,
$$
(\delta_1,\,\delta_3)\le 6
\tag{4.2.12}
$$
and
$$
(\delta_1,\,\delta_4)\ \le
$$
$$
\split
4\,\max{\left(\left(\sqrt{2}+\cos{\alpha_3\over 2}\right)^2,
{\left(2+\sqrt{2\cos^2{\alpha_3\over 2}+{1\over 4}}\,\right)^2-
{1\over 4}\over 2}\right)}-2\le
 &10+8\sqrt{2}=\\
 &=21.313708...\ .
\endsplit
\tag{4.2.13}
$$
Moreover, the set $\{\delta_1,\delta_3,\delta_4\}$ generates the
$3$-dimensional hyperbolic vector space and has
a connected Gram graph.

{\bf Type (III)}.
There exist its six consecutive vertices
$A_0$, $A_1$, $A_2$, $A_3$, $A_4$, $A_5$
(where $A_0=A_5$ if $\M$ is a pentagon) with right angles
$\alpha_1=A_0A_1A_2={\pi\over 2}$, $\alpha_2=A_1A_2A_3={\pi\over 2}$,
$\alpha_3=A_2A_3A_4={\pi\over 2}$, $\alpha_4=A_3A_4A_5={\pi\over 2}$
such that for orthogonal vectors $\delta_1$, $\delta_2$, $\delta_3$,
$\delta_4$ and $\delta_5$ to lines $(A_0A_1)$, $(A_1A_2)$, $(A_2A_3)$,
$(A_3A_4)$ and $(A_4A_5)$ respectively directed
outwards of $\M$ and with
$\delta_1^2=\delta_2^2=\delta_3^2=\delta_4^2=\delta_5^2=-2$
one has $(\delta_1,\,\delta_2)=(\delta_2,\,\delta_3)=
(\delta_3,\,\delta_4)=(\delta_4,\,\delta_5)=0$ and
$$
(\delta_1,\,\delta_3)\le 6,
\tag{4.2.14}
$$
$$
(\delta_3,\,\delta_5)\le 6,
\tag{4.2.15}
$$
and
$$
(\delta_1,\,\delta_5)\le 30.
\tag{4.2.16}
$$
Moreover, the set $\{\delta_1,\,\delta_3,\,\delta_5\}$ generates the
$3$-dimensional hyperbolic vector space and has a connected Gram graph.
\endproclaim

\demo{Proofs of Theorems 4.2.1, 4.2.2, 4.2.3} If $\M$ is parabolic
relative to $O\notin \M$, then $\M$ is finite and we can use
Theorems 4.1.4, 4.1.5 and 4.1.8. Thus, we can suppose that
$O\in \M$ and $\M$ is infinite.

To prove Theorems 4.1.4, 4.1.5 and
4.1.8, we were taken a point $O$ inside of $\M$ and
used the formulae of Lemma 4.1.2. To prove Theorems 4.2.1, 4.2.2, 4.2.3,
we use the point $O$ at infinity of $\M$ such that $\M$ is restricted
parabolic relative to $O$. We use an analog of Lemma 4.1.2 which uses
the infinite point $O$.

Let $O=\br_{++}c$ where $c^2=0$.
For a line $(AB)$ with terminals $A$ and $B$ at infinity
and $\delta\in \Phi$ orthogonal to $(AB)$ with $\delta^2=-2$ we
introduce an `angle'
$$
\theta(\delta)={1\over (c,\,\delta)}
\tag{4.2.17}
$$
where one can replace $\delta$ by $-\delta$ according to the orientation
of the angle $AOB$: e. g. one should take $\delta$ such that
$\Ha^+_\delta$
contains  $O$ if $AOB$ is oriented correctly, and one should take
$-\delta$ if not.

The `angle' $\theta(\delta)$ really behaves
like an angle. For three points $A$, $B$ and $C$ at infinity and
three lines $(AB)$, $(BC)$ and $(AC)$ with the corresponding vectors
$\delta_1$, $\delta_2$ and $\delta_3$ orthogonal to lines $(AB)$,
$(BC)$ and $(AC)$ respectively and with
$\delta_1^2=\delta_2^2=\delta_3^2=-2$, one has
$$
\theta(\delta_1)+\theta(\delta_2)=\theta(\delta_3).
\tag{4.2.18}
$$
We leave an elementary proof of \thetag{4.2.18} to a reader.

We have the following analog of Lemma 4.1.2.

\proclaim{Lemma 4.2.4} Let $(AB)$ and $(CD)$ are two lines on a
hyperbolic plane with terminals $A$, $B$, $C$, $D$ at infinity,
and $O=\br_{++}c$ an infinite point on the hyperbolic plane which
does not belong to each line $(AB)$ and $(CD)$ and orientations of the
triangles $AOB$ and $COD$ coincide.
Let $\delta_1$ and $\delta_2$ are orthogonal vectors with
square $-2$ to lines $(AB)$ and $(CD)$ respectively such that
$O$ is contained in both half-planes $\Ha_{\delta_1}^+$ and
$\Ha_{\delta_2}^+$. Let $\delta_{12}$ be the orthogonal vector with square
$-2$ to the line $(BC)$.
Then
$$
(\delta_1,\delta_2)=
4{(\theta(\delta_1)+\theta(\delta_{12}))
(\theta(\delta_2)+\theta(\delta_{12}))
\over \theta(\delta_1)\theta(\delta_2)}-2.
$$

As a corollary, we get :

1) If lines $(AB)$ and $(CD)$ do not intersect each other, then
$$
2 \cosh{\rho }=
(\delta_1,\delta_2)=
4{(\theta(\delta_1)+\theta(\delta_{12}))
(\theta(\delta_2)+\theta(\delta_{12}))
\over \theta(\delta_1)\theta(\delta_2)}-2
$$
where $\rho$ is the distance between lines $(AB)$ and $(CD)$
(here and in what follows we normalize the curvature $\kappa=-1$).

2) If lines $(AB)$ and $(CD)$ define an angle $\alpha$
containing $O$, then
$$
2 \cos{\alpha}=
(\delta_1,\delta_2)=
4{(\theta(\delta_1)-\theta(\delta_{21}))(\theta(\delta_2)-
\theta(\delta_{21}))
\over \theta(\delta_1)\theta(\delta_2)}-2
$$
where $\theta(\delta_{21})=-\theta(\delta_{12})=\theta(-\delta_{12})$.
\endproclaim

\demo{Proof of Lemma 4.2.4}
One can prove it similarly to Lemma 4.1.2.
One can also prove it as follows. Take a finite point $O^\prime$
and move $O^\prime$ to the infinite point $O$. Lemma 4.2.4
is the limit of the Lemma 4.1.2 applied to $O^\prime$
when $O^\prime$ tends to $O$. We leave details to a reader.
\enddemo

Now the proof of Theorems 4.2.1, 4.2.2 and 4.2.3 is the same as proof of
the corresponding theorems 4.1.4, 4.1.5 and 4.1.8 if one uses
the `angles' \thetag{4.2.17} and  Lemma 4.2.4. One has more:
almost all inequalities become equalities, and the proof
is even simpler.
\enddemo

\subhead
4.3. Description of narrow places of fundamental polygons $\M$ of
reflection subgroups $W\subset W(S)$ of elliptic and parabolic type
where $\rk S=3$.
Application to reflective lattices of elliptic and parabolic type
\endsubhead

Let $S$ be a primitive hyperbolic lattice of $\rk S=3$ and
$W\subset W(S)$ its reflection subgroup  of elliptic or parabolic type
with a fundamental polygon $\M$ (see Sects. 1.3, 1.4).
Remind that it means that $W$ and $P(\M)_{\pr}$ have restricted arithmetic
type and $P(\M)_{\pr}$ has a generalized lattice Weyl vector $\rho$ with
$\rho^2\ge 0$. It is known (e. g. see \cite{N11}, \cite{N9}) that
the polygon $\M$ is elliptic if $\rho^2>0$, and it is parabolic relative
to $\br_{++}\rho$ if $\rho^2=0$. We remind that the lattice $S$ having
reflection subgroups $W\subset W(S)$ of elliptic or parabolic type
is called reflective of elliptic or parabolic type.
Here we apply Theorem 4.2.3 to describe narrow places of $\M$. Using
this description, we shall get a finite list of lattices such that
any elliptically or parabolically reflective lattice $S$ belongs to
the list. Let $a(S^\ast/S)$ be the
exponent of the discriminant group $S^\ast/S$: i. e. $a(S)$ is the
least natural $a$ such that $aS^\ast/S=0$. We shall get an
estimate of $a(S^\ast /S)$.

Let $K$ be a lattice. We denote by $K_0$ the primitive lattice defined by
$K$.
Thus, $K=K_0(\lambda)$ where $\lambda\in \bn$ and $K_0$ is primitive.
If $K$ is generated by elements with the Gram matrix $A$, then $K_0$
is generated by the same elements with the Gram matrix $A/\lambda$ where
$\lambda$ is the greatest common divisor of all elements of $A$.
We denote by $a(A)$ the exponent of a finite Abelian group $A$: i. e.
$a(A)$ is the least natural $a$ such that $aA=0$.

We have the following useful statement (compare the proof of
Theorem 1 in \cite{N5, Appendix}).

\proclaim{Proposition 4.3.1} Let $L$ be a primitive lattice and
$\alpha_1,\dots,\alpha_k$ are primitive
roots of $L$ which generate a sublattice
$G\subset L$ of a finite index. Let $G_0$ be the
primitive lattice defined by $G$. Then $a(L^\ast/L)\vert
8a(G_0^\ast/G_0)^2$.
\endproclaim

\demo{Proof} We have $G=G_0(\lambda)\subset L$ for some natural
$\lambda \in \bn$. Since $G$ is generated by roots $\alpha_i$ of $L$
and $\alpha_i^2\vert 2(\alpha_i,L)$, we get
$2G_0(\lambda) \subset 2L\subset \lambda G_0(\lambda)^\ast$.
Identifying (naturally) modules of the lattices $G_0$ and $G_0(\lambda)$,
we obviously have $\lambda G_0(\lambda)^\ast=G_0^\ast$.
It follows that
$$
L=M(\lambda)\ \text{where}\ G_0\subset M\subset {1\over 2}G_0^\ast.
\tag{4.3.1}
$$
Here $M$ is any intermediate module which is invariant with
respect to reflections in roots $\alpha_i$ defining the lattice $G_0$,
and the roots $\alpha_i$ should be primitive in $M$.
The $\lambda \in \bn$ is the smallest natural number such that
$M(\lambda)$ is a lattice (otherwise, the lattice $L$ is not primitive).
Using \thetag{4.3.1}, we get
$$
2G_0\subset M^\ast \subset G_0^\ast.
\tag{4.3.2}
$$
If $tM\subset M^\ast$, $t\in \bn$, then $M(t)$ is a lattice.
Really, for any $m_1,\,m_2\in M$ we have
$t(m_1,m_2)=(tm_1,m_2)\in \bz$ because $tm_1\in M^\ast$.
Using \thetag{4.3.1} and \thetag{4.3.2}, we get
$4aM\subset M^\ast$ where $a=a(G_0^\ast/G_0)$.
It follows that $\lambda\vert 4a$.
Identifying
modules of $M$ and $M(\lambda )$, we have
$M(\lambda)^\ast=(1/\lambda)M^\ast\subset (1/4a)G_0^\ast$. It follows that
the exponent of $M(\lambda)^\ast/M(\lambda)$ divides the
exponent of $(1/4a)G_0^\ast /2G_0$ which is equal to
$8a^2$. This finishes the proof.
\enddemo

Below we describe narrow places of the fundamental polygons $\M$.
According to Theorem 4.2.3, a narrow place of $\M$ is defined by
vectors $\delta_1,\dots ,\delta_k$ with $\delta_i^2=-2$ where
$k\le 5$. We denote by $r_i$ corresponding primitive roots of
$S$ such that
$$
\delta_i={2r_i\over \sqrt{-2r_i^2}}.
\tag{4.3.3}
$$
with the Gram matrix
$$
\Gamma=\left(\gamma_{ij}\right)=\left((\delta_i,\delta_j)\right).
\tag{4.3.4}
$$
Since $r_i$ are primitive roots of $S$, we have
$r_i^2\vert 2(S,r_i)$. It then follows that
$$
\alpha_{ij}=\gamma_{ij}^2=(\delta_i,\,\delta_j)^2=
{4(r_i,r_j)^2\over r_i^2r_j^2}\in \bz_{+}
\tag{4.3.5}
$$
are non-negative integers.

We want to describe all possible matrices
$$
B=\left((r_i,r_j)\right)_{\pr}
\tag{4.3.6}
$$
where for an integral matrix $T$ we denote by $T_{\pr}$ the corresponding
primitive integral matrix $T/t$ where $t=\text{g.c.d}(T)$ denote
the greatest common divisor of all elements of $T$.
We shall make it in three steps.

First we describe all possible symmetric
$(k\times k)$-matrices
$$
{\Cal A}=(\alpha_{ij})=\left((\delta_i,\,\delta_j)^2\right)=
\left({4(r_i,r_j)^2\over r_i^2r_j^2}\right).
\tag{4.3.7}
$$
The matrix $\Cal A$ has non-negative integral coefficients, all
its diagonal coefficients are equal to $4$. All cyclic products
$$
\alpha_{i_1i_2}\alpha_{i_2i_3}\cdots \alpha_{i_{r-1}i_r}\alpha_{{i_r}i_1}
\tag{4.3.8}
$$
are perfect squares.

In the second place, we describe all possible $k\times k$ symmetrizable
generalized Cartan matrices (see \cite{Kac} about generalized Cartan
matrices, but remember that we use the opposite sign)
$$
A=(a_{ij})=\left({2(r_i,r_j)\over -r_i^2}\right)
\tag{4.3.9}
$$
using relations
$$
a_{ii}=-2,\ a_{ij}\in \bz_+\ \text{if}\ i\not=j,\
a_{ij}a_{ji}=\alpha_{ij},\
a_{ij}=0\ \text{iff}\ a_{ji}=0,
\tag{4.3.10}
$$
and
$$
a_{i_1i_2}a_{i_2i_3}\cdots a_{i_{r-1}i_r}a_{i_ri_1}=
a_{i_1i_r}a_{i_ri_{r-1}} \cdots a_{i_3i_2}a_{i_2i_1}.
\tag{4.3.11}
$$

In the third place, we find a diagonal matrix
$$
\Lambda =\text{diag}(\lambda_1,\dots,\lambda_k)
\tag{4.3.12}
$$
with $\lambda_i\in \bn$ such that the matrix
$$
B^\prime= A \Lambda
\tag{4.3.13}
$$
is symmetric. This defines the matrix $\Lambda$ uniquely up to
multiplication by a scalar rational matrix (since the Gram graph of
$r_1,\dots,r_k$ is connected and the matrix $A$ is
indecomposable). The matrix
$\Lambda=s\,\text{diag}(-r_1^2,\dots,-r_k^2)$
where $s\in \bq_{++}$. Then we
calculate
$$
B={B^\prime\over \text{g.c.d}(B^\prime)}
\tag{4.3.14}
$$
which gives the matrix $\left((r_i,r_j)\right)_{\pr}$,
see \thetag{4.3.6}.
These procedure gives a finite set of possible matrices $B$.
See the corresponding general considerations in \cite{N5, Appendix})

For the lattice $G_0=[r_1,\dots , r_k]_{\pr}$ defined by $B$
we calculate the invariant $a(B)=a(G_0^\ast/G_0)$, the invariant
$a_1(B)$ which is the product of all different odd prime
divisors of $a(B)$, and the invariant
$a_2(B)$ which is the greatest odd prime
divisor of $a(B)$. By Proposition 4.3.1, we get an estimate
of the similar invariants $a(S)$, $a_1(S)$ and $a_2(S)$:
$$
a(S)\le 8a(B)^2,\ a_1(S)\le a_1(B),\ a_2(S)\le a_2(B).
\tag{4.3.15}
$$
Calculating the invariants $a(B)$, $a_1(B)$ and $a_2(B)$
for all possible matrices $B$ and taking their maximum,
we estimate the invariants $a(S)$, $a_1(S)$ and $a_2(S)$ for all
elliptically or parabolically reflective lattices $S$.

Below we describe this procedure for each type of the narrow place of
Theorem 4.2.3.

\subsubhead
4.3.1. Matrices $B$ of the narrow places of the type (I1)
\endsubsubhead
 It is a particular case of the type (I) of Theorem 4.2.3 when
additionally the angle $\alpha_1\not=\pi/2$. For this case
$k=3$ and $\delta_1,\,\delta_2,\,\delta_3$ give a bases of
the $3$-dimensional hyperbolic vector space. We get that the
matrix $\Cal A$ is a symmetric matrix
$$
{\Cal A}=\pmatrix
4 & \alpha_{12} &\alpha_{13}\\
\alpha_{21}&4&\alpha_{23}\\
\alpha_{31}&\alpha_{32}&4
\endpmatrix
\tag{4.3.1.1}
$$
with integral non-negative coefficients
where from the condition of narrow places we have
(after changing numeration if necessary) that
$$
\split
&1\le \alpha_{12}=\alpha_{21}\le 4,\ \
\alpha_{12}\le \alpha_{23}=\alpha_{32}\le 4,\\
&\alpha_{23}\le \alpha_{13}=\alpha_{31}\le
\left[\left(\sqrt{2+\sqrt{\alpha_{12}}}+
\sqrt{2+\sqrt{\alpha_{23}}}\,\right)^2-
2\right]^2
\endsplit
\tag{4.3.1.2}
$$
(here we also use that $2\cos{\alpha\over 2}=\sqrt{2+2\cos{\alpha}}$),
and
$$
d=det(\Gamma )=-8+2\sqrt{\alpha_{12}\alpha_{23}\alpha_{31}}+
2\alpha_{12}+2\alpha_{23}+2\alpha_{13}>0
\tag{4.3.1.3}
$$
where
$$
\alpha_{12}\alpha_{23}\alpha_{31}\ \ \text{is a perfect square.}
\tag{4.3.1.4}
$$
It is easy to enumerate the finite set of all matrices $\Cal A$
satisfying \thetag{4.3.1.2}, \thetag{4.3.1.3} and \thetag{4.3.1.4}.
For each $\Cal A$ we then find all symmetrizable generalized
Cartan matrices
$$
A=
\pmatrix
-2     & a_{12}   &a_{13}\\
a_{21} & -2       & a_{23}\\
a_{31} & a_{32}   & -2
\endpmatrix
\tag{4.3.1.5}
$$
using the relations
$$
a_{12}a_{21}=\alpha_{12},\ a_{23}a_{32}=\alpha_{23},\
a_{13}a_{31}=\alpha_{13},\
a_{12}a_{23}a_{31}=a_{13}a_{32}a_{21}.
\tag{4.3.1.6}
$$
The diagonal matrix
$$
\Lambda=\text{diag}(a_{13}a_{32},\ a_{23}a_{31},\ a_{31}a_{32}).
\tag{4.3.1.7}
$$
Finally, we get that
$$
B=(A\Lambda)_{\pr}.
\tag{4.3.1.8}
$$

In Appendix, we give the Program 4: fund11.gen which uses this algorithm
to enumerate all the matrices $B$. For each of them it calculates
the invariants $a(B)$, $a_1(B)$, $a_2(B)$ and finds the number
$nI1$ of all the matrices $B$, and the numbers
$$
aI1=\max_{B}{a(B)},\ aI1_1=\max_{B}{a_1(B)},\ aI1_2=\max_B{a_2(B)}.
\tag{4.3.1.9}
$$
Calculation using this program gives
$$
nI1=272,\ aI1=3528,\ aI1_1=543,\ aI1_2=181.
\tag{4.3.1.10}
$$

\subsubhead
4.3.2. Matrices $B$ of the narrow places of the type (I0)
\endsubsubhead
 It is a particular case of Type (I) of Theorem 4.2.3 when
additionally the angle $\alpha_1=\pi/2$. For this case
$k=3$ and $\delta_1,\,\delta_2,\,\delta_3$ give a bases of
the $3$-dimensional hyperbolic vector space. We get that the
matrix $\Cal A$ is a symmetric matrix
$$
{\Cal A}=\pmatrix
4 & 0 &\alpha_{13}\\
0 &4&\alpha_{23}\\
\alpha_{31}&\alpha_{32}&4
\endpmatrix
\tag{4.3.2.1}
$$
with integral non-negative coefficients
where from the condition of narrow places we have
(after changing numeration if necessary) that
$$
1\le \alpha_{23}=\alpha_{32}\le 4,\ \
\alpha_{23}\le \alpha_{13}=\alpha_{31}\le
\left[\left(\sqrt{2}+\sqrt{2+\sqrt{\alpha_{23}}}\right)^2-2\right]^2,
\tag{4.3.2.2}
$$
and
$$
d=det(\Gamma )=-8+2\alpha_{23}+2\alpha_{13}>0.
\tag{4.3.2.3}
$$
For each $\Cal A$ we then
find all symmetrizable generalized Cartan matrices
$$
A=
\pmatrix
-2     & 0   &a_{13}\\
0      & -2       & a_{23}\\
a_{31} & a_{32}   & -2
\endpmatrix
\tag{4.3.2.4}
$$
using the relations
$$
a_{23}a_{32}=\alpha_{23},\ a_{13}a_{31}=\alpha_{13}.
\tag{4.3.2.5}
$$
The diagonal matrix
$$
\Lambda=\text{diag}(a_{13}a_{32},\ a_{23}a_{31},\ a_{31}a_{32}).
\tag{4.3.2.6}
$$
Finally, we get that
$$
B=(A\Lambda)_{\pr}.
\tag{4.3.2.7}
$$
In Appendix, we give the Program 5: fund10.gen which uses this algorithm
to enumerate all the matrices $B$. For each of them it calculates
the invariants $a(B)$, $a_1(B)$, $a_2(B)$ and finds the number
$nI0$ of all the matrices $B$, and the numbers
$$
aI0=\max_{B}{a(B)},\ aI0_1=\max_{B}{a_1(B)},\ aI0_2=\max_B{a_2(B)}.
\tag{4.3.2.8}
$$
Calculation using this program gives
$$
nI0=2998,\ aI0=69192,\ aI0_1=10209,\ aI0_2=89.
\tag{4.3.2.9}
$$

\subsubhead
4.3.3. Matrices $B$ of the narrow places of the type (II1)
\endsubsubhead
 It is a particular case of Type (II) of Theorem 4.2.3 when
additionally the angle $\alpha_3\not=\pi/2$. For this case
$k=4$ and $\delta_1,\,\delta_2,\,\delta_3,\,\delta_4$
generate the $3$-dimensional hyperbolic vector space and any three of them
give a bases of the space. We get that the
matrix $\Cal A$ is a symmetric matrix
$$
{\Cal A}=\pmatrix
4          & 0         &\alpha_{13}&\alpha_{14}\\
0          & 4         &   0       & \alpha_{24}\\
\alpha_{31}& 0         &   4       & \alpha_{34}\\
\alpha_{41}&\alpha_{42}&\alpha_{43}&  4
\endpmatrix
\tag{4.3.3.1}
$$
with integral non-negative coefficients where
$$
\split
&1\le \alpha_{34}=\alpha_{43}\le 4,\ \
4< \alpha_{13}=\alpha_{31}\le 36,\\
&0\le \alpha_{14}=\alpha_{41}\le
\left[4\,\max{
\left(\left(\sqrt{2}\hskip-1pt+\hskip-1pt
\sqrt{{\sqrt{\alpha_{34}}\over 4}\hskip-1pt+\hskip-1pt
{1\over 2}}\,\right)^2,
{\left(2\hskip-1pt+\hskip-1pt\sqrt{{\sqrt{\alpha_{34}}\over 2}
\hskip-1pt+\hskip-1pt
{5\over 4}}\,\right)^2
\hskip-4pt - \hskip-3pt {1\over 4}\over 2}\right)}-2\right]^2\hskip-6pt,\\
&\left[\left(\sqrt{2}+\sqrt{2+\sqrt{\alpha_{34}}}\right)^2
\hskip-2pt-\hskip-2pt 2\right]^2
<\alpha_{24}=\alpha_{42}.
\endsplit
\tag{4.3.3.2}
$$
Here we use the conditions (II) of Theorem 4.2.3 and add some
inequalities to avoid repeating of cases we have considered in
Sect. 4.3.3. We also have
$$
det(\Gamma )=4(4-\sqrt{\alpha_{13}\alpha_{34}\alpha_{14}}-
\alpha_{14}-\alpha_{34}-\alpha_{13})+\alpha_{24}\alpha_{13}-4\alpha_{24}=0.
\tag{4.3.3.3}
$$
where
$$
\alpha_{13}\alpha_{34}\alpha_{14}\ \text{is a perfect square}.
\tag{4.3.3.4}
$$
The condition \thetag{4.3.3.3} is equivalent
$$
\alpha_{24}={4(\sqrt{\alpha_{13}\alpha_{34}\alpha_{14}}+
\alpha_{14}+\alpha_{34}+\alpha_{13}-4)
\over
\alpha_{13}-4}.
\tag{4.3.3.5}
$$
We can easily enumerate the finite set of all possible
matrices $\Cal A$ satisfying these conditions. For
each $\Cal A$ we find all symmetrizable generalized Cartan
matrices
$$
A=
\pmatrix
-2      &  0      &a_{13}   & a_{14}\\
0       & -2      &   0     & a_{24}\\
a_{31}  &  0      &   -2    & a_{34}\\
a_{41}  &  a_{42} &a_{43}   &  -2
\endpmatrix
\tag{4.3.3.6}
$$
using relations
$$
\split
&a_{34}a_{43}=\alpha_{34},\ a_{13}a_{31}=\alpha_{13},\
a_{24}a_{42}=\alpha_{24},\ a_{13}a_{34}a_{41}=
\sqrt{\alpha_{13}\alpha_{34}\alpha_{14}},\\
&a_{14}a_{41}=\alpha_{14},\ a_{14}=a_{41}=0\
\text{if}\ \alpha_{14}=0.
\endsplit
\tag{4.3.3.7}
$$
The diagonal matrix
$$
\Lambda=\text{diag}(a_{13}a_{34}a_{42},\ a_{31}a_{43}a_{24},\
a_{31}a_{34}a_{42},\ a_{31}a_{43}a_{42}).
\tag{4.3.3.8}
$$
Finally, we get that
$$
B=(A\Lambda)_{\pr}.
\tag{4.3.3.9}
$$
In Appendix, we give the Program 6: fund21.gen which uses this algorithm
to enumerate all the matrices $B$. For each of them it calculates
the invariants $a(B)$, $a_1(B)$, $a_2(B)$ and finds the number
$nII1$ of all the matrices $B$, and the numbers
$$
aII1=\max_{B}{a(B)},\ aII1_1=\max_{B}{a_1(B)},\ aII1_2=\max_B{a_2(B)}.
\tag{4.3.3.10}
$$
Calculation using this program gives
$$
nII1=9818,\ aII1=47432,\ aII1_1=10965,\ aII1_2=487.
\tag{4.3.3.11}
$$

\subsubhead
4.3.4. Matrices $B$ of the narrow places of the type (II0)
\endsubsubhead
 It is a particular case of Type (II) of Theorem 4.2.3 when
additionally the angle $\alpha_3=\pi/2$. For this case
$k=4$ and $\delta_1,\,\delta_2,\,\delta_3,\,\delta_4$
generate the $3$-dimensional hyperbolic vector space and any three of them
give a bases of the space. We get that the
matrix $\Cal A$ is a symmetric matrix
$$
{\Cal A}=\pmatrix
4          & 0         &\alpha_{13}&\alpha_{14}\\
0          & 4         &   0       & \alpha_{24}\\
\alpha_{31}& 0         &   4       &  0\\
\alpha_{41}&\alpha_{42}&   0       &  4
\endpmatrix
\tag{4.3.4.1}
$$
with integral non-negative coefficients where
$$
4< \alpha_{13}=\alpha_{31}\le 36,\ \
0< \alpha_{14}=\alpha_{41}\le (8+4\sqrt{5})^2=287.108350...\,,\
\alpha_{13}\le \alpha_{24}.
\tag{4.3.4.2}
$$
($8+4\sqrt{5}=16.94427190...$) Here we use the conditions (II)
of Theorem 4.2.3 and that
$\M$ is a fundamental polygon having at least $4$ sides.
We also have
$$
det(\Gamma )=4(4-\alpha_{14}-\alpha_{13})+\alpha_{24}\alpha_{13}-
4\alpha_{24}=0
\tag{4.3.4.3}
$$
which is equivalent
$$
(\alpha_{13}-4)(\alpha_{24}-4)=4\alpha_{14}.
\tag{4.3.4.4}
$$
It is easy to enumerate the finite set of all possible
matrices $\Cal A$ satisfying these conditions. For
each $\Cal A$ we find all symmetrizable generalized Cartan
matrices
$$
A=
\pmatrix
-2      &  0      &a_{13}   & a_{14}\\
0       & -2      &   0     & a_{24}\\
a_{31}  &  0      &   -2    & 0     \\
a_{41}  &  a_{42} &   0     & -2
\endpmatrix
\tag{4.3.4.5}
$$
using relations
$$
a_{13}a_{31}=\alpha_{13},\
a_{24}a_{42}=\alpha_{24},\
a_{14}a_{41}=\alpha_{14},\
\tag{4.3.4.6}
$$
The diagonal matrix
$$
\Lambda=\text{diag}(a_{13}a_{14}a_{42},\ a_{13}a_{41}a_{24},\
a_{31}a_{14}a_{42},\ a_{13}a_{42}a_{41}).
\tag{4.3.4.7}
$$
Finally, we get that
$$
B=(A\Lambda)_{\pr}.
\tag{4.3.4.8}
$$
In Appendix, we give the Program 7: fund20.gen which uses this algorithm
to enumerate all the matrices $B$. For each of them it calculates
the invariants $a(B)$, $a_1(B)$, $a_2(B)$ and finds the number
$nII0$ of all the matrices $B$, and the numbers
$$
aII0=\max_{B}{a(B)},\ aII0_1=\max_{B}{a_1(B)},\ aII0_2=\max_B{a_2(B)}.
\tag{4.3.4.9}
$$
Calculation using this program gives
$$
nII0=376208,\ aII0=995316,\ aII0_1=238569,\ aII0_2=283.
\tag{4.3.4.10}
$$

\subsubhead
4.3.5. Matrices $B$ of the narrow places of the type (III)
\endsubsubhead
For this case
$k=5$ and $\delta_1,\,\delta_2,\,\delta_3,\,\delta_4,\,\delta_5$
generate the $3$-dimensional hyperbolic vector space, and any three of
them
give a bases of the space. The matrix $\Cal A$ is a symmetric matrix
$$
{\Cal A}=\pmatrix
4          & 0         &\alpha_{13}&\alpha_{14} &\alpha_{15}\\
0          & 4         &   0       & \alpha_{24}&\alpha_{25}\\
\alpha_{31}& 0         &   4       &  0         &\alpha_{35}\\
\alpha_{41}&\alpha_{42}&   0       &  4         & 0         \\
\alpha_{51}&\alpha_{52}&\alpha_{53}&  0         & 4
\endpmatrix
\tag{4.3.5.1}
$$
with integral non-negative coefficients where
$$
\split
&4< \alpha_{13}=\alpha_{31}\le 36,\ \
\alpha_{31}\le \alpha_{35}=\alpha_{53}\le 36,\ \
0\le \alpha_{15}=\alpha_{51}\le 30^2=900,\\
&287.108350<\alpha_{14}=\alpha_{41},\ \
287.108350<\alpha_{25}=\alpha_{52},
\endsplit
\tag{4.3.5.2}
$$
Here we use the conditions (III) of Theorem 4.2.3 and that
$\M$ is a fundamental polygon having at least $5$ sides. The
last two inequalities were added to avoid repetition with the
previous case. The Gram matrix
$\Gamma=\Gamma(\delta_1,\,\delta_2,\,\delta_3,\,\delta_4,\,\delta_5)$
has the rank $3$. Coefficients of the matrix $\Cal A$
are defined by the coefficients $\alpha_{13}$, $\alpha_{15}$,
$\alpha_{35}$ defining the Gram matrix
$\Gamma(\delta_1,\,\delta_3,\,\delta_5)$. Let
$$
d=4(\alpha_{13}+\alpha_{35}+\alpha_{15}+
\sqrt{\alpha_{13}\alpha_{35}\alpha_{15}}-4),
\tag{4.3.5.3}
$$
one can see that
$d/2=det(\Gamma(\delta_1,\,\delta_3,\,\delta_5))$. We have
$$
\alpha_{14}={d\over \alpha_{35}-4},\ \
\alpha_{25}={d\over \alpha_{13}-4}
\tag{4.3.5.4}
$$
and
$$
\alpha_{24}={4(\alpha_{13}\alpha_{35}+
4\sqrt{\alpha_{13}\alpha_{35}\alpha_{15}}+4\alpha_{15})
\over
(\alpha_{13}-4)(\alpha_{35}-4)}
\tag{4.3.5.5}
$$
Here \thetag{4.3.5.4} follows from
$\det(\Gamma (\delta_1,\,\delta_3,\,\delta_5,\,\delta_2))=
\det(\Gamma (\delta_1,\,\delta_3,\,\delta_5,\,\delta_4))=0$.
To get \thetag{4.3.5.5}, one should remark that
$\delta_2=\gamma_{25}\delta_5^\ast$ and
$\delta_4=\gamma_{14}\delta_1^\ast$.
It follows that $\gamma_{24}=\gamma_{25}\gamma_{14}(g^{-1})_{13}$
where $g=\Gamma(\delta_1,\,\delta_3,\,\delta_5)$. All together,
\thetag{4.3.5.4} and \thetag{4.3.5.5} are equivalent to
$\rk(\Gamma(\delta_1,\,\delta_2,\,\delta_3,\,\delta_4,\,\delta_5))=3$,
or that $\delta_1,\,\delta_2,\,\delta_3,\,\delta_4,\,\delta_5$ generate
a $3$-dimensional hyperbolic form. Moreover, cyclic products
$$
\alpha_{13}\alpha_{35}\alpha_{52}\alpha_{24}\alpha_{41},\ \
\alpha_{13}\alpha_{35}\alpha_{51}\ \ \text{are perfect squares}.
\tag{4.3.5.6}
$$
It is easy to enumerate the finite set of all possible
matrices $\Cal A$ satisfying these conditions. For
each $\Cal A$ we find all symmetrizable generalized Cartan
matrices
$$
A=\pmatrix
-2           & 0         &a_{13}    & a_{14} & a_{15}\\
 0           & -2        &   0      & a_{24} & a_{25}\\
 a_{31}      & 0         &   -2     &  0     & a_{35}\\
 a_{41}      & a_{42}    &   0      &  -2    & 0     \\
 a_{51}      & a_{52}    &   a_{53} &  0     & -2
\endpmatrix
\tag{4.3.5.7}
$$
using relations
$$
\split
&a_{13}a_{31}=\alpha_{13},\ a_{35}a_{53}=\alpha_{35},\
a_{14}a_{41}=\alpha_{41},\ a_{25}a_{52}=\alpha_{25},\
a_{15}a_{51}=\alpha_{15},\\
&a_{15}=a_{51}=0\ \text{if}\ \alpha_{15}=0,\
a_{24}a_{42}=\alpha_{24},\ a_{25}a_{52}=\alpha_{25},\
a_{35}a_{53}=\alpha_{35},\\
&a_{13}a_{35}a_{52}a_{24}a_{41}=
\sqrt{\alpha_{13}\alpha_{35}\alpha_{52}\alpha_{24}\alpha_{41}},\
a_{13}a_{35}a_{51}=\sqrt{\alpha_{13}\alpha_{35}\alpha_{51}},
\endsplit
\tag{4.3.5.8}
$$
The diagonal matrix
$$
\Lambda=\text{diag}
(a_{14}a_{13}a_{42}a_{25},\
a_{41}a_{13}a_{24}a_{25},\
a_{14}a_{31}a_{42}a_{25},\
a_{41}a_{13}a_{42}a_{25},\
a_{41}a_{13}a_{24}a_{52})
\tag{4.3.5.9}
$$
Finally, we get that
$$
B=(A\Lambda)_{\pr}.
\tag{4.3.5.10}
$$
In Appendix, we give the Program 8: fund30.gen which uses this algorithm
to enumerate all the matrices $B$. For each of them it calculates
the invariants $a(B)$, $a_1(B)$, $a_2(B)$ and finds the number
$nIII$ of all the matrices $B$, and the numbers
$$
aIII=\max_{B}{a(B)},\ aIII_1=\max_{B}{a_1(B)},\ aIII_2=\max_B{a_2(B)}.
\tag{4.3.5.11}
$$
Calculation using this program gives
$$
nIII=200539,\ aIII=324900,\ aIII_1=26565,\ aIII_2=907.
\tag{4.3.5.12}
$$

\subsubhead
4.3.6. The global estimate of invariants of primitive reflective
hyperbolic lattices of the rank $3$ having elliptic or parabolic type
\endsubsubhead

For a lattice $L$ we denote by $a(L^\ast/L)$ the exponent of
the discriminant group $L^\ast/L$, we denote by $a_1(L^\ast/L)$
the product of all different odd prime divisors of $a(L^\ast/L)$, and
we denote by $a_2(L^\ast/L)$ the greatest prime divisor of
$a(L^\ast/L)$.

Using Proposition 4.3.1 and calculations \thetag{4.3.1.10},
\thetag{4.3.2.9}, \thetag{4.3.3.11}, \linebreak \thetag{4.3.4.10},
\thetag{4.3.5.12} we get

\proclaim{Theorem 4.3.6.1} For any primitive reflective hyperbolic
lattice $S$ of $\rk S=3$ having elliptic or parabolic type we
have estimates:
$$
a(S^\ast/S)\le 8(995316)^2,\ \
a_1(S^\ast/S)\le 238569,\ \
a_2(S^\ast/S)\le 907.
$$
\endproclaim

Since $\det (S)\le a(S^\ast/S)^2$ and number of lattices with the
fixed rank and determinant is finite (e. g. see \cite{C}),
Theorem 4.3.6.1 gives a finite list of lattices which contains
all the reflective lattices $S$.

The estimates of Theorem 4.3.6.1 are very preliminary,
and we shall significantly improve them below.

\head
5. Classification of reflective hyperbolic lattices of the rank 3 and
of elliptic or parabolic type: proofs
\endhead

Now we are ready to prove classification results of Sect. 2.

\subhead
5.1. Proof of Basic Theorem 2.3.2.1
\endsubhead

We first prove

\proclaim{Theorem 5.1.1} Any elliptically or parabolically
reflective main hyperbolic lattice $S$ of
rank $3$ and with square-free determinant belongs to the list of
Table 3 containing all main hyperbolic lattices
$S$ of the rank $3$ with square-free $d=\det(S)\le 100000$ and
$h=hnr(S)\le 1$.
\endproclaim

\demo{Proof}
Let $S$ be an elliptically or parabolically reflective main
hyperbolic lattice of the rank $3$ and with square-free
determinant $d=\det (S)$. Let $\M$ be a fundamental polygon of
$W(S)$.

First we will show that $h(S)\le 1$ (it gives the proof of
Lemma 2.3.1.2).
If $S$ is elliptically reflective, then $\M$ is an elliptic
(finite and of finite volume) polygon. This polygon has $\le 1$
central symmetries. Really, otherwise, a composition of two
different central symmetries gives an automorphism of infinite order
of $\M$ which is impossible. It follows that $h(S)\le 1$.
Assume that $S$ is parabolically reflective and $\rho$ is
a generalized lattice Weyl vector for $\M$.
Let $A(\M)\subset O^+(S)$ be the group of symmetries of $\M$. By
definition of a generalized lattice Weyl vector $\rho \in S$,
then $\rho$ is preserved by a subgroup $A\subset A(\M)$ of finite
index, $\rho^2=0$ and $\rho\not=0$.
For any $\phi \in A(\M)$, the element $\phi (\rho)$ also has all these
properties and is then a generalized
lattice Weyl vector for $\M$. If $\phi (\rho)\not=\rho$, the group
$A(\M)$ has a subgroup of finite index which is trivial
on the hyperbolic sublattice $\bz\rho+\bz\phi(\rho)\subset S$.
It then follows that $A(\M)$ is finite and $S$ is
elliptically reflective. We get a contradiction.
Thus $\phi (\rho)=\rho$ for any $\phi \in A(\M)$.
If $u\in A(\M)$ is a central symmetry, then the fixed part
$S^u=\{x\in S\ |\ u(x)=x\}$ of $u$
is negative definite and does not have non-zero elements $\rho$ with
$\rho^2=0$. This shows that $h=h(S)=0$.

Now we apply to $S$ and $\M$ results of Sect. 4.3.
By Proposition 4.3.1 and \thetag{4.3.1.10}, \thetag{4.3.2.9},
\thetag{4.3.3.11}, \thetag{4.3.5.12}, we have
$d\le 100000$ if $\M$ has a narrow place of types
$I1$, $I0$, $II1$ or $III$. Really, for all this cases
the invariant $a_1< 50000$ and then $d\le 2a_1(S)\le 2a_1< 100000$.
Thus, Theorem 5.1.1 is valid for these lattices $S$.

Only if $\M$ has a narrow place of type $II0$, our
estimate \thetag{4.3.4.10} of $a_1$ is not good enough:
$$
a_1(S)\le aII0_1=238569,\ \ a_2(S)\le 238.
\tag{5.1.1}
$$
Let us consider this case.
In Sect. 4.3.4, for all primitive elliptically or parabolically
reflective hyperbolic lattices $S$ or rank $3$
we considered the primitive Gram matrices
$B=\Gamma_{\pr}$ where
$\Gamma=\Gamma(\alpha_1,\,\alpha_2,\,\alpha_3,\,\alpha_4)$ is
the Gram matrix of primitive roots of $S$ corresponding to
a narrow place of type $II0$ of $\M$. We now consider them only
for main lattices $S$ with square-free $d=\det(S)$.
Then we should add some additional
conditions to Sect. 4.3.4. We introduce these
conditions below.

Let $\beta \in S$ be a primitive root and $K=\delta^\perp$.
Then
$$
\text{either\ } S=\bz\beta\oplus K\ \text{or\ }
S=[\beta,\ K, (\beta+k)/2]\ \text{where\ } k\in K.
\tag{5.1.2}
$$
It follows (by a simple consideration over $\bz_2$)
that
$$
\beta^2\ \text{is square-free}.
\tag{5.1.3}
$$
If $\beta_1,\,\beta_2\in S$ are two primitive orthogonal roots
and $\bz h=[\beta_1,\,\beta_2]^\perp$, then for any odd $p$
$$
S\otimes \bz_p=\bz_p\beta_1\oplus \bz_p\beta_2\oplus \bz_p h\
\text{and\ }h^2\ \text{is square-free.}
\tag{5.1.4}
$$
It follows
$$
\text{g.c.d.}(\beta_1^2,\beta_2^2)\le 2,\ \
\text{g.c.d.}(\beta_1^2,h^2)\le 2,\ \
\text{g.c.d.}(\beta_2^2,h^2)\le 2,
\tag{5.1.5}
$$
and
$$
d=\beta_1^2\beta_2^2h^2/ 2^t,\ \text{where\ }t\equiv 0\mod 2.
\tag{5.1.6}
$$
It defines the $t$ uniquely because $d$ is square-free. From
\thetag{5.1.6}, for any odd $p\vert d$, we also have
$$
(-1)^{\eta_p}=
\cases
\left({\beta_1^2/p\over p}\right) &\text{if $p\vert\beta_1^2$,}\\
\left({\beta_2^2/p\over p}\right) &\text{if $p\vert\beta_2^2$,}\\
\left({h^2/p\over p}\right)=
\left({d\beta_1^2\beta_2^2/p\over p}\right)
&\text{if $p\vert h^2$,}
\endcases
\tag{5.1.7}
$$
where $\eta$ is the invariant of $S$. The Gram matrix
$\Gamma=\Gamma(\alpha_1,\,\alpha_2,\,\alpha_3,\,\alpha_4)$
has three pairs of orthogonal roots:
$(\alpha_1,\,\alpha_2)=(\alpha_2,\,\alpha_3)=
(\alpha_3,\alpha_4)=0$. It then follows (it is sufficient to have only
one pair ) that
$$
\text{either\ } \Gamma=B\
\text{or\ }
\Gamma=2B.
\tag{5.1.8}
$$
By \thetag{5.1.3}, in \thetag{5.1.8}, the second case is possible
only if $b_{11}\equiv b_{22}\equiv b_{33}\equiv b_{44}\equiv 1\mod 2$.
We denote (like in Sect. 4.3) by $G_0$ the lattice defined
by the primitive matrix $B$. Suppose that
$b_{11}\equiv b_{22}\equiv b_{33}\equiv b_{44}\equiv 1\mod 2$
and $\det(G_0))=2^k m$ where $m$ is odd and $k\equiv 1\mod 2$. We
then have the second case $\Gamma=2B$ in \thetag{5.1.8},
because otherwise $S$ is odd but $d$ is even.
Thus, for a fixed $B$, we may have both cases in
\thetag{5.1.8} only
if $b_{11}\equiv b_{22}\equiv b_{33}\equiv b_{44}\equiv 1\mod 2$
and $\det(G_0)=2^k m$ where $m$ is odd and $k\equiv 0\mod 2$.
In all other cases the matrix $B$ prescribes the Gram matrix
$\Gamma$. For both cases in \thetag{5.1.8},
we can calculate $d$ using the matrix $B$
since $d$ is square-free. If one of $\alpha_i^2$ is odd but $d$ is even,
then $\Gamma=B$, and this case is impossible for a main $S$ (the
lattice $S$ should be even for even $d$).

Summarizing considerations above, we get the following additional
conditions for $B$ and we calculate the invariants $(d,\eta)$ of the
lattice $S$ using the matrix $B$. We denote by $G_0$ the lattice
with the matrix $B$.
We have
$$
b_{11},\,b_{22},\,b_{33},\,b_{44}\ \text{are square-free}
\tag{5.1.9}
$$
and
$$
\text{g.c.d}(b_{11},b_{22})\le 2,\ \text{g.c.d}(b_{22},b_{33})\le 2,\
\text{g.c.d}(b_{33},b_{44})\le 2.
\tag{5.1.10}
$$
We denote by $\nu_p(a)$ the order of $p$ in factorisation of $a$.
We have
$$
\text{if odd\ }p\vert b_{11}b_{22}b_{33}b_{44},\ \text{then\ }
\nu_p(\det(G_0))\equiv 1\mod 2.
\tag{5.1.11}
$$
We have
$$
\Gamma=B,\ \text{if not all}\ b_{11},\,b_{22},\,b_{33},\,b_{44}\
\text{are odd};
\tag{5.1.12}
$$
$$
\Gamma=2B,\ \text{if\ }b_{11}\equiv b_{22}\equiv b_{33}\equiv b_{44}
\equiv 1\mod 2\ \text{and\ }\nu_2(det(G_0))\equiv 1\mod 2;
\tag{5.1.13}
$$
$$
\Gamma=B\ \text{or}\ \Gamma=2B,\
\text{if\ }b_{11}\equiv b_{22}\equiv b_{33}\equiv b_{44}
\equiv 1\hskip-2pt\mod \hskip-2pt2\
\text{and\ }\nu_2(det(G_0))\equiv 0\hskip-2pt\mod \hskip-2pt 2.
\tag{5.1.14}
$$
We denote by $G$ the lattice generated by $\alpha_i$ and
defined by the Gram matrix $\Gamma=\left((\alpha_i,\,\alpha_j)\right)$.
Let $d$ be the product of all $p$ such that $p\vert \det (G)$ and
$\nu_p(det(G))\equiv 1\mod 2$. We have
$$
\nu_2(d)\equiv 0\mod 2\ \text{if one of\ }\alpha_i^2\ \text{is odd.}
\tag{5.1.15}
$$
For odd $p\vert d$ we have
$$
(-1)^{\eta_p}=
\cases
\left({\alpha_1^2/p\over p}\right) &\text{if $p\vert\alpha_1^2$,}\\
\left({\alpha_2^2/p\over p}\right) &\text{if $p\vert\alpha_2^2$,}\\
\left({d\alpha_1^2\alpha_2^2/p\over p}\right)
&\text{otherwise.}
\endcases
\tag{5.1.16}
$$
This defines the invariant $\eta$. At last, we have the most delicate
condition which enormously drops the finite number of possibilities:
$$
h=hnr(d,\eta)\le 1
\tag{5.1.17}
$$
(see Theorem 3.2.1 about $hnr(d,\eta)$).
To calculate $hnr(d,\eta)$ using Theorem 3.2.1,
we should calculate class-numbers of imaginary quadratic fields
of the discriminant $D$ where $-4d\le D<0$. By \thetag{5.1.1},
$d\le 2\cdot 238569$ and $0<-D\le 8\cdot 238569=1908552$.
Moreover, by \thetag{5.1.1}, prime divisors of $D$
are not more than $238$. Thus, checking the condition \thetag{5.1.17},
we work with reasonable (not too big) integers. Or, if one wants,
we should use a program which correctly calculates
class-numbers of discriminants $D$ where $-1908552\le D<0$ and
$D$ is product of primes $p\le 238$.

In Appendix: Program 9: fund20.main, we give a program for
GP/PARI calculator, which enumerates all the matrices $B$
satisfying the conditions of Sect. 4.3 (this part of Program 9
is the same as the Program 7: fund20.gen) and conditions
\thetag{5.1.9} --- \thetag{5.1.17},
and gives all triplets
of invariants $(d,\eta,h)$ of the matrices $B$
(it gives 132 triplets).
We can see that all these triplets belong to the Table 3.
The GP/PARI calculator calculates
class numbers of negative discriminants $D$ for $|D|<10^{25}$ (see
user's guide to the calculator).
This finishes the proof of Theorem 5.1.1.
\enddemo

\smallpagebreak

By Theorem 5.1.1, to prove Theorem 2.3.2.1, we
should now check reflective type of all lattices of Table 3 and
calculate the sets
$P(\M)_{\pr}$ and the Gram matrices $G(P(\M)_{\pr})$ if they
are elliptically or parabolically reflective.
In Table 3, for each
pair of invariants $(d,\eta)$ we give a main lattice $S$ with
these invariants. If
$$
S=U\oplus \langle -d \rangle,
\tag{5.1.18}
$$
we have
$$
(-1)^{\eta_p}=\left({-d/p\over p}\right),\
\text{for any odd\ } p\vert d.
\tag{5.1.19}
$$
If
$$
S=\langle n_1\rangle \oplus \langle -n_2 \rangle\oplus
\langle -n_3 \rangle (\epsilon_1/2,\epsilon_2/2,\epsilon_3/2)
\tag{5.1.20}
$$
where $\epsilon_1,\epsilon_2,\epsilon_3 \in \{0,\,1\}$,
we have $d=n_1n_2n_3$ if
$\epsilon_1=\epsilon_2=\epsilon_3=0$, and
$d=n_1n_2n_3/4$ otherwise. Moreover, $n_i\equiv 0\mod 2$ if
$\epsilon_i=1$;
 $(n_1\epsilon_1-n_2\epsilon_2-n_3\epsilon_3)/4\in \bz$;
$(n_1\epsilon_1-n_2\epsilon_2-n_3\epsilon_3)/4\in 2\bz$ if
$d\equiv 0\mod 2$; for any odd $p\vert d$ we have
$$
(-1)^{\eta_p}=
\cases
\left({n_1/p\over p}\right) &\text{if $p\vert n_1$,}\\
\left({-n_2/p\over p}\right) &\text{if $p\vert n_2$,}\\
\left({-n_3/p\over p}\right) &\text{if $p\vert n_3$.}
\endcases
\tag{5.1.21}
$$
Checking all these conditions, one can prove that our calculation of main
lattices $S$ corresponding to the invariants $(d,\eta)$ of
Table 3 are correct. Using known
criteria (e. g. see \cite{Se}), it is easy to prove that
$S=U\oplus \langle -d \rangle$ if and only if $S$ represents $0$.
The condition \thetag{5.1.19} is equivalent to this property.
If \thetag{5.1.19} is valid, we always give $S$ in the form
$S=U\oplus \langle -d \rangle$.

Thus, fortunately, all our lattices $S$ have one of two forms
\thetag{5.1.18} or \thetag{5.1.20}.

We use the Vinberg's algorithm \cite{V2} to calculate the sets
$P(\M)_{\pr}$ for lattices of the forms \thetag{5.1.18} and
\thetag{5.1.20}. Below we describe this algorithm.

First, we remark that the
lattice $U\oplus \langle -d \rangle$ is equivariantly equivalent
to its maximal even sublattice $U\oplus \langle -4d \rangle$ if
$d$ is odd. These lattices have naturally isomorphic groups of
automorphisms and the reflection groups. Thus, it is sufficient
to consider only the lattices
$$
U\oplus \langle -2k \rangle,\ \ k\in \bn.
\tag{5.1.22}
$$

For lattices $S$ of the form \thetag{5.1.22},
we use the isotropic vector $c=(1,0,0)$
as the center of Vinberg's algorithm.
We find $\M$ which contains $\br_{++}c$
as an infinite vertex and $v_1=(0,0,1)$,
$v_2=(n,0,-1)$ as roots orthogonal to faces of
$\M$ containing $\br_{++}c$
(equivalently, they are roots of the height $0$).
It is easy to see that $v_1,v_2$ give the orthogonal primitive
roots to the fundamental polyhedron of the stabilizer subgroup of
$c$ in the reflection group $W(S)$.

For lattices $S$ of the form \thetag{5.1.20} we
assume that either $\epsilon=(\epsilon_1,\epsilon_2,\epsilon_3)/2=0$ or
$\epsilon_i\not=0$ at least for two different $i$.
We take $c=(1,0,0)$ as
the center of Vinberg's algorithm.
We take as roots of the height $0$
\newline
$v_1=(0,1,0)$ and $v_2=(0,0,1)$ if $n_2/n3\not=1/3,\,1,\,3$;
\newline
$v_1=(0,1,0)$, $v_2=(0,-1,1)$ if
$n_2=n_3$ and $\epsilon=(0,0,0),\,(1/2,1/2,1/2)$;
\newline
$v_1=(0,1,0)$, $v_2=(0,-1/2,1/2)$ if $n_2=n_3$ 	and
$\epsilon=(0,1/2,1/2)$;
\newline
$v_1=(0,1,0)$, $v_2=(0,0,1)$ if $n_2=n_3$ and
$\epsilon\not=(0,0,0),\,(1/2,1/2,1/2),\,(0,1/2,1/2)$;
\newline
$v_1=(0,1,0)$, $v_2=(0,-1/2,1/2)$ if $n_2/n_3=3$ and
$\epsilon=(0,1/2,1/2)$;
\newline
$v_1=(0,1,0)$, $v_2=(0,0,1)$ if $n_2/n_3=3$ and
$\epsilon\not=(0,1/2,1/2)$;
\newline
$v_1=(0,0,1)$, $v_2=(0,1/2,-1/2)$ if $n_2/n_3=1/3$ and
$\epsilon=(0,1/2,1/2)$;
\newline
$v_1=(0,1,0)$, $v_2=(0,0,1)$ if $n_2/n_3=1/3$ and
$\epsilon\not=(0,1/2,1/2)$.
\newline
One can check that $v_1,v_2$ are the orthogonal primitive roots
to the fundamental polyhedron of the stabilizer subgroup of $c$ in $W(S)$.

Further steps of the Vinberg's
algorithm are prescribed canonically (see \cite{V2}). One introduces
the height
$$
\text{height}={2(\delta,\,c)^2\over -\delta^2} \in \bz
\tag{5.1.23}
$$
of primitive roots $\delta \in S$. The roots
$\{v_1,\,v_2\}$ above give all roots of $P(\M)_{\pr}$ of the height $0$.
If one knows all roots of $P(\M)_{\pr}^{\ \le n}\subset P(\M)_{\pr}$
of the height $\le n$,
all primitive roots $P(\M)_{\pr}^{\ n+1}\subset P(\M)_{\pr}$
of the height $n+1$ are given by the condition
$$
{2(\delta,\,c)^2\over -\delta^2}=n+1\ \text{and\ }
(\delta,\,P(\M)_{\pr}^{\ \le n})\ge 0.
\tag{5.1.24}
$$
In Appendix, Program 10: refl0.1, we give this algorithm for
lattices \thetag{5.1.22}. In Appendix, Program 12: refl0.13,
we give this algorithm for lattices \thetag{5.1.20}. We also
give Program 11: refl0.12 which calculates only for
$\epsilon_1=\epsilon_2=\epsilon_3=0$ but is faster. These programs
calculate a sequence $v_i$ of all elements of the subset
$P(\M)_{\pr}^{\ \le n}\subset P(\M)_{\pr}$ (i. e. all roots of the
height $\le n$).

After calculation of $P(\M)_{\pr}^{\ \le n}$ for a sufficiently
large height $n$, next steps of our algorithm are as
follows. A sequence $r=r_1,\dots,r_k$ of different elements
from $P(\M)_{\pr}$ is called
a {\it chain} if
$$
0\le {4(r_i,r_{i+1})^2\over r_i^2r_{i+1}^2}\le 2,
\ \text{for any\ } 1\le i\le k-1.
\tag{5.1.25}
$$
Geometrically it means that $\Ha_{r_1},\dots,\Ha_{r_k}$ are
lines of consecutive sides (i. e. defining vertices) of the polygon
$\M$. We find the maximal chain
$e=(e_1,\dots, e_k)$ in $P(\M)_{\pr}^{\ \le n}$ containing $v_1$.
Suppose that (for a sufficiently large $n$) we also have
$$
{4(e_1,e_k)^2\over e_1^2e_k^2}\le 2,
\tag{5.1.26}
$$
i. e. the lines $\Ha_{e_k}$ and $\Ha_{e_1}$ also
define a vertex of $\M$. Then the polygon $\M$ is elliptic,
$P(\M)_{\pr}=e$, and the chain $e$ gives the orthogonal primitive roots to
consecutive sides of the elliptic polygon $\M$. This situation
takes place for $122$ invariants $(d,\eta)$ of the Table 3 marked by $er$
(elliptically reflective cases). The result of our calculation of the
chain $e$ and the Gram matrix $\left((e_i,e_j)\right)$, $e_i,\,e_j\in e$,
for these $122$ cases is given in Table 1.

For all other $206-122=84$ cases we should prove that $S$ is
not reflective of elliptic or parabolic type.

To prove that, for small determinants $d$ we use the following arguments.
For a large height $n$, the chain $e\subset P(\M)_{\pr}^{\ \le n}$
contains a {\it period:} there exists $1<q<k$ such that
the Gram matrices $\Gamma(e_1,\,e_2)$ and
$\Gamma(e_{q},\,e_{q+1})$ coincide and there exists
$C\in O^+(S)$ such that $C(e_1)=e_q$, $C(e_2)=e_{q+1}$.
It follows that $C\in A(\M)$. We
find all integral (i.e. from $S$) eigenvectors $u$ of $C$,
and show that all of them have $u^2<0$. Clearly, these eigenvectors
have the eigenvalue $\pm 1$.
It follows that $S$ is not elliptically or parabolically reflective.
Really, if $S$ is elliptically reflective, the group $A(\M)$
is finite, $C$ has finite order and has
an integral eigenvector with the eigenvalue $1$ and
positive square. If $S$ is parabolically
reflective, then $A$ has an eigenvector with the eigenvalue $1$ and
with zero square (the generalized lattice Weyl vector $\rho$).
See the proof above of Lemma 2.3.1.2. If $S$ is hyperbolically
reflective, the eigenvector $u$ defines a generalized
lattice Weyl vector $\rho$.

Let us consider an example of these calculations for
$(d,\eta)=(114,2)$. Then $S=U\oplus \langle -114 \rangle$.
Using Program 10 we calculate up to the height $n=50000$ and
find the chain
$e$ in $P(\M)_{\pr}^{\ \le n}$ which is equal to
$$
\split
&e_1=(321,30,-13),\ e_2=(28,2,-1),\ e_3=(57,0,-1),\ e_4=(0, 0, 1),\\
&e_5=(-1, 1, 0),\ e_6=(9, 6, -1),\ e_7=(292,122,-25)
\endsplit
\tag{5.1.27}
$$
with the Gram matrix
$$
ge=
\pmatrix
-6& 0& 228& 1482& 291& 714& 10872\\
0 &-2& 0 &  114 & 26 &  72& 1150\\
228 &0 &-114& 114& 57 &228& 4104\\
1482& 114 &114& -114& 0 & 114& 2850\\
291 & 26  &57 &   0 &-2 &  3 &  170\\
714 &72& 228& 114& 3 &-6& 0\\
10872& 1150 &4104& 2850& 170 &0 &-2
\endpmatrix .
\tag{5.1.28}
$$
We see that the Gram matrices of $e_1,e_2$ and $e_6,e_7$ coincide.
We calculate that
$$
C=
\pmatrix
2209 & 22800 & 107160\\
912 &  9409  & 44232\\
-188& -1940  &-9119
\endpmatrix
\tag{5.1.29}
$$
belongs to $O^+(S)$ and $C(e_1^t,e_2^t)=(e_6^t,e_7^t)$. It follows that
$C\in A(\M)$. The matrix $C$ has the only integral eigenvector
$w=(95,19,-6)$, it has the eigenvalue $1$. We have $w^2=-494$.
It follows that $S$ is not elliptically or parabolically reflective:
the automorphism $C$ has infinite order and does not have integral
eigenvectors with non-negative square.

Since for $(d,\eta)=(114,2)$ the invariant $h=0$, it is possible that
$S$ is hyperbolically reflective with the generalized lattice
Weyl vector $\rho=w$. Let us prove that this is the case.
We find another chain $f$ in $P(\M)_{\pr}^{\ \le n}$ which is
equal to
$$
\split
&f_1=(1766,172,-73),\  f_2=(6384,627,-265),\  f_3=(4560, 456, -191),\\
&f_4=(283,29,-12),\ f_5=(18, 3, -1),\ f_6=(14, 4, -1),\\
&f_7=(456, 171, -37),\ f_8=(2280,912,-191),\ f_9=(427, 173, -36)
\endsplit
\tag{5.1.30}
$$
with the Gram matrix
$$
gf=
\left(\smallmatrix
-2 &0 &114& 26& 72 &1150& 72504 &413250& 79370\\
0& -114 &114& 57& 228 &4104& 259806& 1481658& 284601\\
114& 114& -114& 0& 114& 2850 &182058& 1039566 &199728\\
26& 57 & 0& -2& 3& 170 &11001& 62928& 12094\\
72& 228& 114& 3 &-6& 0& 228 &1482& 291\\
1150& 4104& 2850& 170& 0 &-2& 0& 114& 26\\
72504& 259806& 182058& 11001& 228& 0& -114& 114& 57\\
413250& 1481658& 1039566& 62928 &1482& 114 &114& -114& 0\\
79370& 284601& 199728& 12094& 291& 26 &57& 0& -2
\endsmallmatrix\right).
\tag{5.1.31}
$$
We see that the Gram matrix of $f_1,f_2$ is equal to the Gram matrix
of $f_6,f_7$. We calculate that $C(f_1^t,f_2^t)=(f_6^t,f_7^t)$.
Moreover, we can see that $(e_i,w)>0$ and $(f_i,w)<0$.

It follows that $P(\M)_{\pr}=H(e)\cup H(f)$ where $H=[C]$ is
the infinite cyclic group generated by the $C$. The infinite polygon
$\M$ contains the  line (the axis) $\Ha_w$ which is
preserved by $H$. All sides of $\M$ define two infinite
chains $\Ha_\delta,\,\delta\in H(e)$ and $\Ha_\delta,\,\delta \in H(f)$,
these two chains are contained in
two different half-planes bounded by the line $\Ha_w$.
The polygon $\M$ is finite in every orthogonal cylinder over a compact
base in $\Ha_w$ (it is restricted hyperbolic relative to $\Ha_w$).
The polygon $\M$ has an infinite symmetry
group $H=[C]$. It is not the full symmetry group $A(\M)$ of $\M$.
After considerations above, it is easy to see
that $A(\M)$ is the cyclic group generated by the glide reflection
$C_1$ with axis $\Ha_w$ where
$$
C_1=\pmatrix
{48}&{475}&{2280}\cr{19}&{192}&{912}\cr{-4}&{-40}&{-191}
\endpmatrix.
\tag{5.1.32}
$$
The vector $w$ is the eigenvector of $C_1$ with the eigenvalue $-1$.
We have $C_1^2=C$. The glide reflection $C_1$ changes places
the infinite chains $H(e)$ and $H(f)$. The lattice $S$ is hyperbolically
reflective with the generalized lattice Weyl vector $\rho=w$.

\smallpagebreak

For large $d$ following arguments are very useful.
We calculate the set $P(\M)_{\pr}^{\ \le n}$ for a sufficiently large
height $n$ and find two pairs of elements
$r_1,r_2\in P(\M)_{\pr}^{\ \le n}$,
$s_1,s_2\in P(\M)_{\pr}^{\ \le n}$ such that
$4(r_1,r_2)^2/r_1^2r_2^2\le 2$,
$4(s_1,s_2)^2/s_1^2s_2^2\le 2$ and Gram matrices of
these pairs coincide. In many cases it is sufficient to
consider pairs with $(r_1,r_2)=(s_1,s_2)=0$ and
$r_1^2=s_1^2$, $r_2^2=s_2^2$. Thus, orthogonal sides of these pairs
define two vertices of $\M$. We find such pairs that there exists
an automorphism $B\in O^+(S)$ such that $B(r_i)=s_i$.
Then $B\in A(\M)$. Similarly, considering two other pairs
in $P(\M)_{\pr}^{\ \le n}$, we find another $C\in A(\M)$. If
$B^2C^2\not=C^2B^2$, the lattice $S$ is not reflective of any type:
elliptic, parabolic or hyperbolic. In many cases it is sufficient
to find only one $B\in A(\M)$ and calculate that $B^{12}\not=E$.
Really, then $B$ has infinite order and the lattice $S$ cannot be
elliptically reflective. If $h=1$, the lattice $S$ also
cannot be parabolically or hyperbolically reflective.

We consider an example of these calculations for
$(d,\eta)=(3990,4)$. Then
$S=\langle 30\rangle\oplus \langle -38\rangle \oplus
\langle-14 \rangle(1/2,1/2,0)$. Using Program 12 from Appendix,
we calculate up to the height $n=500000$. It gives $33$ elements
$v_i\in P(\M)_{\pr}^{\ \le n}$. They are
$$
\split
&v_1=(0,1,0),\ v_2=(0,0,1),\ v_3=(1/2,-1/2,0),\\
&v_4=(2,0,-3),\ v_5=(13/2,-9/2,-6),\ v_6=(7,-3,-9),\\
&v_7=(57/2,-21/2,-38),\ v_8=(17/2,-9/2,-10),\ v_9=(28,-21,-22),\\
&v_{10}=(63/2,-35/2,-36),\ v_{11}=(57,-37,-57),\ v_{12}=(25,-15,-27),\\
&v_{13}=(42,-14,-57),\ v_{14}=(17,-4,-24),\ v_{15}=(76,-58,-57),\\
&v_{16}=(19,-6,-26),\ v_{17}=(52,-15,-72),\ v_{18}=(84, -21, -118),\\
&v_{19}=(1729/2, -1015/2, -950),\ v_{20}=(119/2, -69/2, -66),\\
&v_{21}=(69/2,-43/2,-36),\\
&v_{22}=(73,-57,-51),\ v_{23}=(74, -57, -54),\
v_{24}=(231/2,-147/2,-118),\\
&v_{25}=(2261/2, -1645/2, -950),\ v_{26}=(101,-27,-141),\
v_{27}=(266,-132,-323),\\
&v_{28}=(119,-87,-99),\ v_{29}=(128,-63,-156),\
v_{30}=(399/2,-315/2,-134),\\
&v_{31}=(342,-96,-475),\ v_{32}=(361,-77,-513),\ v_{33}=(238,-119,-288).
\endsplit
\tag{5.1.33}
$$
We look at all the pairs $v_i,v_j$ with $(v_i,v_j)=0$, and we
find that there are five such pairs $v_1,v_2$; $v_7,v_{13}$; $v_{15},v_9$;
$v_{11},v_{24}$; $v_{27},v_{33}$ having squares $v_i^2=-38$, $v_j^2=-14$.
We calculate a primitive orthogonal element $w\in S$ to each of
these five pairs $v_i,v_j$, and we find that between these five pairs
there are two pairs
$v_{11},v_{24}$ and $v_{27},v_{33}$ such that
$(w+v_j)/2\in S$. We then find the matrix
$$
C=\pmatrix
6863/2 & 5339/2 & 1694\\
-3345/2 & -2601/2  &-826\\
-4200  &-3268  &-2073
\endpmatrix
\tag{5.1.34}
$$
such that $C\in O^+(S)$ and $C(v_{11}^t,v_{24}^t)=(v_{27}^t,v_{33}^t)$.
It follows that $C\in A(\M)$. We have $C^{12}\not=E$.
It follows that $C$ has infinite order.
Thus, the lattice $S$ is not
elliptically reflective. For $(d,\eta)=(3990,4)$,
the invariant $h=1$, thus
the lattice $S$ cannot be also parabolically or hyperbolically reflective.
It follows that $S$ is not reflective.

These arguments permit to prove that the rest $84$ cases of Table 3
which are not contained in the Table 1, are not reflective of elliptic
or parabolic type. It finishes the proof.

\remark{Remark 5.1.2} Here we want to outline another way
which helps to study reflective type of lattices of Table 3.
We can write similar programs as Program 9: fund20.main for all
types of narrow places of $\M$ of main lattices $S$ having
the invariant $h\le 1$. They would be programs
fund11.main, fund10.main, fund21.main and fund 30.main specializing
the programs fund11.gen, fund10.gen, fund21.gen and fund30.gen
(one can write them similarly to Program 9: fund20.main).
Together with Program 9: fund20.main, they give a list
L of invariants $(d,\eta,h)$ containing in Table 3.
All invariants $(d,\eta,h)$ of Table 3 which are not in the list
L cannot be elliptically or parabolically reflective. If $h=1$,
they cannot be hyperbolically reflective either.

We did this calculations and we found that the list L is much smaller than
Table 3 (e. g.  Program 9: fund20.main gives only $132$ triplets
$(d,\eta,h)$).
It shows that the narrow places of polyhedra arguments 
are sometimes stronger than 
arithmetic arguments of studying the invariant $h$ (even if we forget 
about the problem with infinity). Both these arguments
surprisingly fit together.

Here is the list of invariants $(d,\eta,h)$ which are contained in
Table 3 but are not in the list L (39 pairs):
$$
\split
&(d,\eta,h)=\\
&(57,1,1);\ (65,3,1);\ (71,0,1);\ (119,3,1);\ (161,1,1);\ (182,3,1);\\
&(194,0,1);\ (246,0,1);\ (259,3,1);\ (266,0,1);\ (266,3,1);\ (285,1,1);\\
&(299,3,1);\ (326,0,1);\ (335,0,1);\ (354,2,1);\ (386,0,1);\ (407,0,1);\\
&(506,0,1);\ (530,0,1);\ (534,0,1);\ (546,2,0);\ (602,3,1);\ (645,6,1);\\
&(714,6,1);\ (777,6,1);\ (854,3,1);\ (897,5,1);\ (897,7,1);\ (935,6,1);\\
&(966,2,1);\ (1106,1,1);\ (1254,4,1);\ (1394,3,1);\ (1659,2,1);\
(2210,6,1);\\
&(3311,1,1);\ (3990,4,1);\ (4466,1,1).
\endsplit
\tag{5.1.35}
$$
These $(d,\eta,h)$ give lattices $S$ which are not elliptically or
parabolically reflective: their fundamental polygon $\M$ does not have a 
narrow place satisfying Theorem 4.2.3. Of course, the list
\thetag{5.1.35} is in complete agreement with calculations above
using Vinberg's algorithm.
\endremark

\subhead
5.2. Proof of Theorem 2.3.3.1
\endsubhead

We should check reflective type of non-main lattices $\widetilde{S}$
corresponding to main lattices $S$ of Table 1 with odd $d$ (there
are $97$ these cases). See Proposition 2.2.6. Let $S$ be one of lattices
of Table 1 (or of Table 3 marked by er) with invariants $(d,\eta)$
where $d$ is odd. Then
$\widetilde{S}$ has invariants $(2d,odd,\eta+\omega(p))$.
If $S$ and $\widetilde{S}$ are equivariantly equivalent, i. e.
$$
\sum_{p\mid d}
{(1-p+4\eta_p+4\omega(p))} \equiv 0\ \text{or}\  6 \mod 8,
\tag{5.2.1}
$$
then lattices $S$ and $\widetilde{S}$ have the same reflective type
and calculation of $P(\M)_{\pr}$ and its Gram matrix for the lattice
$\widetilde{S}$ follows from similar
calculation for $S$ (see Remark 2.3.3.2). There are $21$ these cases.
For example if $S=U\oplus \langle -d \rangle$, where $d$ is odd, then
$\widetilde{S}=\langle 1\rangle\oplus \langle -1 \rangle\oplus
\langle -2d\rangle$, and these lattices are equivariantly equivalent.
Thus, we need to study only cases when \thetag{5.2.1} is not valid.
There are $97-21=76$ these cases. Let $S$ be one of these lattices.

We have: if $S=\langle n_1\rangle \oplus \langle -n_2\rangle
\oplus \langle -n_3 \rangle$, where all $n_i$ are odd, then
$\widetilde{S}=\langle 2n_1\rangle \oplus \langle -2n_2\rangle
\oplus \langle -2n_3 \rangle(\epsilon_1,\epsilon_2,\epsilon_3)$
where one of $\epsilon_i$ is equal to $0$, two of them are $1/2$
and $(\epsilon_1n_1-\epsilon_2n_2-\epsilon_3n_3)/2$ is odd.
If $S=\langle n_1\rangle \oplus \langle -n_2\rangle
\oplus \langle -n_3 \rangle(\epsilon_1,\epsilon_2,\epsilon_3)$
where for example $\epsilon_1=\epsilon_2=1/2$ and $\epsilon_3=0$
(it follows that $n_1\equiv n_2\equiv 2 \mod 4$ and $n_3$ is odd),
then $\widetilde{S}=\langle n_1/2\rangle \oplus \langle -n_2/2\rangle
\oplus \langle -2n_3 \rangle$.
We see that the lattice $\widetilde{S}$ which we should check for
reflective type, has the form \thetag{5.1.20}. Thus, we need
to make similar calculations as in Sect. 5.1. They are of the same
difficulty. These calculations finish the proof.

\head
6. Appendix: Programs for GP/PARI calculator
\endhead

\centerline{{\bf Program 1:} h2}

\smallpagebreak

\noindent
$\backslash\backslash$hclass(d,muu)
calculates h=(hrI(d,muu),hrII(d,muu),hnr(d,muu))
\newline
$\backslash\backslash$here d$<$0 and 0$\backslash$le muu$<$2$^\wedge$k
are integers,
\newline
$\backslash\backslash$k is the number of all odd prime divisors of d.
\newline
$\backslash\backslash$Assume d is a fundamental discriminant
of fundamental binary
\newline
$\backslash\backslash$positive lattices
(i.e. with a square free determinant),
\newline
$\backslash\backslash$then d$\backslash$equiv 1 $\backslash$mod 4 or
\newline
$\backslash\backslash$d $\backslash$equiv
$\backslash$pm 4, 8$\backslash$mod
16.
\newline
$\backslash\backslash$ Assume (d,muu) is the genus of them. Then
\newline
$\backslash\backslash$hrI(d,muu), hrII(d,muu), hnr(d,muu)
are numbers of ambiguous
\newline
$\backslash\backslash$classes of the types I, II and
\newline
$\backslash\backslash$hnr(d,muu)=(h(d)/2$^\wedge$\{tau(d)\}--
hrI(d,muu)--hrII(d,muu))/2
\newline
$\backslash\backslash$of non-ambiguous classes of the general
\newline
$\backslash\backslash$equivalence respectively of the genus (d,muu);
\newline
$\backslash\backslash$hclass(d,muu)=[0,0,0] if
\newline
$\backslash\backslash$d is not a fundamental discriminant.
\newline
$\backslash\backslash$if d is a fundamental discriminant but
0$\backslash$le muu $<$2$^\wedge$k
\newline
$\backslash\backslash$does not correspond to a
\newline
$\backslash\backslash$genus,
then hnr(d,muu)=[0,0,h(d)/2$^\wedge$\{tau(d)+1\}]
\newline
hclass(d,muu,dd,fdd,dd1,beta,alpha,k,hr,hrI,hrII,t,h)=$\backslash$
\newline
h=[0,0,0];hr=0;hrI=0;hrII=0;$\backslash$
\newline
if(mod(d,4)!=mod(1,4),,dd=--d;$\backslash$
\newline
fdd=factor(dd);$\backslash$
\newline
fordiv(dd,dd1,if(dd1$>$dd/dd1,,$\backslash$
\newline
beta=1;alpha=1;k=1;$\backslash$
\newline
while(alpha,$\backslash$
\newline
if(type(dd1/fdd[k,1])==1,$\backslash$
\newline
if(kro(2$\ast$dd1/fdd[k,1],fdd[k,1])==(--1)$^\wedge$bittest(muu,k--1),
$\backslash$
\newline
if(k$>$=matsize(fdd)[1],alpha=0,k=k+1),beta=0;alpha=0),$\backslash$
\newline
if(kro(2$\ast$(dd/dd1)/fdd[k,1],fdd[k,1])==(--1)$
^\wedge$bittest(muu,k--1),$\backslash$
\newline
if(k$>$=matsize(fdd)[1],alpha=0,k=k+1),beta=0;alpha=0)));$\backslash$
\newline
hr=hr+beta));t=matsize(fdd)[1]--1;$\backslash$
\newline
h=[0,hr,(classno(d)/2$^\wedge$t--hr)/2]);$\backslash$
\newline
if(mod(d,16)!=mod(4,16),,dd=--d/4;$\backslash$
\newline
fdd=factor(dd);$\backslash$
\newline
fordiv(dd,dd1,if(dd1$>$dd/dd1,,$\backslash$
\newline
beta=1;alpha=1;k=1;$\backslash$
\newline
while(alpha,$\backslash$
\newline
if(type(dd1/fdd[k,1])==1,$\backslash$
\newline
if(kro(dd1/fdd[k,1],fdd[k,1])==(--1)$^\wedge$bittest(muu,k--1),$\backslash$
\newline
if(k$>$=matsize(fdd)[1],alpha=0,k=k+1),beta=0;alpha=0),$\backslash$
\newline
if(kro((dd/dd1)/fdd[k,1],fdd[k,1])==(--1)$^\wedge$bittest(muu,k--1),
$\backslash$
\newline
if(k$>$=matsize(fdd)[1],alpha=0,k=k+1),beta=0;alpha=0)));$\backslash$
\newline
hr=hr+beta));t=matsize(fdd)[1]--1;$\backslash$
\newline
h=[hr,0,(classno(d)/2$^\wedge$t--hr)/2]);$\backslash$
\newline
if(mod(d,16)!=mod(8,16),,dd=--d/4;$\backslash$
\newline
if(dd/2==1,hr=1;h=[1,0,0],fdd=factor(dd/2);$\backslash$
\newline
fordiv(dd/2,dd1,$\backslash$
\newline
beta=1;alpha=1;k=1;$\backslash$
\newline
while(alpha,$\backslash$
\newline
if(type(dd1/fdd[k,1])==1,$\backslash$
\newline
if(kro(dd1/fdd[k,1],fdd[k,1])==(--1)$^\wedge$bittest(muu,k--1),$\backslash$
\newline
if(k$>$=matsize(fdd)[1],alpha=0,k=k+1),beta=0;alpha=0),$\backslash$
\newline
if(kro((dd/dd1)/fdd[k,1],fdd[k,1])==(--1)$^\wedge
$bittest(muu,k--1),$\backslash$
\newline
if(k$>$=matsize(fdd)[1],alpha=0,k=k+1),beta=0;alpha=0)));$\backslash$
\newline
hr=hr+beta);t=matsize(fdd)[1];$\backslash$
\newline
h=[hr,0,(classno(d)/2$^\wedge$t--hr)/2]));$\backslash$
\newline
if(mod(d,16)!=mod(--4,16),,dd=--d/4;$\backslash$
\newline
if(dd==1,h=[1,1,0],fdd=factor(dd);$\backslash$
\newline
fordiv(dd,dd1,if(dd1$>$dd/dd1,,$\backslash$
\newline
beta=1;alpha=1;k=1;$\backslash$
\newline
while(alpha,$\backslash$
\newline
if(type(dd1/fdd[k,1])==1,$\backslash$
\newline
if(kro(dd1/fdd[k,1],fdd[k,1])==(--1)$^\wedge$bittest(muu,k--1),$\backslash$
\newline
if(k$>$=matsize(fdd)[1],alpha=0,k=k+1),beta=0;alpha=0),$\backslash$
\newline
if(kro((dd/dd1)/fdd[k,1],fdd[k,1])==(--1)$^\wedge
$bittest(muu,k--1),$\backslash$
\newline
if(k$>$=matsize(fdd)[1],alpha=0,k=k+1),beta=0;alpha=0)));$\backslash$
\newline
hrI=hrI+beta));$\backslash$
\newline
fordiv(dd,dd1,if(dd1$>$dd/dd1,,$\backslash$
\newline
beta=1;alpha=1;k=1;$\backslash$
\newline
while(alpha,$\backslash$
\newline
if(type(dd1/fdd[k,1])==1,$\backslash$
\newline
if(kro(2$\ast$dd1/fdd[k,1],fdd[k,1])==(--1)$^\wedge
$bittest(muu,k--1),$\backslash$
\newline
if(k$>$=matsize(fdd)[1],alpha=0,k=k+1),beta=0;alpha=0),$\backslash$
\newline
if(kro(2$\ast$(dd/dd1)/fdd[k,1],
fdd[k,1])==(--1)$^\wedge$bittest(muu,k--1),$\backslash$
\newline
if(k$>$=matsize(fdd)[1],alpha=0,k=k+1),beta=0;alpha=0)));$\backslash$
\newline
hrII=hrII+beta));$\backslash$
\newline
hr=hrI+hrII;t=matsize(fdd)[1];$\backslash$\
\newline
h=[hrI,hrII,(classno(d)/2$^\wedge$t--hr)/2]));h;

\smallpagebreak

\centerline{{\bf Program 2:} h3}

\smallpagebreak

\noindent
$\backslash\backslash$the program h3
\newline
$\backslash\backslash$hnr(d,et) calculates the number of classes of
\newline
$\backslash\backslash$non-reflective central
symmetries of a 3-dimensional main
\newline
$\backslash\backslash$hyperbolic lattice with the square-free determinant
\newline
$\backslash\backslash$d and the invariant et (a non-negative integer)
\newline
$\backslash\backslash$whose binary
decomposition et\_$\{$p\_k$\}$....et\_$\{$p\_1$\}$
\newline
$\backslash\backslash$gives the map of all odd prime
divisors p\_1,...,p\_k
\newline
$\backslash\backslash$of the d in increasing order to $\{$0,1$\}$

\smallpagebreak

\noindent
$\backslash$r h2
\smallpagebreak

\noindent
$\backslash\backslash$checking the condition (5)=(3.2.7)
\newline
beta5(d,et,n,fd,sfd,alpha,k,b,etap)=$\backslash$
\newline
fd=factor(d);sfd=matsize(fd)[1];alpha=1;k=1;$\backslash$
\newline
while(alpha,$\backslash$
\newline
if(k$>$sfd,b=1;alpha=0,$\backslash$
\newline
if(type(n/fd[k,1])!=1,k=k+1,$\backslash$
\newline
if(fd[k,1]==2,k=k+1,$\backslash$
\newline
if(fd[1,1]==2,etap=bittest(et,k--2),etap=bittest(et,k--1));$\backslash$
\newline
if(kro(n/fd[k,1],fd[k,1])==(--1)$^\wedge$etap,k=k+1,b=0;alpha=0)))));b;

\smallpagebreak

\noindent
$\backslash\backslash$checking the condition (6)=(3.2.8)
\newline
beta6(d,et,n,fd,sfd,u,k,b)=$\backslash$
\newline
fd=factor(d);sfd=matsize(fd)[1];$\backslash$
\newline
u=mod(0,8);$\backslash$
\newline
for(k=1,sfd,$\backslash$
\newline
if(type((d/n)/fd[k,1])!=1,,u=u+mod(1--fd[k,1]+4$\ast
$bittest(et,k--1),8)));$\backslash$
\newline
if(u!=mod(--2,8),b=0,b=1);b;

\smallpagebreak

\noindent
$\backslash\backslash$checking the condition (8)=(3.2.10)
\newline
beta8(d,et,n,fd,sfd,u,k,b)=$\backslash$
\newline
fd=factor(d);sfd=matsize(fd)[1];$\backslash$
\newline
u=mod(0,8);$\backslash$
\newline
for(k=1,sfd,$\backslash$
\newline
if(type((d/n)/fd[k,1])!=1,,u=u+mod(1--fd[k,1]+
4$\ast$bittest(et,k--1),8)));$\backslash$
\newline
u=u+mod(((d/n)$^\wedge$2--1)/2,8);$\backslash$
\newline
if(u!=mod(2,8),b=1,b=0);b;

\smallpagebreak

\noindent
$\backslash\backslash$checking the condition (11)=(3.2.13)
\newline
beta11(d,et,n,fd,sfd,u,k,b)=$\backslash$
\newline
fd=factor(d);sfd=matsize(fd)[1];$\backslash$
\newline
u=mod(0,4);$\backslash$
\newline
for(k=2,sfd,u=u+mod(1--fd[k,1]+4$\ast$bittest(et,k--2),4));$\backslash$
\newline
u=--u--mod(1,4);$\backslash$
\newline
if(u!=mod(n/2,4),b=0,b=1);b;

\smallpagebreak

\noindent
$\backslash\backslash$calculation of the numbers et\_p+epsilon(p)
if odd p$|$d
\newline
eps(d,et,e,fd,sfd,sfd1,e1,et1,eet1,k)=$\backslash$
\newline
if(d==1,e=0,$\backslash$
\newline
fd=factor(d);sfd=matsize(fd)[1];$\backslash$
\newline
if(fd[1,1]$>$2,sfd1=sfd,sfd1=sfd--1);$\backslash$
\newline
e1=vector(sfd1,k,mod((fd[sfd--k+1,1]--1)/2,2));$\backslash$
\newline
et1=vector(sfd1,k,mod(bittest(et,sfd1--k),2));$\backslash$
\newline
eet1=e1+et1;eet1=lift(eet1);$\backslash$
\newline
e=0;for(k=0,sfd1--1,e=e+2$^\wedge$k$\ast$eet1[sfd1--k]));e;

\smallpagebreak

\noindent
$\backslash\backslash$calculation of
the numbers et\_p+epsilon(p)+omega(p) if odd p$|$d
\newline
epsomeg(d,et,e,fd,sfd,sfd1,e1,et1,eet1,k)=$\backslash$
\newline
if(d==1,e=0,$\backslash$
\newline
fd=factor(d);sfd=matsize(fd)[1];$\backslash$
\newline
if(fd[1,1]$>$2,sfd1=sfd,sfd1=sfd--1);$\backslash$
\newline
e1=vector(sfd1,k,mod((fd[sfd--k+1,1]--1)/2+(fd[sfd--k+1,1]$^\wedge
$2--1)/8,2));$\backslash$
\newline
et1=vector(sfd1,k,mod(bittest(et,sfd1--k),2));$\backslash$
\newline
eet1=e1+et1;eet1=lift(eet1);$\backslash$
\newline
e=0;for(k=0,sfd1--1,e=e+2$^\wedge$k$\ast$eet1[sfd1--k]));e;

\smallpagebreak

\noindent
$\backslash\backslash$calculation of the numbers mu\_p=eta\_p
if odd p$|$t$|$2d
\newline
muuu(d,et,t,m,fd,sfd,sfd1,k1,k)=$\backslash$
\newline
if(d<=2,m=0,m=0;fd=factor(d);sfd=matsize(fd)[1];$\backslash$
\newline
if(fd[1,1]$>$2,sfd1=1,sfd1=2);k1=0;$\backslash$
\newline
for(k=0,sfd--sfd1,$\backslash$
\newline
if(type(t/fd[k+sfd1,1])==1,m=m+2$^\wedge$k1$\ast$bittest(et,k);k1=k1+1,)));m;

\smallpagebreak

\noindent
$\backslash\backslash$number of classes of non-reflective central
symmetries
\newline
hnr(d,et,h,n,n1)=$\backslash$
\newline
h=0;$\backslash$
\newline
if(mod(d,2)==mod(1,2),$\backslash$
\newline
fordiv(d,n,$\backslash$
\newline
if(beta5(d,et,n)==1\&\&beta6(d,et,n)==1,$\backslash$
\newline
h=h+hclass(--d/n,muuu(d,eps(d,et),d/n))[3],));$\backslash$
\newline
fordiv(d,n,$\backslash$
\newline
if(beta5(d,et,n)==1\&\&beta8(d,et,n)==1,$\backslash$
\newline
h=h+hclass(--4$\ast$d/n,muuu(d,eps(d,et),d/n))[3],));$\backslash$
\newline
fordiv(d,n1,$\backslash$
\newline
n=2$\ast$n1;$\backslash$
\newline
if(beta5(d,et,n)==1,$\backslash$
\newline
h=h+hclass(--16$\ast$d/n,muuu(d,eps(d,et),d/n1))[3],)),$\backslash$
\newline
fordiv(d/2,n1,$\backslash$
\newline
n=2$\ast$n1;$\backslash$
\newline
if(beta5(d,et,n)==1\&\&beta11(d,et,n)==1,$\backslash$
\newline
h=h+hclass(--d/n,muuu(d,eps(d,et),d/n))[3],));$\backslash$
\newline
fordiv(d/2,n1,$\backslash$
\newline
n=2$\ast$n1;$\backslash$
\newline
if(beta5(d,et,n)==1\&\&mod(n1,4)==mod(--d/2,4),$\backslash$
\newline
h=h+hclass(--4$\ast$d/n,muuu(d,epsomeg(d,et),d/n))[3],));$\backslash$
\newline
fordiv(d/2,n1,$\backslash$
\newline
n=2$\ast$n1;$\backslash$
\newline
if(beta5(d,et,n)==1\&\&mod(n1,4)==mod(d/2,4)\&\&n<d,$\backslash$
\newline
h=h+2$\ast$hclass(--4$\ast$d/n,muuu(d,epsomeg(d,et),d/n))[3]+$\backslash$
\newline
hclass(--4$\ast$d/n,muuu(d,epsomeg(d,et),d/n))[2],)));h;

\smallpagebreak

\centerline{{\bf Program 3:} refh3}

\smallpagebreak

\noindent
$\backslash\backslash$refh3(N) gives the list of invariants (d,et) of
\newline
$\backslash\backslash$main hyperbolic lattices with the square-free
\newline
$\backslash\backslash$determinant d$\backslash$le N and of
the rank 3 having the
\newline
$\backslash\backslash$number h=hnr(d,et)$\backslash$le 1 of
classes of non-reflective
\newline
$\backslash\backslash$central symmetries. Here et is a non-negative
integer
\newline
$\backslash\backslash$having the binary decomposition et=p\_k,...,p\_1
where
\newline
$\backslash\backslash$p\_1,...,p\_k are all
odd prime divisors of d in icreasing order

\smallpagebreak

\noindent
$\backslash$r h2
\newline
$\backslash$r h3

\smallpagebreak

\noindent
refh3(m,n,d,fd,sfd,sfd1,sig,n,et,h)=$\backslash$
\newline
n=0;$\backslash$
\newline
for(d=1,m,$\backslash$
\newline
if(issqfree(d)!=1,,$\backslash$
\newline
if(d$<$=2,n=n+1;h=0;et=0;pprint("n=",n," d=",d," et=",et," h=",h),$
\backslash$
\newline
fd=factor(d);sfd=matsize(fd)[1];$\backslash$
\newline
if(fd[1,1]==2,sfd1=sfd--1,sfd1=sfd);$\backslash$
\newline
for(et=0,2$^\wedge$sfd1--1,$\backslash$
\newline
if(fd[1,1]!=2,,$\backslash$
\newline
sig=mod(0,8);for(k=2,sfd,sig=sig+mod(1--fd[k,1]+4$\ast
$bittest(et,k--2),8)));$\backslash$
\newline
if(fd[1,1]==2\&\&(sig==mod(0,8)$||$sig==mod(--2,8)),$\backslash$
\newline
h=hnr(d,et);if(h$<$=1,n=n+1;$\backslash$
\newline
pprint("n=",n," d=",d," et=",et," h=",h),),$\backslash$
\newline
if(fd[1,1]!=2,h=hnr(d,et);$\backslash$
\newline
if(h$<$=1,n=n+1;pprint("n=",n," d=",d," et=",et," h=",h),),))))))

\smallpagebreak

\centerline{{\bf Program 4:} fund11.gen}

\smallpagebreak

\noindent
$\backslash$l;$\backslash$
\newline
epsilon=1;epsilon1=1;epsilon2=1;n=0;$\backslash$
\newline
for(alpha12=1,4,for(alpha23=alpha12,4,$\backslash$
\newline
u=((sqrt(2+sqrt(alpha12))+sqrt(2+sqrt(alpha23)))$^\wedge$2--
2)$^\wedge$2+0.000001;$\backslash$
\newline
for(alpha13=alpha23,u,$\backslash$
\newline
if(alpha23==0$||$issquare(alpha12$\ast
$alpha23$\ast$alpha13)!=1$||$$\backslash$
\newline
--8+2$\ast$isqrt(alpha12$\ast$alpha23$\ast$alpha13)+2$\ast
$alpha12+2$\ast$alpha13+2$\ast$alpha23$<$=0,,$\backslash$
\newline
alpha=4$\ast$idmat(3);alpha[1,2]=alpha12;alpha[2,1]=alpha12;$\backslash$
\newline
alpha[2,3]=alpha23;alpha[3,2]=alpha23;alpha[1,3]=alpha13;$\backslash$
\newline
alpha[3,1]=alpha13;dalpha=$\backslash$
\newline
--8+2$\ast$isqrt(alpha12$\ast$alpha23$\ast$alpha13)+2$\ast
$alpha12+2$\ast$alpha13+2$\ast$alpha23;$\backslash$
\newline
a=--2$\ast$idmat(3);$\backslash$
\newline
fordiv(alpha[1,2],a12,fordiv(alpha[2,3],a23,fordiv(alpha[1,3],a13,$
\backslash$
\newline
a21=alpha[1,2]/a12;a32=alpha[2,3]/a23;a31=alpha[1,3]/a13;$\backslash$
\newline
if(a12$\ast$a23$\ast$a31!=a21$\ast$a13$\ast$a32,,$\backslash$
\newline
a[1,2]=a12;a[2,1]=a21;a[2,3]=a23;a[3,2]=a32;a[1,3]=a13;a[3,1]=a31;$
\backslash$
\newline
dd=idmat(3);dd[1,1]=a13$\ast$a32;dd[2,2]=a23$\ast
$a31;dd[3,3]=a31$\ast$a32;$\backslash$
\newline
b=a$\ast$dd;b=b/content(b);n=n+1;$\backslash$
\newline
db=smith(b);r=db[1];$\backslash$
\newline
if(r$>$epsilon,epsilon=r,);$\backslash$
\newline
fr=factor(r);tfr=matsize(fr)[1];$\backslash$
\newline
r1=1;for(j=1,tfr,r1=r1$\ast$fr[j,1]);$\backslash$
\newline
if(type(r1/2)==1,r1=r1/2,);$\backslash$
\newline
if(r1$>$epsilon1,epsilon1=r1,);$\backslash$
\newline
if(tfr$<$=0,,if(fr[tfr,1]$>$epsilon2,epsilon2=fr[tfr,1],));$\backslash$
\newline
))))))));pprint("nI1=",n);pprint("aI1=",epsilon);$\backslash$
\newline
pprint("aI1\_1=",epsilon1);pprint("aI1\_2=",epsilon2);

\smallpagebreak

\centerline{{\bf Program 5:} fund10.gen}

\smallpagebreak

\noindent
$\backslash$l;$\backslash$
epsilon=1;epsilon1=1;epsilon2=1;$\backslash$
\newline
n=0;$\backslash$
\newline
alpha12=0;for(alpha23=1,4,$\backslash$
\newline
w=((sqrt(2)+sqrt(2+sqrt(alpha23)))$^\wedge$2--2)$^\wedge
$2+0.000001;$\backslash$
\newline
for(alpha13=alpha23,w,$\backslash$
\newline
if(alpha23==0$||$issquare(alpha12$\ast$alpha23$\ast
$alpha13)!=1$||$$\backslash$
\newline
--8+2$\ast$isqrt(alpha12$\ast$alpha23$\ast$alpha13)+2$\ast
$alpha12+2$\ast$alpha13+2$\ast$alpha23$<$=0,,$\backslash$
\newline
alpha=4$\ast$idmat(3);alpha[1,2]=alpha12;alpha[2,1]=alpha12;$\backslash$
\newline
alpha[2,3]=alpha23;alpha[3,2]=alpha23;alpha[1,3]=alpha13;$\backslash$
\newline
alpha[3,1]=alpha13;dalpha=$\backslash$
\newline
--8+2$\ast$isqrt(alpha12$\ast$alpha23$\ast$alpha13)+2$\ast
$alpha12+2$\ast$alpha13+2$\ast$alpha23;$\backslash$
\newline
a=--2$\ast$idmat(3);a12=0;a21=0;$\backslash$
\newline
fordiv(alpha[2,3],a23,fordiv(alpha[1,3],a13,$\backslash$
\newline
a32=alpha[2,3]/a23;a31=alpha[1,3]/a13;$\backslash$
\newline
if(a12$\ast$a23$\ast$a31!=a21$\ast$a13$\ast$a32,,$\backslash$
\newline
a[1,2]=a12;a[2,1]=a21;a[2,3]=a23;a[3,2]=a32;a[1,3]=a13;a[3,1]=a31;$\backslash$
\newline
dd=idmat(3);dd[1,1]=a13$\ast$a32;dd[2,2]=a23$\ast
$a31;dd[3,3]=a31$\ast$a32;$\backslash$
\newline
b=a$\ast$dd;b=b/content(b);n=n+1;$\backslash$
\newline
db=smith(b);r=db[1];$\backslash$
\newline
if(r$>$epsilon,epsilon=r,);$\backslash$
\newline
fr=factor(r);tfr=matsize(fr)[1];$\backslash$
\newline
r1=1;for(j=1,tfr,r1=r1$\ast$fr[j,1]);$\backslash$
\newline
if(type(r1/2)==1,r1=r1/2,);$\backslash$
\newline
if(r1$>$epsilon1,epsilon1=r1,);$\backslash$
\newline
if(tfr$<$=0,,if(fr[tfr,1]$>$epsilon2,epsilon2=fr[tfr,1],));$\backslash$
\newline
))))));pprint("nI0=",n);pprint("aI0=",epsilon);$\backslash$
\newline
pprint("aI0\_1=",epsilon1);pprint("aI0\_2=",epsilon2);

\smallpagebreak

\centerline{{\bf Program 6:} fund21.gen}

\smallpagebreak

\noindent
$\backslash$l;$\backslash$
\newline
epsilon=1;epsilon1=1;epsilon2=1;n=0;$\backslash$
\newline
alpha12=alpha23=0;$\backslash$
\newline
for(alpha34=1,4,$\backslash$
\newline
w=(4$\ast$max((sqrt(2)+sqrt(sqrt(alpha34)/4+1/2))$^\wedge$2,$\backslash$
\newline
((2+sqrt(sqrt(alpha34)/2+5/4))$^\wedge$2--1/4)/2)--2)$
^\wedge$2+0.000001;$\backslash$
\newline
for(alpha14=0,w,$\backslash$
\newline
for(alpha13=5,36,$\backslash$
\newline
if(issquare(u=alpha13$\ast$alpha34$\ast$alpha14)!=1,,$\backslash$
\newline
if(type(alpha24=4+4$\ast$(alpha14+alpha34+isqrt(u))/(alpha13--4))!=1$
||$$\backslash$
\newline
alpha24$<$=((sqrt(2)+sqrt(2+sqrt(alpha34)))$^\wedge
$2--2)$^\wedge$2--0.0000001,,$\backslash$
\newline
alpha=4$\ast$idmat(4);$\backslash$
\newline
alpha[1,3]=alpha[3,1]=alpha13;alpha[3,4]=alpha[4,3]=alpha34;$\backslash$
\newline
alpha[2,4]=alpha[4,2]=alpha24;alpha[1,4]=alpha[4,1]=alpha14;$\backslash$
\newline
fordiv(alpha[3,4],a34,$\backslash$
\newline
fordiv(alpha[1,3],a13,$\backslash$
\newline
fordiv(alpha[2,4],a24,$\backslash$
\newline
if(type(a41=isqrt(u)/(a13$\ast$a34))!=1,,$\backslash$
\newline
a=--2$\ast$idmat(4);$\backslash$
\newline
a[3,4]=a34;a[4,3]=alpha[3,4]/a34;$\backslash$
\newline
a[2,4]=a24;a[4,2]=alpha[2,4]/a24;$\backslash$
\newline
a[1,3]=a13;a[3,1]=alpha[1,3]/a13;$\backslash$
\newline
a[4,1]=a41;if(a41==0,a[1,4]=0,a[1,4]=alpha[1,4]/a41);$\backslash$
\newline
if(type(a[1,4])!=1,,$\backslash$
\newline
diag=idmat(4);diag[1,1]=a[1,3]$\ast$a[3,4]$\ast$a[4,2];$\backslash$
\newline
diag[2,2]=a[3,1]$\ast$a[4,3]$\ast$a[2,4];diag[3,3]=a[3,1]$\ast
$a[3,4]$\ast$a[4,2];$\backslash$
\newline
diag[4,4]=a[3,1]$\ast$a[4,3]$\ast$a[4,2];b=a$\ast
$diag;b=b/content(b);n=n+1;$\backslash$
\newline
db=smith(b);r=db[2];$\backslash$
\newline
if(r$>$epsilon,epsilon=r,);$\backslash$
\newline
fr=factor(r);tfr=matsize(fr)[1];$\backslash$
\newline
r1=1;for(j=1,tfr,r1=r1$\ast$fr[j,1]);$\backslash$
\newline
if(type(r1/2)==1,r1=r1/2,);$\backslash$
\newline
if(r1$>$epsilon1,epsilon1=r1,);$\backslash$
\newline
if(tfr$<$=0,,if(fr[tfr,1]$>$epsilon2,epsilon2=fr[tfr,1],));$\backslash$
\newline
))))))))));pprint("nII1=",n);pprint("aII1=",epsilon);$\backslash$
\newline
pprint("aII1\_1=",epsilon1);pprint("aII1\_2=",epsilon2);

\smallpagebreak

\centerline{{\bf Program 7:} fund20.gen}

\smallpagebreak

\noindent
$\backslash$l;$\backslash$
\newline
epsilon=1;epsilon1=1;epsilon2=1;n=0;$\backslash$
\newline
alpha12=0;alpha23=0;alpha34=0;$\backslash$
\newline
for(alpha14=1,287.10,$\backslash$
\newline
fordiv(4$\ast$alpha14,aa,$\backslash$
\newline
if(aa$^\wedge$2$>$4$\ast$alpha14,,alpha13=4+aa;alpha24=4$\ast
$alpha14/aa+4;$\backslash$
\newline
if(alpha13$>$36,,$\backslash$
\newline
alpha=4$\ast$idmat(4);$\backslash$
\newline
alpha[1,3]=alpha13;alpha[3,1]=alpha13;alpha[3,4]=alpha34;$\backslash$
\newline
alpha[4,3]=alpha34;alpha[2,4]=alpha24;alpha[4,2]=alpha24;$\backslash$
\newline
alpha[1,4]=alpha14;alpha[4,1]=alpha14;$\backslash$
\newline
fordiv(alpha[1,3],a13,fordiv(alpha[2,4],a24,fordiv(alpha[1,4],a14,$
\backslash$
\newline
a=--2$\ast$idmat(4);a[1,3]=a13;a[3,1]=alpha[1,3]/a13;a[1,4]=a14;$
\backslash$
\newline
a[4,1]=alpha[1,4]/a14;a[1,3]=a13;a[3,1]=alpha[1,3]/a13;$\backslash$
\newline
a[2,4]=a24;a[4,2]=alpha[2,4]/a24;$\backslash$
\newline
diag=idmat(4);diag[1,1]=a[1,3]$\ast$a[1,4]$\ast$a[4,2];$\backslash$
\newline
diag[2,2]=a[1,3]$\ast$a[4,1]$\ast$a[2,4];diag[3,3]=a[3,1]$\ast
$a[1,4]$\ast$a[4,2];$\backslash$
\newline
diag[4,4]=a[1,3]$\ast$a[4,2]$\ast
$a[4,1];b=a$\ast$diag;b=b/content(b);n=n+1;$\backslash$
\newline
db=smith(b);r=db[2];$\backslash$
\newline
if(r$>$epsilon,epsilon=r,);$\backslash$
\newline
fr=factor(r);tfr=matsize(fr)[1];$\backslash$
\newline
r1=1;for(j=1,tfr,r1=r1$\ast$fr[j,1]);$\backslash$
\newline
if(type(r1/2)==1,r1=r1/2,);$\backslash$
\newline
if(r1$>$epsilon1,epsilon1=r1,);$\backslash$
\newline
if(tfr$<$=0,,if(fr[tfr,1]$>$epsilon2,epsilon2=fr[tfr,1],));$\backslash$
\newline
)))))));pprint("nII0=",n);pprint("aII0=",epsilon);$\backslash$
\newline
pprint("aII0\_1=",epsilon1);pprint("aII0\_2=",epsilon2);

\smallpagebreak

\centerline{{\bf Program 8:} fund30.gen}

\smallpagebreak

\noindent
$\backslash$l;$\backslash$
\newline
epsilon=1;epsilon1=1;epsilon2=1;n=0;$\backslash$
\newline
al12=al23=al34=al45=0;$\backslash$
\newline
for(al15=0,900,for(al13=5,36,for(al35=al13,36,$\backslash$
\newline
if(issquare(q=al13$\ast$al35$\ast
$al15)==0,,d=(al13+al35+al15--4+isqrt(q))$\ast$4;$\backslash$
\newline
if(type(al14=d/(al35--4))!=1$||$al14$<$=287.108350$||
$type(al25=d/(al13--4))!=1$||$$\backslash$
\newline
al25$<$=287.108350$||$$\backslash$
\newline
type(al24=(al13$\ast$al35+4$\ast$al15+4$\ast
$isqrt(q))$\ast$4/((al35--4)$\ast$(al13--4)))!=1$||$$\backslash$
\newline
issquare(q1=al13$\ast$al35$\ast$al25$\ast$al24$\ast$al14)==0,,$\backslash$
\newline
al=idmat(5)$\ast$4;al[1,5]=al[5,1]=al15;al[1,3]=al[3,1]=al13;$\backslash$
\newline
al[1,4]=al[4,1]=al14;al[2,4]=al[4,2]=al24;al[2,5]=al[5,2]=al25;$\backslash$
\newline
al[3,5]=al[5,3]=al35;$\backslash$
\newline
fordiv(al13,a13,$\backslash$
\newline
fordiv(al35,a35,a51=isqrt(q)/a13/a35;$\backslash$
\newline
if(a51==0,a15==0,if(type(a15=al15/a51)!=1,,$\backslash$
\newline
fordiv(al14,a14,$\backslash$
\newline
fordiv(al24,a24,a52=isqrt(q1)/a13/a35/a24$\ast$a14/al14;$\backslash$
\newline
if(type(a25=al25/a52)!=1,,$\backslash$
\newline
a=idmat(5)$\ast$--2;$\backslash$
\newline
a[1,3]=a13;a[3,1]=al13/a13;$\backslash$
\newline
a[1,4]=a14;a[4,1]=al14/a14;$\backslash$
\newline
a[1,5]=a15;a[5,1]=a51;$\backslash$
\newline
a[2,4]=a24;a[4,2]=al24/a24;$\backslash$
\newline
a[2,5]=a25;a[5,2]=a52;$\backslash$
\newline
a[3,5]=a35;a[5,3]=al35/a35;$\backslash$
\newline
diag=idmat(5);diag[1,1]=a[1,4]$\ast$a[1,3]$\ast
$a[4,2]$\ast$a[2,5];$\backslash$
\newline
diag[2,2]=a[4,1]$\ast$a[1,3]$\ast$a[2,4]$\ast$a[2,5];$\backslash$
\newline
diag[3,3]=a[1,4]$\ast$a[3,1]$\ast$a[4,2]$\ast$a[2,5];$\backslash$
\newline
diag[4,4]=a[4,1]$\ast$a[1,3]$\ast$a[4,2]$\ast$a[2,5];$\backslash$
\newline
diag[5,5]=a[4,1]$\ast$a[1,3]$\ast$a[2,4]$\ast$a[5,2];$\backslash$
\newline
b=a$\ast$diag;b=b/content(b);n=n+1;$\backslash$
\newline
db=smith(b);r=db[3];$\backslash$
\newline
if(r$>$epsilon,epsilon=r,);$\backslash$
\newline
fr=factor(r);tfr=matsize(fr)[1];$\backslash$
\newline
r1=1;for(j=1,tfr,r1=r1$\ast$fr[j,1]);$\backslash$
\newline
if(type(r1/2)==1,r1=r1/2,);$\backslash$
\newline
if(r1$>$epsilon1,epsilon1=r1,);$\backslash$
\newline
if(fr[tfr,1]$>$epsilon2,epsilon2=fr[tfr,1],);$\backslash$
\newline
))))))))))));pprint("nIII0=",n);$\backslash$
\newline
pprint("aIII0=",epsilon);pprint("aIII0\_1=",epsilon1);$\backslash$
\newline
pprint("aIII0\_2=",epsilon2);

\smallpagebreak

\centerline{{\bf Program 9:} fund20.main}

\smallpagebreak

\noindent
$\backslash$r h2
\newline
$\backslash$r h3
\newline
$\backslash$l;$\backslash$
\newline
n=0;$\backslash$
\newline
alpha12=0;alpha23=0;alpha34=0;$\backslash$
\newline
for(alpha14=1,287.10,$\backslash$
\newline
fordiv(4$\ast$alpha14,aa,$\backslash$
\newline
if(aa$^\wedge$2$>$4$\ast$alpha14,,alpha13=4+aa;alpha24=4$\ast
$alpha14/aa+4;$\backslash$
\newline
alpha=4$\ast$idmat(4);$\backslash$
\newline
alpha[1,3]=alpha13;alpha[3,1]=alpha13;alpha[3,4]=alpha34;$\backslash$
\newline
alpha[4,3]=alpha34;alpha[2,4]=alpha24;alpha[4,2]=alpha24;$\backslash$
\newline
alpha[1,4]=alpha14;alpha[4,1]=alpha14;$\backslash$
\newline
fordiv(alpha[1,3],a13,fordiv(alpha[2,4],a24,fordiv(alpha[1,4],a14,$
\backslash$
\newline
a=--2$\ast$idmat(4);a[1,3]=a13;a[3,1]=alpha[1,3]/a13;a[1,4]=a14;$\backslash$
\newline
a[4,1]=alpha[1,4]/a14;a[1,3]=a13;a[3,1]=alpha[1,3]/a13;$\backslash$
\newline
a[2,4]=a24;a[4,2]=alpha[2,4]/a24;$\backslash$
\newline
diag=idmat(4);diag[1,1]=a[1,3]$\ast$a[1,4]$\ast$a[4,2];$\backslash$
\newline
diag[2,2]=a[1,3]$\ast$a[4,1]$\ast$a[2,4];diag[3,3]=a[3,1]$\ast
$a[1,4]$\ast$a[4,2];$\backslash$
\newline
diag[4,4]=a[1,3]$\ast$a[4,2]$\ast$a[4,1];b=a$\ast
$diag;b=b/content(b);$\backslash$
\newline
m=m+1;db=smith(b);dbb=db[2];fdbb=factor(dbb);$\backslash$
\newline
if(issqfree(b[1,1])==0$||$issqfree(b[2,2])==0$||
$issqfree(b[3,3])==0$||$$\backslash$
\newline
issqfree(b[4,4])==0$||$content([--b[1,1],--b[2,2]])$>$2$||$$\backslash$
\newline
content([--b[2,2],--b[3,3]])$>$2$||
$content([--b[3,3],--b[4,4]])$>$2,,$\backslash$
\newline
detb=db[2]$\ast$db[3]$\ast$db[4];$\backslash$
\newline
if(detb==1,,$\backslash$
\newline
fdetb=factor(detb);$\backslash$
\newline
if(content([--b[1,1]$\ast
$--b[2,2]$\ast$--b[3,3]$\ast$--b[4,4],16])$>$1\&\&$\backslash$
\newline
content([--b[1,1]$\ast$--b[2,2]$\ast
$--b[3,3]$\ast$--b[4,4],16])$<$16\&\&$\backslash$
\newline
fdetb[1,1]==2\&\&mod(fdetb[1,2],2)==mod(1,2),,$\backslash$
\newline
gam=0;$\backslash$
\newline
for(j=1,matsize(fdetb)[1],$\backslash$
\newline
if(fdetb[j,1]!=2\&\&$\backslash$
\newline
type(--b[1,1]$\ast$--b[2,2]$\ast$--b[3,3]$\ast
$--b[4,4]/fdetb[j,1])==1\&\&$\backslash$
\newline
mod(fdetb[j,2],2)==mod(0,2),gam=1,));$\backslash$
\newline
if(gam==1,,$\backslash$
\newline
for(k=0,1,$\backslash$
\newline
if(k==0\&\&content(--[b[1,1]$\ast$--b[2,2]$\ast
$--b[3,3]$\ast$--b[4,4],2])==1\&\&$\backslash$
\newline
fdetb[1,1]==2\&\&mod(fdetb[1,2],2)==mod(1,2),,$\backslash$
\newline
if(k==1\&\&content(--[b[1,1]$\ast$--b[2,2]$\ast$--b[3,3]$\ast
$--b[4,4],2])==1,$\backslash$
\newline
b=2$\ast$b;db=smith(b);dbb=db[2];fdbb=factor(dbb);gam1=0,gam1=1);$
\backslash$
\newline
if(k==1\&\&gam1=1,,$\backslash$
\newline
detb=db[2]$\ast$db[3]$\ast$db[4];fdetb=factor(detb);$\backslash$
\newline
d=1;for(k=1,matsize(fdetb)[1],$\backslash$
\newline
if(mod(fdetb[k,2],2)==mod(1,2),d=d$\ast$fdetb[k,1],));$\backslash$
\newline
if(d$<$=2,et=0,$\backslash$
\newline
fd=factor(d);if(fd[1,1]==2,d1=d/2,d1=d);fd1=factor(d1);$\backslash$
\newline
et=0;for(k=1,matsize(fd1)[1],$\backslash$
\newline
if(type(b[1,1]/fd1[k,1])==1,$\backslash$
\newline
if(kro(b[1,1]/fd1[k,1],fd1[k,1])==1,,et=et+2$^\wedge$(k--1)),$\backslash$
\newline
if(type(b[2,2]/fd1[k,1])==1,$\backslash$
\newline
if(kro(b[2,2]/fd1[k,1],fd1[k,1])==1,,et=et+2$^\wedge$(k--1)),$\backslash$
\newline
if(kro(d$\ast$b[1,1]$\ast$b[2,2]/fd1[k,1],fd1[k,1])==1,,et=et+2$
^\wedge$(k--1))))));$\backslash$
\newline
hhh=hnr(d,et);if(hhh$>$1,,$\backslash$
\newline
n=n+1;pprint("n=",n);pprint(a);pprint(b);pprint(db);$\backslash$
\newline
pprint(fdbb);pprint("d=",d," et=",et," h=",hhh);$\backslash$
\newline
if(n==1,u=matrix(1,3,j,k,0);u[1,]=[d,et,hhh],$\backslash$
\newline
aaa=0;nnn=matsize(u)[1];$\backslash$
\newline
for(j=1,nnn,if(u[j,]==[d,et,hhh],aaa=1,));$\backslash$
\newline
if(aaa==1,,u1=matrix(nnn+1,3,j,k,0);$\backslash$
\newline
for(t=1,nnn,$\backslash$
\newline
if(u[t,1]$<$d$||$(u[t,1]==d\&\&u[t,2]$<$et),u1[t,]=u[t,],$\backslash$
\newline
u1[t+1,]=u[t,]));$\backslash$
\newline
for(t=1,nnn+1,if(u1[t,]==[0,0,0],u1[t,]=[d,et,hhh],));u=u1));$\backslash$
\newline
))))))))))))));pprint("u=",u);pprint("matsize u=",matsize(u));

\smallpagebreak

\centerline{{\bf Program 10:} refl0.1}

\smallpagebreak

\noindent
$\backslash$$\backslash$main(n,h) calculates for U+$<$--2n$>$ (it is
given by the matrix g)
\newline
$\backslash$$\backslash$the vectors v of the height $\backslash$le h
from P($\backslash$M)\_\{pr\}
\newline
$\backslash$$\backslash$and calculates chains e and f and
their Gram matrices ge, gf
\newline
refl(n,h,m,m1,j,k,v1,y1,y2,y3,z)=$\backslash$
\newline
g=[0,1,0;1,0,0;0,0,--2$\ast$n];$\backslash$
\newline
m=3;v=matrix(m,3,j,k,0);v[1,]=[0,0,1];v[2,]=[n,0,--1];v[3,]=[--1,1,0];$
\backslash$
\newline
for(h1=2,h,$\backslash$
\newline
fordiv(h1,y2,if(type(y2$^\wedge$2/h1)!=1$||$type(n$\ast$h1/y2$^\wedge
$2)!=1,,$\backslash$
\newline
d=--2$\ast$y2$^\wedge$2/h1;$\backslash$
\newline
for(z=floor(sqrt((2$\ast$y2$^\wedge$2--d)/(2$\ast$n))),y2,$\backslash$
\newline
y3=--z;y1=(2$\ast$n$\ast$y3$^\wedge$2+d)/(2$\ast$y2);$\backslash$
\newline
if(type(2$\ast$y1/d)!=1$||$content([y1,y2,y3])!=1,,$\backslash$
\newline
u=[y1,y2,y3];$\backslash$
\newline
alpha=1;m1=1;$\backslash$
\newline
while(alpha,if(m1$>$m,alpha=0,if(v[m1,]$\ast$g$\ast$u$\widetilde{\ }$$>
$=0,m1=m1+1,alpha=0));$\backslash$
\newline
if(m1$<$=m,,m=m1;v1=matrix(m,3,j,k,0);$\backslash$
\newline
for(j=1,m--1,v1[j,]=v[j,]);v1[m,]=[y1,y2,y3];v=v1;kill(v1)$\backslash$
\newline
)))))));v;

\noindent
$\backslash$$\backslash$$\ast$$\ast$

\noindent
$\backslash$$\backslash$getting e1--matrix from v
\newline
e1fromv(v,s,s1,alpha,s2,ex)=$\backslash$
\newline
s=matsize(v)[1];$\backslash$
\newline
e1=matrix(s,3,j,k,0);e1[1,]=v[2,];e1[2,]=v[1,];$\backslash$
\newline
alpha=1;s1=2;$\backslash$
\newline
while(alpha,$\backslash$
\newline
if(s1==s,alpha=0,$\backslash$
\newline
for(j=1,s,$\backslash$
\newline
if(v[j,]==e1[s1--1,]$||$v[j,]==e1[s1,]$||$$\backslash$
\newline
(v[j,]$\ast$g$\ast$e1[s1,]$\widetilde{\ }$)$^\wedge$2/((v[j,]$\ast
$g$\ast$v[j,]$\widetilde{\ }$)$\ast$(e1[s1,]$\ast
$g$\ast$e1[s1,]$\widetilde{\ }$))$>$1,,$\backslash$
\newline
s2=s1+1;e1[s2,]=v[j,]));if(s2$>$s1,s1=s2,alpha=0)));$\backslash$
\newline
ex=matrix(s1,3,j,k,e1[j,k]);e1=ex;

\noindent
$\backslash$$\backslash$getting e2--matrix from v
\newline
e2fromv(v,s,s1,alpha,s2,ex)=$\backslash$
\newline
s=matsize(v)[1];$\backslash$
\newline
e2=matrix(s,3,j,k,0);e2[1,]=v[1,];e2[2,]=v[2,];$\backslash$
\newline
alpha=1;s1=2;$\backslash$
\newline
while(alpha,$\backslash$
\newline
if(s1==s,alpha=0,$\backslash$
\newline
for(j=1,s,$\backslash$
\newline
if(v[j,]==e2[s1--1,]$||$v[j,]==e2[s1,]$||$$\backslash$
\newline
(v[j,]$\ast$g$\ast$e2[s1,]$\widetilde{\ }$)$^\wedge
$2/((v[j,]$\ast$g$\ast$v[j,]$\widetilde{\ }$)$\ast
$(e2[s1,]$\ast$g$\ast$e2[s1,]$\widetilde{\ }$))$>$1,,$\backslash$
\newline
s2=s1+1;e2[s2,]=v[j,]));if(s2$>$s1,s1=s2,alpha=0)));$\backslash$
\newline
ex=matrix(s1,3,j,k,e2[j,k]);e2=ex;

\noindent
$\backslash$$\backslash$getting e--matrix from v
\newline
efromv(v,e1,e2,s,s1,s2)=$\backslash$
\newline
s=matsize(v)[1];s1=matsize(e1)[1];s2=matsize(e2)[1];$\backslash$
\newline
if(s1==s\&\&s2==s,e=e1,$\backslash$
\newline
e=matrix(s1+s2--2,3,j,k,0);$\backslash$
\newline
for(j=1,s2,e[j,]=e2[s2+1--j,]);$\backslash$
\newline
for(j=1,s1--2,e[s2+j,]=e1[2+j,]));$\backslash$
\newline
e;ge=e$\ast$g$\ast$e$\widetilde{\ }$;

\noindent
$\backslash$$\backslash$getting f1--matrix from v
\newline
f1fromv(v,e,s,s1,s2,q,alpha,f1x)=$\backslash$
\newline
s=matsize(v)[1];s1=matsize(e)[1];$\backslash$
\newline
if(s==s1,f1=0;v1=0,$\backslash$
\newline
f1=matrix(s--s1,3,j,k,0);v1=v;$\backslash$
\newline
for(j=1,s,$\backslash$
\newline
for(k=1,s1,$\backslash$
\newline
if(v1[j,]!=e[k,],,v1[j,]=[0,0,0])));$\backslash$
\newline
q=1;s2=0;alpha=1;$\backslash$
\newline
while(alpha,$\backslash$
\newline
if(q$>$s,alpha=0,$\backslash$
\newline
if(v1[q,]==[0,0,0],q=q+1,s2=s2+1;f1[s2,]=v1[q,];$\backslash$
\newline
v1[q,]=[0,0,0];alpha=0)));$\backslash$
\newline
if(s2==s--s1,,$\backslash$
\newline
for(t=1,s,$\backslash$
\newline
for(q=1,s,$\backslash$
\newline
if(v1[q,]==[0,0,0],,$\backslash$
\newline
if((v1[q,]$\ast$g$\ast$f1[s2,]$\widetilde{\ }$)$^\wedge
$2/((v1[q,]$\ast$g$\ast$v1[q,]$\widetilde{\ }$)$\ast
$(f1[s2,]$\ast$g$\ast$f1[s2,]$\widetilde{\ }$))$>$1,,$\backslash$
\newline
s2=s2+1;f1[s2,]=v1[q,];v1[q,]=[0,0,0])))));$\backslash$
\newline
if(s2==0,,f1x=matrix(s2,3,j,k,f1[j,k]);f1=f1x);f1);

\noindent
$\backslash$$\backslash$getting f2--matrix from v
\newline
f2fromv(f1,s,s1,s2,s3,q,f2x)=$\backslash$
\newline
if(f1==0,f2=0;s3=0,$\backslash$
\newline
s=matsize(v1)[1];s1=matsize(e)[1];s2=matsize(f1)[1];$\backslash$
\newline
f2=matrix(s--s1--s2+1,3,j,k,0);$\backslash$
\newline
s3=1;f2[1,]=f1[1,];$\backslash$
\newline
for(t=1,s,$\backslash$
\newline
for(q=1,s,$\backslash$
\newline
if(v1[q,]==[0,0,0],,$\backslash$
\newline
if((v1[q,]$\ast$g$\ast$f2[s3,]$\widetilde{\ }$)$^\wedge
$2/((v1[q,]$\ast$g$\ast$v1[q,]$\widetilde{\ }$)$\ast
$(f2[s3,]$\ast$g$\ast$f2[s3,]$\widetilde{\ }$))$>$1,,$\backslash$
\newline
s3=s3+1;f2[s3,]=v1[q,];v1[q,]=[0,0,0]))));$\backslash$
\newline
if(s3==0,,f2x=matrix(s3,3,j,k,f2[j,k]);f2=f2x));f2;

\noindent
$\backslash$$\backslash$getting f--matrix from v
\newline
ffromv(f1,f2,s1,s2)=$\backslash$
\newline
if(f2==0,f=f1,$\backslash$
\newline
s1=matsize(f1)[1];s2=matsize(f2)[1];$\backslash$
\newline
if(s2==0,f=f1,$\backslash$
\newline
f=matrix(s1+s2--1,3,j,k,0);$\backslash$
\newline
for(j=1,s2,f[j,]=f2[s2+1--j,]);$\backslash$
\newline
for(j=1,s1--1,f[s2+j,]=f1[1+j,])));$\backslash$
\newline
if(f==0,gf=0,gf=f$\ast$g$\ast$f$\widetilde{\ }$);

\noindent
$\backslash$$\backslash$$\ast$$\ast$$\ast$

\noindent
$\backslash$l;$\backslash$
\newline
main(n,h)=refl(n,h);e1fromv(v);e2fromv(v,e1);efromv(v,e1,e2);$\backslash$
\newline
f1fromv(v,e);f2fromv(f1);ffromv(f1,f2);

\smallpagebreak

\centerline{{\bf Program 11:} refl0.12}

\smallpagebreak

\noindent
$\backslash$$\backslash$main(n1,n2,n3,h) calculates for $<$n1$>
$$\backslash$oplus $<$--n2$>$$\backslash$oplus $<$--n3$>$
\newline
$\backslash$$\backslash$(it is given by the matrix g)
\newline
$\backslash$$\backslash$the vectors v of the height $\backslash
$le h from P($\backslash$M)\_\{pr\}
\newline
$\backslash$$\backslash$and calculates chains e and f and
their Gram matrices ge, gf
\newline
refl(n1,n2,n3,h,m,m1,j,k,v1,y1,y2,y3,z,d,dd,h1,u,u1,w,w1)=$\backslash$
\newline
g=[n1,0,0;0,--n2,0;0,0,--n3];$\backslash$
\newline
m=2;v=matrix(m,3,j,k,0);v[1,]=[0,1,0];$\backslash$
\newline
if(n2==n3,v[2,]=[0,--1,1],v[2,]=[0,0,1]);$\backslash$
\newline
for(h1=1,h,$\backslash$
\newline
fordiv(h1,y1,dd=gcd(2$\ast$lcm(lcm(n1,n2),n3),2$\ast$n1$\ast$y1);$\backslash$
\newline
fordiv(dd,d,if(h1!=2$\ast$n1$\ast$y1$^\wedge$2/d,,$\backslash$
\newline
for(z=0,floor(sqrt((w=n1$\ast$y1$^\wedge$2+d)/n2)+0.000001),$\backslash$
\newline
if(type(2$\ast$n2$\ast$z/d)!=1$||$type(w1=(w--n2$\ast$z$^\wedge
$2)/n3)!=1,,$\backslash$
\newline
if(issquare(w1)!=1,,$\backslash$
\newline
y2=--z;y3=--isqrt(w1);if(type(2$\ast$n3$\ast$y3/d)!=1$||
$content([y1,y2,y3])!=1,,$\backslash$
\newline
u=[y1,y2,y3];$\backslash$
\newline
alpha=1;m1=1;$\backslash$
\newline
while(alpha,if(m1$>$m,alpha=0,if(v[m1,]$\ast$g$\ast$u$\widetilde{\ }$$>
$=0,m1=m1+1,alpha=0));$\backslash$
\newline
if(m1$<$=m,,m=m1;v1=matrix(m,3,j,k,0);$\backslash$
\newline
for(j=1,m--1,v1[j,]=v[j,]);v1[m,]=[y1,y2,y3];v=v1;kill(v1)$\backslash$
\newline
))))))))));v;

\noindent
$\backslash$$\backslash$$\ast$$\ast$
\newline
This part is the same as in Program 10: refl0.1 between $\backslash
$$\backslash$$\ast$$\ast$ and $\backslash$$\backslash$$\ast$$\ast$$\ast$

\noindent
$\backslash$$\backslash$$\ast$$\ast$$\ast$

\noindent
$\backslash$l;$\backslash$
\newline
main(n1,n2,n3,h)=$\backslash$
\newline
refl(n1,n2,n3,h);e1fromv(v);e2fromv(v,e1);efromv(v,e1,e2);$\backslash$
\newline
f1fromv(v,e);f2fromv(f1);ffromv(f1,f2);

\smallpagebreak

\centerline{{\bf Program 12:} refl0.13}

\smallpagebreak

\noindent
$\backslash$$\backslash$main(n1,n2,n3,eps1,eps2,eps3,h) calculates for
\newline
$\backslash$$\backslash$$<$n1$>$$\backslash$oplus $<$--n2$>
$$\backslash$oplus $<$--n3$>$(eps1/2,eps2/2,eps3/2)
\newline
$\backslash$$\backslash$(it is given by the matrix g)
where ni$\backslash$ge 0 and epsi=0 or 1 and
\newline
$\backslash$$\backslash$either all of them are 0 or at least
two of them are not zero,
\newline
$\backslash$$\backslash$the vectors v of the height $\backslash$le
h from P($\backslash$M)\_\{pr\}
\newline
$\backslash$$\backslash$and calculates chains e and f and
their Gram matrices ge, gf
\newline
refl(n1,n2,n3,eps1,eps2,eps3,h,m,m1,j,k,v1,y1,y2,y3,y1t,y2t,y3t,$
\backslash$
\newline
z,d,dt,dd,h1,u,u1,w,w1)=$\backslash$
\newline
g=[n1,0,0;0,--n2,0;0,0,--n3];eps=[eps1,eps2,eps3]/2;$\backslash$
\newline
m=2;v=matrix(m,3,j,k,0);$\backslash$
\newline
if(n2/n3!=1\&\&n2/n3!=3\&\&n2/n3!=1/3,v[1,]=[0,1,0];v[2,]=[0,0,1],$
\backslash$
\newline
if(n2==n3,if(eps==[0,0,0]$||
$eps==[1/2,1/2,1/2],v[1,]=[0,1,0];v[2,]=[0,--1,1],$\backslash$
\newline
if(eps==[0,1/2,1/2],v[1,]=[0,1,0];v[2,]=[0,--1/2,1/2],$\backslash$
\newline
v[1,]=[0,1,0];v[2,]=[0,0,1])),$\backslash$
\newline
if(n2/n3==3,if(eps==[0,1/2,1/2],v[1,]=[0,1,0];v[2,]=[0,--1/2,1/2],$
\backslash$
\newline
v[1,]=[0,1,0];v[2,]=[0,0,1]),$\backslash$
\newline
if(n2/n3==1/3,if(eps==[0,1/2,1/2],v[1,]=[0,0,1];v[2,]=[0,1/2,--1/2],$
\backslash$
\newline
v[1,]=[0,1,0];v[2,]=[0,0,1]),))));$\backslash$
\newline
for(h1=1,h,$\backslash$
\newline
fordiv(h1,y1t,dd=gcd(2$\ast$lcm(lcm(n1,n2),n3),n1$\ast$y1t);$\backslash$
\newline
fordiv(dd,d,if(d$\ast$h1!=n1$\ast$y1t$^\wedge$2,,dt=4$\ast$d;$\backslash$
\newline
for(z=0,floor(sqrt((w=n1$\ast$y1t$^\wedge$2+dt)/n2)+0.000001),$\backslash$
\newline
if(type(n2$\ast$z/d)!=1$||$type(w1=(w--n2$\ast$z$^\wedge
$2)/n3)!=1,,$\backslash$
\newline
if(issquare(w1)!=1,,$\backslash$
\newline
y2t=--z;y3t=--isqrt(w1);if(type(n3$\ast$y3t/d)!=1,,$\backslash$
\newline
y1=y1t/2;y2=y2t/2;y3=y3t/2;u=[y1,y2,y3];$\backslash$
\newline
if(type(2$\ast$u$\ast$g$\ast$eps$\widetilde{\ }$/d)!=1,,$\backslash$
\newline
if(mod(2$\ast$u,2)!=mod([0,0,0],2)\&\&mod(2$\ast$u,2)!=mod(2$\ast
$eps,2),,$\backslash$
\newline
if((mod(2$\ast$u,2)==mod([0,0,0],2)\&\&(content(u)$>$1$||$$\backslash$
\newline
mod(2$\ast$u,4)==mod(4$\ast$eps,4)))$||$$\backslash$
\newline
(eps!=[0,0,0]\&\&mod(2$\ast$u,2)==mod(2$\ast$eps,2)\&\&content(2$\ast
$u)$>$1),,$\backslash$
\newline
alpha=1;m1=1;$\backslash$
\newline
while(alpha,if(m1$>$m,alpha=0,if(v[m1,]$\ast$g$\ast$u$\widetilde{\ }$$>
$=0,m1=m1+1,alpha=0));$\backslash$
\newline
if(m1$<$=m,,m=m1;v1=matrix(m,3,j,k,0);$\backslash$
\newline
for(j=1,m--1,v1[j,]=v[j,]);v1[m,]=[y1,y2,y3];v=v1;kill(v1)$\backslash$
\newline
)))))))))))));v;

\noindent
$\backslash$$\backslash$$\ast$$\ast$
\newline
This part is the same as in Program 10: refl0.1 between $\backslash
$$\backslash$$\ast$$\ast$ and $\backslash$$\backslash$$\ast$$\ast$$\ast$

\noindent
$\backslash$$\backslash$$\ast$$\ast$$\ast$

\noindent
$\backslash$l;$\backslash$
\newline
main(n1,n2,n3,eps1,eps2,eps3,h)=$\backslash$
\newline
refl(n1,n2,n3,eps1,eps2,eps3,h);$\backslash$
\newline
e1fromv(v);e2fromv(v,e1);efromv(v,e1,e2);$\backslash$
\newline
f1fromv(v,e);f2fromv(f1);ffromv(f1,f2);


\Refs
\widestnumber\key{vedG2}

\ref
\key AN1
\by V.A. Alexeev and V.V. Nikulin
\paper The classification of Del Pezzo surfaces with log terminal
singularities of the index $\le 2$, involutions of K3 surfaces
and reflection groups in Lobachevsky spaces (Russian)
\jour Doklady po matematike i prilogeniyam, MIAN
\vol 2 \issue 2 \yr 1988 \pages 51--150
\endref

\ref
\key AN2
\by V.A. Alexeev and V.V. Nikulin
\paper The classification of Del Pezzo surfaces with log terminal
singularities of the index $\le 2$ and involutions of K3 surfaces
\jour Dokl. AN SSSR \vol 306 \issue 3 \yr 1989 \pages 525--528
\transl\nofrills English transl. in
Soviet Math. Dokl. \yr 1989 \vol 39
\endref

\ref
\key B1
\by R. Borcherds
\paper Generalized Kac--Moody algebras
\jour J. of Algebra
\vol 115
\yr 1988
\pages 501--512
\endref

\ref
\key B2
\by R. Borcherds
\paper The monster Lie algebra
\jour Adv. Math.
\vol 83
\yr 1990
\pages 30--47
\endref

\ref
\key B3
\by R. Borcherds
\paper The monstrous moonshine and monstrous Lie superalgebras
\jour Invent. Math.
\vol 109
\yr 1992
\pages 405--444
\endref

\ref
\key B4
\by R. Borcherds
\paper Sporadic groups and string theory
\inbook Proc. European Congress of Mathem. 1992
\pages 411--421
\endref

\ref
\key B5
\by R. Borcherds
\paper Automorphic forms on $O_{s+2,2}$ and
infinite products
\jour Invent. Math. \vol 120
\yr 1995
\pages 161--213
\endref

\ref
\key B6
\by R. Borcherds
\paper The moduli space of Enriques surfaces and the fake monster Lie
superalgebra
\jour Topology
\yr 1996
\vol 35 \issue 3
\pages 699--710
\endref

\ref
\key B-Sh
\by Z.I. Borevich and I.R. Shafarevich
\book Number Theory
\publ Nauka
\publaddr Moscow
\yr 1985
\endref

\ref
\key C
\by J.W.S. Cassels
\book Rational quadratic forms
\publ Academic Press
\yr 1978
\endref

\ref
\key CCL
\by G.L. Cardoso, G. Curio and D. L\"ust
\paper Perturbative coupling and modular forms in $N=2$ string
models with a Wilson line
\jour Preprint
\yr 1996 \moreref hep-th/9608154
\endref

\ref
\key E1
\by M. Eichler
\book Quadratische Formen und orthogonale Gruppen
\publ Springer-Verlag
\yr 1952
\endref

\ref
\key E2
\by F. Esselmann
\paper \"Uber die maximale Dimension von Lorentz-Gittern
mit coendlicher Spiege\-lungsgruppe
\jour Preprint Univ. Bielefeld
\vol 92-023
\yr 1992
\endref

\ref
\key GN1
\by V.A. Gritsenko, V.V. Nikulin
\paper Siegel automorphic form correction of some Lorentzi\-an
Kac--Moody Lie algebras
\jour Amer. J. Math.
\yr 1997 \vol 119 \issue 1 \pages 181--224
\moreref alg-geom/9504006
\endref

\ref
\key GN2
\by V.A. Gritsenko, V.V. Nikulin
\paper Siegel automorphic form correction of a Lorentzian
Kac--Moody algebra
\jour C. R. Acad. Sci. Paris S\'er. A--B
\vol 321
\yr 1995
\pages 1151--1156
\endref

\ref
\key GN3
\by V.A. Gritsenko, V.V. Nikulin
\paper K3 surfaces, Lorentzian Kac--Moody algebras and
mirror symmetry
\jour  Math. Res. Lett. \yr 1996 \vol 3 \issue 2 \pages 211--229
\moreref  alg-geom/9510008
 \endref

\ref
\key GN4
\by V.A. Gritsenko, V.V. Nikulin
\paper The Igusa modular forms and ``the simplest''
Lorentzian Kac--Moody algebras
\jour Matem. Sbornik
\yr 1996 \vol 187 \issue 11  \pages 27--66
\transl\nofrills English transl. in
\jour Sbornik: Mathematics \vol 187
\yr 1996 \issue 11 \pages 1601--1641
\moreref alg-geom/9603010
\endref

\ref
\key GN5
\by V.A. Gritsenko, V.V. Nikulin
\paper Automorphic forms and Lorentzian Kac-Moody algebras.
Part I
\jour Preprint RIMS Kyoto Univ. \yr 1996
\vol RIMS-1116
\moreref alg-geom/9610022
\endref

\ref
\key GN6
\by V.A. Gritsenko, V.V. Nikulin
\paper Automorphic forms and Lorentzian Kac-Moody algebras.
Part II
\jour Preprint RIMS Kyoto Univ.
\yr 1996
\vol RIMS-1122
\moreref alg-geom/9611028
\endref

\ref
\key GN7
\by V.A. Gritsenko, V.V. Nikulin
\paper The arithmetic mirror symmetry and Calabi--Yau manifolds
\jour Preprint RIMS Kyoto Univ.
\yr 1997
\vol RIMS-1129
\moreref alg-geom/9612002
\endref

\ref
\key HM1
\by J. Harvey, G. Moore
\paper Algebras, BPS-states, and strings
\jour Nucl. Physics.
\vol B463
\yr 1996
\pages 315
\moreref hep-th/9510182
\endref

\ref
\key HM2
\by J. Harvey, G. Moore
\paper On the algebras of BPS-states
\jour Preprint
\yr 1996
\moreref hep-th/9609017
\endref

\ref
\key Kac
\by V. Kac
\book Infinite dimensional Lie algebras
\yr 1990
\publ Cambridge Univ. Press
\endref

\ref
\key Kaw1
\by T. Kawai
\paper String duality and modular forms
\jour Preprint
\yr 1996
\moreref hep-th/9607078
\endref


\ref
\key Kaw2
\by T. Kawai
\paper K3 surfaces, Igusa cusp forms and string theory
\jour Preprint
\yr 1997
\moreref hep-th/ 9710016
\endref

\ref
\key Kn
\by M. Kneser
\paper Klassenzahlen indefiniter quadratischer Formen in drei
oder mehr Ver\"ander\-lichen
\jour Arch. Math. (Basel)
\yr 1956
\vol 7 \pages 323--332
\endref

\ref
\key M
\by G. Moore
\paper String duality, automorphic forms, and generalized
Kac--Moody algebras
\jour Preprint \yr 1997
\moreref hep-th/9710198
\endref

\ref
\key N1
\by V.V. Nikulin
\paper Integral symmetric bilinear forms and some of
their geometric applications
\jour Izv. Akad. Nauk SSSR Ser. Mat.
\vol  43
\yr 1979
\pages 111--177
\transl\nofrills English transl. in
\jour Math. USSR Izv.
\vol 14
\yr 1980
\endref

\ref
\key N2
\by V.V. Nikulin
\paper On factor groups of the automorphism groups of
hyperbolic forms modulo subgroups generated by 2-reflections
\jour Dokl. Akad. Nauk SSSR
\yr 1979
\vol 248
\pages 1307--1309
\transl\nofrills English transl. in
\jour Soviet Math. Dokl.
\yr 1979
\vol 20
\pages 1156--1158
\endref

\ref
\key N3
\by V.V. Nikulin
\paper On the quotient groups of the automorphism groups of
hyperbolic forms by the subgroups generated by 2-reflections,
Algebraic-geometric applications
\jour Current Problems in Math. Vsesoyuz. Inst. Nauchn. i
Tekhn. Informatsii, Moscow
\yr 1981 \vol 18
\pages 3--114
\transl\nofrills English transl. in
\jour J. Soviet Math.
\yr 1983
\vol 22
\pages 1401--1476
\endref

\ref
\key N4
\by V.V. Nikulin
\paper On arithmetic groups generated by
reflections in Lobachevsky spaces
\jour Izv. Akad. Nauk SSSR Ser. Mat.
\vol  44   \yr 1980 \pages 637--669
\transl\nofrills English transl. in
\jour Math. USSR Izv.
\vol 16 \yr 1981
\endref

\ref
\key N5
\by V.V. Nikulin
\paper On the classification of arithmetic groups generated by
reflections in Lobachevsky spaces
\jour Izv. Akad. Nauk SSSR Ser. Mat.
\vol  45
\issue 1
\yr 1981
\pages 113--142
\transl\nofrills English transl. in
\jour Math. USSR Izv.
\vol 18
\yr 1982
\endref

\ref
\key N6
\by V.V. Nikulin
\paper
Surfaces of type K3 with finite automorphism group and Picard group of
rank three
\jour Trudy Inst. Steklov
\yr 1984
\vol 165
\pages 113--142
\transl\nofrills English transl. in
\jour  Proc. Steklov Math. Inst.
\yr 1985
\vol 3
\endref

\ref
\key N7
\by V.V. Nikulin
\paper On a description of the automorphism groups of
Enriques surfaces
\jour Dokl. AN SSSR \vol 277 \yr 1984 \pages 1324--1327
\transl\nofrills English transl. in
\jour  Soviet Math. Dokl.
\yr 1984
\vol 30 \pages 282--285
\endref

\ref
\key N8
\by V.V. Nikulin
\paper Discrete reflection groups in Lobachevsky spaces and
algebraic surfaces
\inbook Proc. Int. Congr. Math. Berkeley 1986
\vol  1
\pages 654--669
\endref

\ref
\key N9
\by V.V. Nikulin
\paper Basis of the diagram method for generalized reflection groups
in Lobachev\-sky spaces and algebraic surfaces with nef anticanonical
class
\jour Intern. J. of Mathem.
\vol  7 \yr 1996  \issue 1
\pages 71--108
\moreref alg-geom/9405011
\endref

\ref
\key N10
\by V.V. Nikulin
\paper A lecture on Kac--Moody Lie algebras of the arithmetic type
\jour Preprint Queen's University, Canada
\vol \#1994-16,
\yr 1994 \moreref alg-geom/9412003
\endref


\ref
\key N11
\by V.V. Nikulin
\paper Reflection groups in Lobachevsky spaces and
the denominator identity for Lorent\-zian Kac--Moody algebras
\jour Izv. Akad. Nauk of Russia. Ser. Mat.
\vol  60
\issue 2
\yr 1996
\pages 73--106
\transl\nofrills English transl. in
\jour Izvestiya Math. \vol 60 \yr 1996 \issue 2
\pages 305--334
\moreref alg-geom/9503003
\endref

\ref
\key N12
\by V.V. Nikulin
\paper The remark on discriminants of K3 surfaces moduli as sets
of zeros of automorphic forms
\jour  J. of Mathematical Sciences, \vol 81 \issue 3
\yr  1996 \pages 2738--2743
\publ Plenum Publishing
\moreref alg-geom/9512018
\endref

\ref
\key N13
\by V.V. Nikulin
\paper K3 surfaces with interesting groups of automorphisms
\jour  Preprint RIMS Kyoto University
\yr 1997 \vol RIMS-1132
\moreref alg-geom/ 9701011
\endref


\ref
\key N14
\by V.V. Nikulin
\paper On the classification of hyperbolic root systems of
the rank three. Part I
\jour  Duke e-prints
\yr 1997
\moreref alg-geom/ 9711032
\endref


\ref
\key P-\u S\u S
\by I.I. Pjatetcki\u i-\u Sapiro, \ I.R. \u Safarevich
\paper A Torelli theorem for algebraic surfaces of type K3
\jour Izv. Akad. Nauk SSSR Ser. Mat.
\vol  35  \yr 1971 \pages 530--572
\transl\nofrills English transl. in
\jour Math. USSR Izv.
\vol 5 \yr 1971
\endref

\ref
\key R
\by M.S. Raghunatan
\book Discrete subgroups of Lie groups
\publ Springer-Verlag
\yr 1972
\endref


\ref
\key SW
\by R. Scharlau and C. Walhorn
\paper Integral lattices and hyperbolic reflection groups
\jour Ast\'erisque
\yr 1992 \vol 209 \pages 279--291
\endref

\ref
\key Se
\by J.-P. Serre
\book Cours d'arithm\'etique
\yr 1970
\publ Presses Universitaires de France
\publaddr Paris
\endref

\ref
\key Sh
\by D. Shanks
\paper Class number, a theory of factorization, and genera
\inbook in Proc. Sympos. Pure Math.
\vol 20
\publ Amer. Math. Soc.
\publaddr Providence, R. I.
\yr 1970
\pages 415--440
\endref

\ref
\key V1
\by \'E.B. Vinberg
\paper Discrete groups generated by reflections in Lobachevsky
spaces
\jour Mat. Sb. (N.S.)
\vol 72 \yr  1967 \pages 471--488
\transl\nofrills English transl. in
\jour Math. USSR Sb. \vol 1 \yr 1967 \pages 429--444
\endref

\ref
\key V2
\by \'E.B. Vinberg
\paper On groups of unit elements of certain quadratic forms
\jour Mat. Sbornik
\yr 1972
\vol 87
\pages 18--36
\transl\nofrills English transl. in
\jour Math USSR Sbornik
\vol 16
\yr 1972
\pages 17--35
\endref

\ref
\key V3
\by \'E.B. Vinberg
\paper The absence of crystallographic reflection groups in
Lobachevsky spaces of large dimension
\jour Trudy Moscow. Mat. Obshch. \yr 1984 \vol 47 \pages 67--102
\transl\nofrills English transl. in
\jour Trans. Moscow Math. Soc. \yr 1985 \vol 47
\endref

\ref
\key V4
\by \'E.B. Vinberg
\paper Hyperbolic reflection groups
\jour Uspekhi Mat. Nauk
\vol 40
\yr 1985
\pages 29--66
\transl\nofrills English transl. in
\jour Russian Math. Surveys
\vol 40
\yr 1985
\endref

\ref
\key VSh
\by \'E.B. Vinberg and O.V. Shvartsman
\paper Discrete groups of motions of spaces of constant curvature
\inbook Sovrem. problemy matem.
Fundam. Napr. Vol. 29, Geometriya 2
\publ VINITI, Moscow
\yr 1988
\pages 147--259
\transl\nofrills English transl. in
\inbook Encyclopaedia of Math. Sciences. Geometry II
\vol 29
\publ Springer-Verlag
\yr 1991
\endref

\ref
\key W
\by C. Walhorn
\paper Arithmetische Spiegelungsgruppen auf dem
4-dimensionalen hyperbolischen Raum
\jour Dissertation zur Erlangung des Doktorgrades der
Fakult\"at f\"ur Mathematik der Universit\"at Bielefeld
\yr 1993
\endref

\endRefs
\enddocument

\end